\documentclass[final,1p,times]{elsarticle}
\usepackage{amssymb}
\usepackage{amsmath, amssymb}
\usepackage{geometry}
\usepackage{mathtools}
\usepackage[T1]{fontenc}
\usepackage{etoolbox}
\usepackage[normalem]{ulem}

\usepackage{url}
\usepackage{braket,amsfonts}
\usepackage[percent]{overpic}
\usepackage{array}
\usepackage{booktabs}
\usepackage{multirow}
\usepackage[caption=false]{subfig}
\usepackage{algorithm}
\usepackage{algpseudocode}

\usepackage{graphicx,epstopdf}

\usepackage{amsopn}

\usepackage{xspace}
\usepackage{bold-extra}
\usepackage[most]{tcolorbox}

\usepackage{color}
\newcommand{\red}[1]{{\color{red} #1}}
\newcommand{\blue}[1]{{\color{blue} #1}}

\newcommand{\yellow}[1]{{\color{yellow} #1}}

 \definecolor{mygray}{gray}{0.2}

\newcommand{\BS}[1]{{\boldsymbol{#1}}}

\def\R{{\mathbb{R}}}

\def\A{{{A}}}

\def\L{{\mathcal{L}}}

\def\z{{\BS{z}}}

\def\nn{{\BS{n}}}

\def\xx{{\BS{\xi}}}
\def\qed{{$\blacksquare$}}
\def\Span{\mbox{span}}

\definecolor{frenchblue}{rgb}{0.0,0.45,0.73}
\definecolor{sinopia}{rgb}{0.8,0.25,0.04}

\newcounter{example}[section]
\renewcommand{\theexample}{\thesection.\arabic{example}} 
\newenvironment{example}[1][] 
{
    \refstepcounter{example} 
    \par\medskip 
    \noindent \textbf{Example~\theexample#1} \rmfamily 
    \rmfamily
}
{
    \medskip 
}

\newcounter{observation}[section]
\renewcommand{\theobservation}{\arabic{observation}} 
\newenvironment{observation}[1][] 
{
    \refstepcounter{observation} 
    \par\medskip 
    \noindent \textit{Observation~\theobservation : #1}   
    \rmfamily
}
{
\medskip 
}

\colorlet{texcscolor}{blue!50!black}
\colorlet{texemcolor}{red!70!black}
\colorlet{texpreamble}{red!70!black}
\colorlet{codebackground}{black!25!white!25}

\usepackage{rotating}
\usepackage{makecell, tabularx}

\usepackage{boondox-cal} 
\def\m{{\mathbcal{m}}}
\def\n{{\mathbcal{n}}}
\def\greedyk{{\mathcal{k}}}
\DeclareMathAlphabet{\mathcalboondox}{OMS}{bodc}{m}{n}
\newcommand{\calo}{\mathcalboondox{O}}
\newcommand{\calk}{\mathcalboondox{K}}
\newcommand{\cals}{\mathcalboondox{S}}
\newcommand{\calu}{\mathcalboondox{U}} 

\def\u{{U}}

\def\onetom{\mathbb{N}_m}
\def\oneton{\mathbb{N}_n}
\def\dt{\triangle t}
\def\mecheps{\varepsilon_{mech}}

\DeclareTotalTCBox{\code}{ v O{} }
{ 
  fontupper=\ttfamily\color{black},
  nobeforeafter,
  tcbox raise base,
  colback=codebackground,colframe=white,
  top=0pt,bottom=0pt,left=0mm,right=0mm,
  leftrule=0pt,rightrule=0pt,toprule=0mm,bottomrule=0mm,
  boxsep=0.5mm,
  #2}{#1}

\patchcmd\newpage{\vfil}{}{}{}
\journal{Mathematics and Computers in Simulation}

\begin{document}

\begin{frontmatter}

\title{Greedy Trial Subspace Selection in Meshfree Time-Stepping Scheme with Applications in Coupled Bulk-Surface Pattern Formations}

\author[1]{Yichen Su\corref{cor1}}
\ead{21481210@life.hkbu.edu.hk}

\author[1]{Leevan Ling}
\ead{lling@hkbu.edu.hk}

\affiliation[1]{organization={Department of Mathematics},
            addressline={Hong Kong Baptist University},
            country={Hong Kong}}
\cortext[cor1]{Corresponding author}

\begin{abstract}
Combining kernel-based collocation methods with time-stepping methods to solve parabolic partial differential equations can potentially introduce challenges in balancing temporal and spatial discretization errors. Typically, using kernels with high orders of smoothness on some sufficiently dense set of trial centers provides high spatial approximation accuracy that can exceed the accuracy of finite difference methods in time. The paper proposes a greedy approach for selecting trial subspaces in the kernel-based collocation method applied to time-stepping to balance errors in both well-conditioned and ill-conditioned scenarios. The approach involves selecting trial centers using a fast block-greedy algorithm with new stopping criteria that aim to balance temporal and spatial errors. Numerical simulations of coupled bulk-surface pattern formations, a system involving two functions in the domain and two on the boundary, illustrate the effectiveness of the proposed method in reducing trial space dimensions while maintaining accuracy.

\end{abstract}

\begin{graphicalabstract}
\end{graphicalabstract}

\begin{highlights}
\item Research highlight 1
\item Research highlight 2
\end{highlights}


\end{frontmatter}


\section{Introduction}
\label{sec:intro}

As kernel-based methods gain popularity, finding appropriate settings remains an open problem: how to stably solve the severely ill-conditioned linear systems that arise. Various strategies have been employed to address instability issues in computations. Regularization, as demonstrated in \cite{AMIRFAKHRIAN2016278, Vito2004SomePO}, is effective in mitigating ill-conditioning problems. In this paper, we consider selecting a trial subspace to achieve more stable solution approximations. The Greedy algorithm was initially explored for symmetric kernel-based interpolation matrix systems by Schaback et al \cite{Schaback+Wendland-Adapgreetechappr:00} and remain an active topic. Traditional methods for adaptive trial subspace selection include the $P$-greedy algorithm, which applies to symmetric PDEs. For the latest developments on power function-based $P$-greedy, see \cite{Campagna+DeMarchi-AdvComputMath:22}. Additionally, \cite{Wenzel-AnalTargDataGree:23} proposes the residual-based $f$-greedy algorithm and combinations of the $P$- and $f$-greedy algorithms to achieve higher convergence rates. The greedy algorithm has advantages for solving both time-dependent and time-independent problems. In the context of time-independent problems, it is directly employed to choose trial subspaces for approximating solutions. Regarding time-dependent problems, we first apply time-stepping methods and subsequently apply the selected trial subspaces for solution approximation at each time step. Following the pioneering work on greedy methods \cite{Hon+SchabackETAL-adapgreealgosolv:03} for solving time-independent partial differential equations (PDEs) using collocation methods, a.k.a. Kansa methods, we proposed various sequential-greedy algorithms \cite{Ling+OpferETAL-Resumeshcolltech:06,Ling+Schaback-imprsubsselealgo:09} to select quasi-optimal sets of trial subspaces that guarantee stable solutions. They built up nonsingular subsystems iteratively by selecting a row and then a column in a greedy fashion. 
Although each version brings improved accuracy, they suffer from high complexity. For instance, in a well-conditioned system, a greedy algorithm might select subsystems that yield unexpectedly high accuracy. For an $n\times n$ asymmetric collocation system, selecting $\greedyk$ columns by any of these algorithms costs $\calo(\greedyk^4+n\greedyk^2)$. The selection of multiple rows and columns results in large computations, clearly reflecting the trade-off between achieving higher accuracy and managing algorithm complexity. The block-greedy algorithm addresses the constraints of sequential algorithms \cite{Ling-fastblocalgoquas:16} and can be implemented in a matrix-free manner, offering a reduced complexity of $\calo(n\greedyk^2)$.
It can also be extended to a fully adaptive algorithm for boundary value problems \cite{Ling+Chiu-Fulladapkernmeth:18} of elliptic type. It allows new data points to be added gradually, dynamically refining the trial subspace without relying on fixed data distribution and enhances the efficiency by solving least-squares problems. \\
\indent A strategy ensures that fill distances of boundary/domain points are similar. Block-greedy to address ill-conditioning from small minimum separating distances.
Specifically, during the subspace selection iterations of the block greedy algorithm, each iteration adds a set of new data points according to the primal/dual residual criterion outlined in \cite{Ling+Schaback-imprsubsselealgo:09}. This criterion helps avoid point clustering and maintains uniform fill distances among boundary and domain points, as detailed in \cite{Ling-fastblocalgoquas:16} Theorem 2.1. As the selection proceeds, the algorithm monitors the condition number of the selected subsystem, and the greedy algorithm also checks whether the expanded submatrix is well-conditioned to allow further expansion. The stability of these kernel-based methods also heavily depends on the kernel's shape parameter. In this paper, a data-driven shape parameter generator fine-tunes this parameter through extensive testing to optimize system stability and numerical accuracy, enhancing the overall performance and reliability. \\
\indent Sobolev kernels work well numerically, producing reasonable approximations with refinement for increased accuracy. This fully automatic meshfree algorithm benefits Kansa method users. While effective for elliptic PDEs, experimental evidence shows that using the greedy algorithm with time-stepping methods for solving parabolic PDEs is not as numerically stable as desired. We provide evidence with unstable numerical experiments in Example \ref{eg: 2D RD PDE}. One significant reason for the instability observed here is the mismatch in the accuracy of spatial and temporal discretizations of PDEs. To address this, we modified the original greedy algorithm by introducing different tolerance values for residuals and condition numbers. This approach efficiently handles issues arising from the accuracy mismatch problem.
In Section~\ref{sec greedy review}, we overview the block-greedy algorithm.
In Section~\ref{sec select tol}, we explore how to properly fine-tune the algorithm's stopping criteria to fix the problem. In this paper, we first discretize PDEs in time and then in space, as a fixed time step is essential to ensure balanced convergence between time and spatial dimensions. We thoroughly studied the least-squares, kernel-based collocation method of lines in \cite{chen2023exploring, Chen_2023}. These approaches have both theoretically and numerically demonstrated the importance of oversampling, which can also be achieved by using the block-greedy algorithm. But this will not be considered in this paper. Numerical experiments verify the effectiveness of the proposed algorithm, followed by conclusions.

\section{Block-greedy algorithm for column subspace selections} \label{sec greedy review}
From here on, \textit{greedy algorithm} refers to the block-version in \cite{Ling-fastblocalgoquas:16}. Although designed for kernel-based collocation methods, it can be viewed purely in terms of linear algebra.
The following overview is not meant to be exhaustive but rather provides sufficient details for this paper.\\
Consider the linear system $\A  \boldsymbol{\lambda} = \mathbf{b}$, where $\A$ is an $m \times n$ matrix with full rank $m$. The greedy algorithm is applied to the linear system and selects a column space of $\A$ that allows accurate and stable computations, ensuring that provides an effective approximation of the vector $\mathbf{b}$.

Suppose at a non-initial step, rows indexed by $\m\subset \onetom :=\{1,2,\ldots,m\}$ and columns indexed by $\n\subset  \oneton$ have been selected. Let $\A(\m,\n)$ denote the corresponding submatrix. The goal is to expand the selections to $\m'$ and $\n'$ in batches rather than one by one, such that:
\begin{itemize}
  \item $|\m'| \gtrsim 2|\m|$  while $|\m'|\leq m$, and
  \item $|\n'|  = 2|\n|$ while $|\n'|< n$,
\end{itemize}
to build an overdetermined subsystem as it iterates.
Here $|\cdot|$ represents the size of the row/ column set, and $|\m'| \gtrsim 2|\m|$ indicates that the number of rows in $\m'$ is greater than or approximately equal to two times of the number of rows in $\m$. When running out of unselected rows or columns, the greedy algorithm will take all $m$ rows or $n$ columns.

    Instead of searching all unselected rows in $\onetom\setminus \m$ and columns in $\oneton\setminus \n$, the greedy algorithm shortlists candidates by sorting the residual magnitudes. To obtain the primal and dual residuals, we formulate a constrained minimization problem. The objective function aims to minimize the Euclidean norm of the solution vector $\boldsymbol{\eta}$, subject to the linear constraint $\A\boldsymbol{\eta} = \mathbf{b}$. To solve this constrained optimization problem, we employ the method of Lagrange multipliers. The optimisation problem is transferred in the form of 

\begin{equation}
    \label{eq: lagrange multi}
    \L(\boldsymbol{\eta}, \boldsymbol{\zeta}) = \frac{1}{2}\boldsymbol{\eta}^{T}\boldsymbol{\eta} + \boldsymbol{\zeta}(\A\boldsymbol{\eta} -  \mathbf{b}),
\end{equation}
where $\boldsymbol{\zeta}$ denotes the Lagrange  multiplier. We take the partial derivatives of the Lagrangian function $\L(\boldsymbol{\eta}, \boldsymbol{\zeta})$ with respect to $\boldsymbol{\eta}$ and $\boldsymbol{\zeta}$ to derive the primal-dual subsystem. The primal-dual subsystem is defined in the form of
      \begin{eqnarray}
          \A(\m,\n) \boldsymbol{\eta} &=& \mathbf{b}(\m),  \label{eq p res}
          \\
          \A(\m,\n)^T \boldsymbol{\zeta} &=& -\boldsymbol{\eta}.\label{eq d res}
      \end{eqnarray}
 Then the primal and dual residuals are defined by
    \begin{equation}\label{eq p/d res}
    \mathbf{r}_{\text{primal}}:= \A(\onetom,\n)\boldsymbol{\eta}-\mathbf{b} \in \R^m,
    \text{\quad and \quad}
    \mathbf{r}_{\text{dual}}:=  \mathcal{E}_{\R^{|\n|}\to\R^n}\boldsymbol{\eta}+\A(\m,\oneton)^T\boldsymbol{\zeta}\in \R^n,
\end{equation}
where $\mathcal{E}_{\R^{|\n|}\to\R^n}\boldsymbol{\eta} \in \R^n$ is the extension of $\boldsymbol{\eta}\in\R^{|\n|}$ with zeros patched into entries not in $\n$.
Then, the greedy algorithm identifies and selects candidate rows and columns from those that remain unselected, choosing those associated with the largest magnitudes of $\mathbf{r}_{\text{primal}}$ and $\mathbf{r}_{\text{dual}}$ in the primal and dual systems, respectively. The final selection is made from shortlisted candidates based on these smaller submatrices.

The greedy algorithm is matrix-free because, as seen in \eqref{eq p/d res}, it only needs to store entries of selected rows and columns, i.e., $A(\m,\oneton)$ and $A(\onetom,\n)$. The computational cost of the greedy algorithm is determined by the size of the shortlisted candidate sets in each iteration.
We refer readers to \cite{Ling-fastblocalgoquas:16} for details on the complexity analysis. As mentioned in the introduction, the greedy algorithm achieves $\calo(n^2 \greedyk)$ complexity, where $\greedyk$ denotes the number of total selected columns.

The greedy algorithm can be initialized with pre-selected rows/columns, i.e., input some $\m$ and $\n$. If none is specified, it begins with rows/columns associated with the largest primal/dual residuals.

After each expansion that (roughly) doubles the number of selected rows/columns, the stopping criteria of the greedy algorithm use two tolerance values, $\tau_{r}$ and $\tau_{\kappa}$.
Besides stopping by a small enough residual, the greedy algorithm will check the condition number of $\A(\m', \n')$ to determine whether it is well-conditioned for further expansion; otherwise, the expansion terminates because of the ill-conditioning. In particular, the greedy algorithm stops when
\begin{itemize}
  \item \textbf{[SC-$1$]} the condition number of the expanded system $\kappa( \A(\m',\n' ) ) > 1/\tau_{\kappa}$ exceeds the tolerance, or
  \item \textbf{[SC-$2$]} the norm of the primal residual $\|  \mathbf{r}_{\text{primal}} \|_\infty < \tau_{r}$ is below the tolerance.
\end{itemize}
It is cheap to compute the residual and check the condition number since the algorithm updates  \cite{Hammarling+Lucas-UpdaFactLeasSqua:08} QR factorizations of $\A(\m,\n)$ to $\A(\m',\n')$ in each iteration.
Recall that in the block-style greedy algorithm, the number of selected columns doubles at each iteration.
When we stop due to SC-$1$, 
a backtracking process is employed to prevent including too many columns, which would lead to a condition number much larger than what the user tolerance allows.
This backtracking process uses the bisection method to select the first $\greedyk$ indices in $\n'$ to generate the largest subset $\n'' \subset \n'$ so that the final column subspace selection yields a condition number just below tolerance, i.e., $\kappa(\A(\m', \n'')) \leqslant 1/\tau_{\kappa}$.
In the case of SC-$2$, the greedy algorithm returns $\n'$ as its final selection. In the original article, a single tolerance $\tau_{r}=\tau_{\kappa} = \varepsilon_{mech}$, the machine epsilon, was used as the default value, which works well in various time-independent elliptic PDEs. To clearly demonstrate the column subspace selection process in the greedy algorithm, the brief framework is outlined in Algorithm \ref{alg:col sel algo}. 

\begin{algorithm}
\caption{Column subspace selection algorithm}\label{alg:col sel algo}
\begin{algorithmic}[1]
\State \textbf{Inputs:} $\A(\onetom, \oneton) \in \R^{m\times n}$,  $\mathbf{b} \in \R^{m}$ $\tau_{\kappa}$, $\tau_{r}$
\State \textbf{Initialize:} 
\State$\m \leftarrow $ row index with the maximum primal residual  
\State$\n \leftarrow $ column index with the maximum dual residual
\While{$|\n| < n$}
    \State Solve primal-dual subsystems (\ref{eq p res}) and (\ref{eq d res}) by QR factorization of $\A(\m, \n)$ 
    \State Calculate $\mathbf{r}_{\text{primal}}$ and $\mathbf{r}_{\text{dual}}$ 
    \State Break if $\mathbf{r}_{\text{primal}} < \tau_{r}$     \hfill   $\triangleright$ \textbf{[SC-$2$]}
    \State Select rows $\m'$ and columns $\n'$ in batches, $\m \leftarrow \m'$, $\n \leftarrow \n'$
    \State Estimate the condition number $\kappa( \A(\m',\n' ) )$  
    \If{$\kappa(\A(\m, \n)) > \tau_{\kappa}$}       \hfill   $\triangleright$ \textbf{[SC-$1$]}
    \State Bisection method to select the first $\greedyk$ indices in $\n'$ to generate $\n''$
    \State $\n \leftarrow \n''$
    \State Break
    \EndIf
    
\EndWhile
\end{algorithmic}
\end{algorithm}

\section{Controlling spatial accuracy in time-stepping methods}
\label{sec select tol}
Before applying the block-greedy algorithm to solve the resultant matrix systems arising from kernel-based collocation and time-stepping methods, we first walk through the main steps of this discretization procedure.
Consider second-order parabolic partial differential equations defined in a bounded domain $\Omega \subset \mathbb{R}^d$ with smooth boundary $\partial\Omega$. The equations are complemented by initial and boundary conditions that we assume are compatible.
 
As a demonstrative example, the heat equation with Dirichlet boundary conditions is defined in the form: 
\begin{equation}\label{eq: GF RDE}
\begin{array}{r@{\;=\;}llrcl}
    \displaystyle
    \frac{\partial u}{\partial t} &  D\nabla^{2}u + f
        & \text{in}& \Omega &\times&[0, T] \\
    u &  g
        &\text{on} & \partial \Omega  &\times&  [0, T]\\
    u &  u_{0}
        &\text{on} &  \Omega  &\times&  \{0\}\\
\end{array}
\end{equation}
with constant diffusion $D>0$ and source function $f$. While this simple example illustrates the approach, the proposed method applies more generally to quasi-linear parabolic PDEs of the above form with variable coefficients.

We define a partition $\{t_k\}_{k=0}^K$ of [0,T] with time steps $\dt _{k} = t_{k} - t_{k-1}$. We denote the approximate solution at time $t_k$ by $\u^{k}\approx u(\boldsymbol{\cdot}, t_k)$ for $k = 0,1,\ldots,K$.
The commonly used Crank-Nicolson scheme semi-discretizes the PDE \eqref{eq: GF RDE} in time and results in a sequence of elliptic PDEs of the form
\begin{equation}
\label{eq: CN}
\frac{\u^{k} - \u^{k-1}}{\Delta t_{k}} = \frac{1}{2}\big(D\nabla^{2}\u^{k} + f^{k} + D\nabla^{2}\u^{k-1} + f^{k-1}\big)
\quad\text{for $k=1,\dots,K$},
\end{equation}
where $f^k := f(\boldsymbol{\cdot}, t_k)$.
This paper also employs semi-implicit backward differentiation formulas (SBDF) in \cite{Ruuth-Implmethreacprob:95}, which use explicit schemes for reaction terms and implicit schemes for diffusion terms.
The system generated by the first-order SBDF (SBDF1) is
\begin{equation}
\label{eq: sbdf1}
\frac{\u^{k} - \u^{k-1}}{\dt _{k}} = f^{k-1} + D\nabla^{2}\u^{k},
\quad\text{for $k=1,\dots,K$},
\end{equation}
and that generated by the second-order SBDF (SBDF2)  is
\begin{equation}
\label{eq: sbdf2}
\frac{1}{2\dt _{k}}\Big(3\u^{k} -4\u^{k-1} + \u^{k-2}\Big) = 2f^{k-1} - f^{k-2} + D\nabla^{2}\u^{k},
\quad\text{for $k=2,\dots,K$},
\end{equation}
with $\u^1$ obtained by \eqref{eq: sbdf1}.

We solve the semi-discretized systems \eqref{eq: CN}, \eqref{eq: sbdf1} or \eqref{eq: sbdf2} using a kernel-based collocation method.
We choose a radius basis kernel $\Phi$ on $\mathbb{R}^d$ giving rise to a positive definite kernel with a smoothness order greater than 2. Common examples include the Gaussian and multiquadrics kernels. Reproducing kernels of Sobolev spaces, like the Whittle-Matern-Sobolev (MS) kernel
and compactly supported Wendland functions, are also suitable choices that come with convergence theory for elliptic PDEs \cite{Chen+Ling-Extrmeshcollmeth:20,Cheung+Ling-Kernembemethconv:18,Cheung+LingETAL-leaskerncollmeth:18}. The theoretical approximation powers of employing these kernels on greedy points were addressed in \cite{Santin_2021}.

Given a set  $\Xi:=\{\xx^{j}\}_{j=1}^{n} \subset \Omega$, we define the finite-dimensional trial space
\begin{equation}
    \label{eq: trial space}
    \calu:=\Span{ \Big\{\Phi(\boldsymbol{\cdot} - \xx )\, \Big|\, \xx  \in
    \Xi 
    \Big\}}.
\end{equation}
We seek approximation to the solution $\u^k$ within $\calu$ in the form of 

\begin{equation}
    \label{eq: appro u}
\sum_{j = 1}^{n} \lambda_{j}^{k} \Phi(\boldsymbol{\cdot}- \xx^{j}),
\qquad\text{for $k=1,\dots,K$},
\end{equation}
where $\boldsymbol{\lambda}^k = \{\lambda_j^k\}_{j=1}^n$ are the unknown coefficients.
We determine $\boldsymbol{\lambda}^k$ by imposing collocation at points $Z= Z_\Omega \cup Z_{\partial \Omega} = \{\mathbf{z}^i\}_{i=1}^m \subset \Omega \cup \partial\Omega$ for the governing equation and boundary conditions respectively. Taking the CN scheme \eqref{eq: CN}, we obtain an $m \times n$ matrix system:
\begin{equation}
    \label{eq: alam = b}
\begin{pmatrix}
        ( 2 \Phi - \dt _k [\nabla^2\Phi])(Z\cap{\Omega}, \Xi) \\
        \Phi(Z\cap{\partial\Omega}, \Xi)
    \end{pmatrix}
    \\ =
    \begin{pmatrix}
        \big( 2\u^{k-1} + \dt (f^{k} + D\nabla^{2}\u^{k-1} + f^{k-1}) \big) (Z\cap{\Omega})\\ g^k(Z\cap{\partial\Omega})
    \end{pmatrix},
\end{equation}
where $g^k := g(\boldsymbol{\cdot}, t_k)$ and the data dependent matrix with entries
\[
   [ \Phi(Z_\Omega, \Xi) ]_{i,j} = \Phi( \z_i - \xx_j)
   \qquad\text{for $\z_i\in Z_\Omega$ and $\xx_i\in\Xi$}.
\]

The other two matrices $\Phi(Z_{\partial\Omega}, \Xi)$ and $[\nabla^2\Phi](Z_{\Omega}, \Xi)$ are defined similarly. The right-hand side vector is known; given $\boldsymbol{\lambda}^{k-1}$, we can evaluate the approximate solution $\u^{k-1}$ and its Laplacian $\nabla^2\u^{k-1}$ using \eqref{eq: appro u}. In linear algebra notation, we can express the fully-discretized  systems \eqref{eq: alam = b} as
\begin{equation}\label{eq Ax=b}
  \A \boldsymbol{\lambda}^k = \mathbf{b}(\boldsymbol{\lambda}^{k-1}),
\quad\text{for $k=1,\dots,K$},
\end{equation}
with $\boldsymbol{\lambda}^0$ or $\mathbf{b}(\boldsymbol{\lambda}^0)$  determined from the initial condition. Using a greedy algorithm to select a column subspace of $\A$ is equivalent to selecting a subset of trial centers from $\Xi$. In the following numerical experiments, we implement the greedy algorithm based on the linear system only at the initial time step, which makes sense for parabolic problems with slowly varying solutions. Running the greedy algorithm at a later time is certainly a possibility for solutions that exhibit more variation over time. However, doing so efficiently requires a reliable error indicator, which is beyond the scope of this work.

\begin{example}\label{eg: 2D RD PDE}\textbf{(Greedy Algorithm in action)}
\begin{figure}
\centering
  \subfloat[PDE error with/without the greedy algorithm]{\includegraphics[width=6.5cm]
 {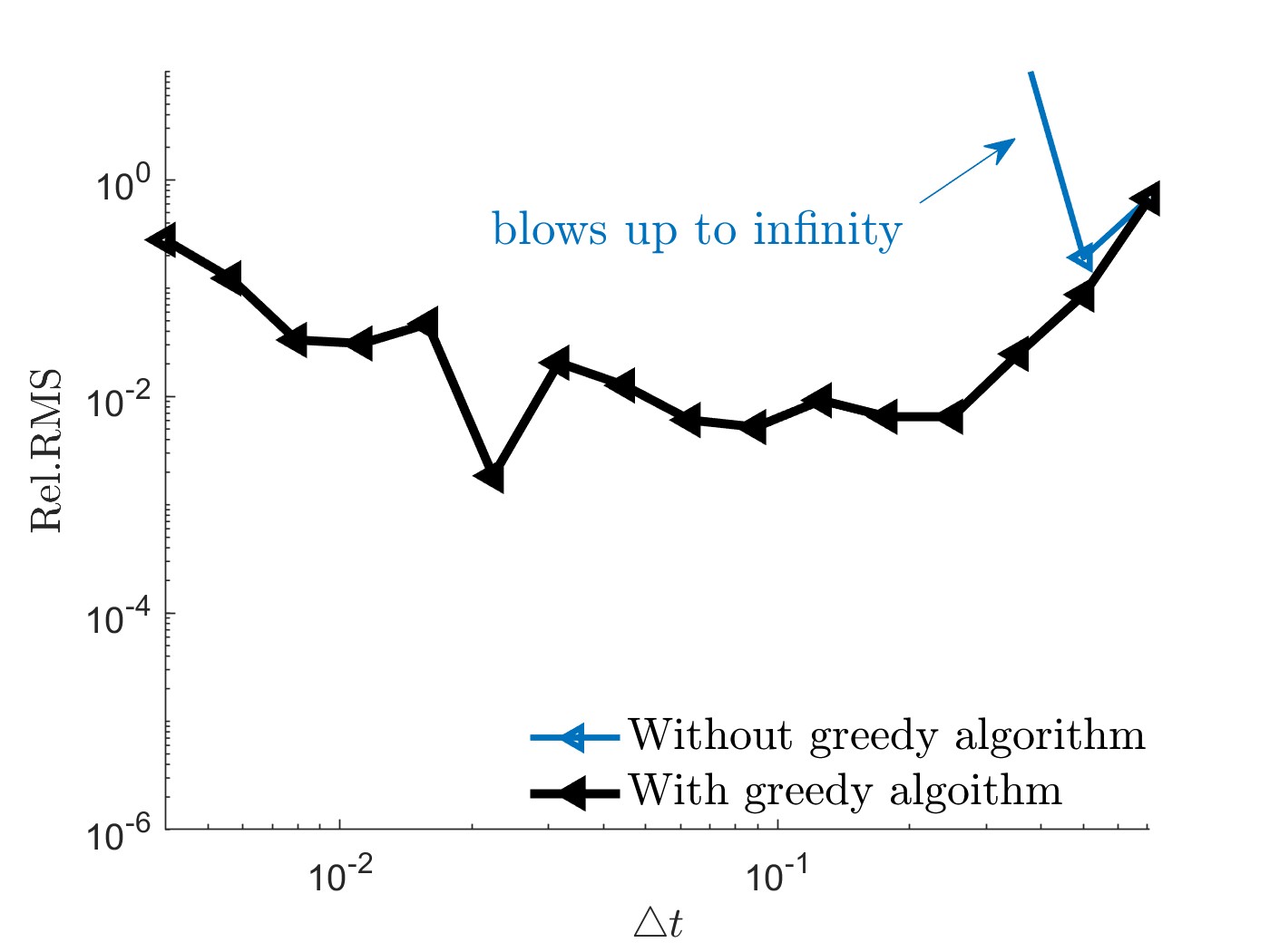}\label{fig:Unstable2Dsqr}}
  \subfloat[Using different tolerances in SC-$2$]{\includegraphics[width=6.5cm]{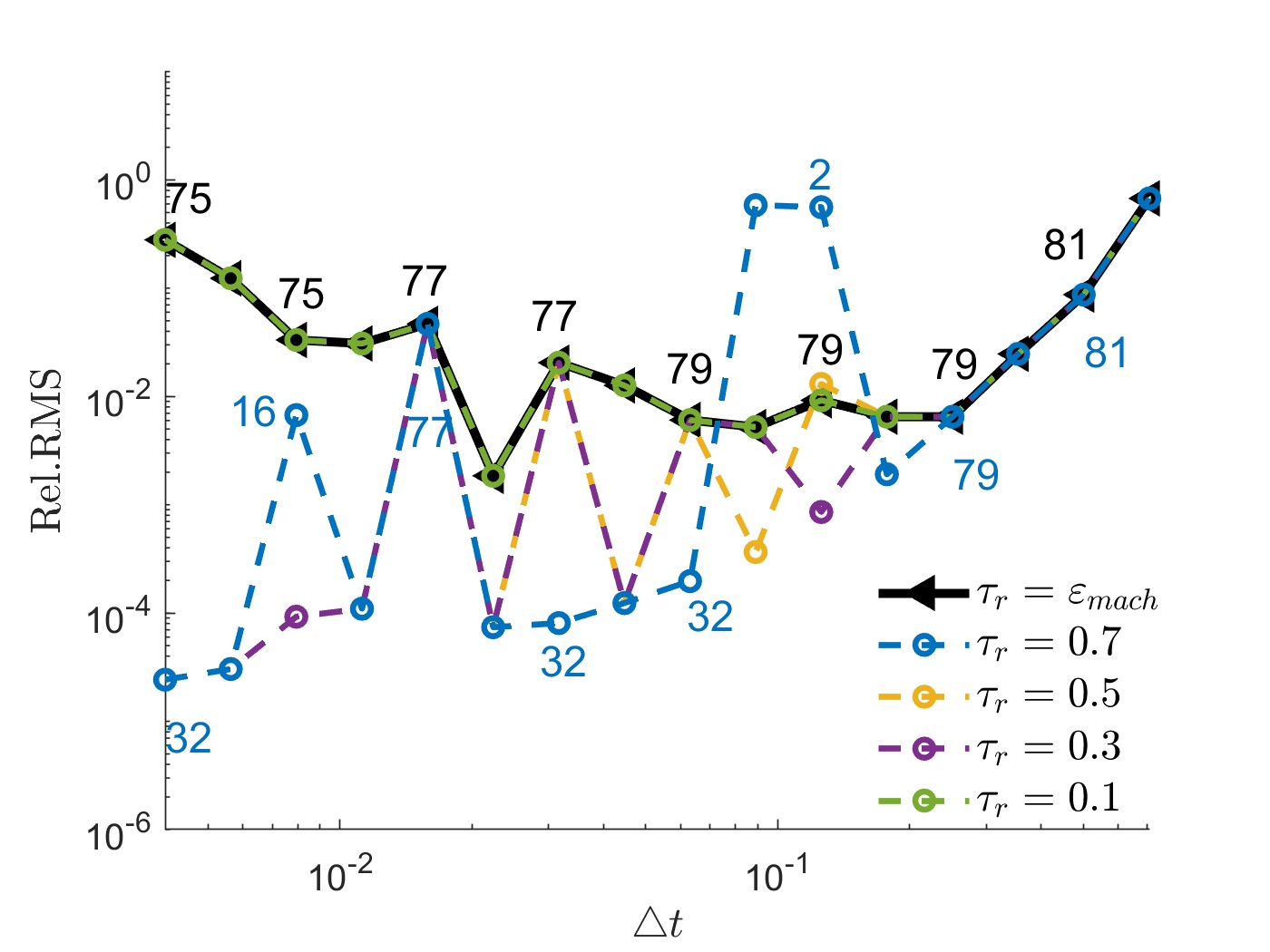}\label{fig:2DComGdyTau}}
   \caption{
   For Example~\ref{eg: 2D RD PDE},
         the relative root mean error profiles of solving a heat equation by meshfree time-stepping method with various $\dt $ (\ref{fig:Unstable2Dsqr}) without and with the greedy algorithm in default settings; (\ref{fig:2DComGdyTau}) with greedy algorithm using different tolerance in stopping criteria. Colored numbers in (\ref{fig:2DComGdyTau}) are the number of selected columns in $\A$ by the greedy algorithm.}
         \label{fig: unstable RDE}
\end{figure}

We compare the performance of the meshfree time-stepping method with and without the greedy algorithm applied. Our comparison is based on a linear heat equation \eqref{eq: GF RDE} with $D = 1$ and $\Omega=[0,1]^2$. The right-hand functions $f$ and $g$ were computed from the exact solution expressed as:
\[
u^{*}([x_1,x_2]^T, t) = \sin(2\pi x_{1}) \sin(\pi x_{2}) \exp{(-2\pi^2t)} + (1 - \exp{(-2\pi^2t)})/2,
\]
which exhibits the standard decay property of heat equations. The Crank-Nicolson scheme was used for temporal discretization, and a fixed time step of $\dt $ was used for all $k$. For spatial discretization, we employed the Halton sequence to generate a set of $n = 300$ data points, which were used as both the trial centers and collocation points. This exactly determined setup will become an overdetermined system if the greedy algorithm selects a proper subset of columns.
The kernel used in this example is the (unscaled) Gaussian kernel, which results in ill-conditioned matrices in the fully discretised system. The relative root mean square error of solutions obtained by the meshfree time-stepping method with various $\dt $ is shown in Figure~\ref{fig: unstable RDE}.

The power of the greedy method is illustrated in Figure \ref{fig:Unstable2Dsqr}. The solution obtained without the greedy algorithm is unstable and diverges for most tested cases. By applying the greedy algorithm (with default parameter values) to the matrix system $\A \boldsymbol{\lambda}^1 = \mathbf{b}(\boldsymbol{\lambda}^{0})$ for the first time step, we obtained a subset of columns of $\A$ indexed by $\n'$.
For all runs, SC-$1$ stopped the algorithm with large condition numbers around $1/\mecheps $.
The resulting submatrix $\A(\mathbb{N}_{m}, \n')$ was then used in the sequence of the fully-discretized system to update solutions in time. This process regularizes the solution but does not fully prevent divergence, and it fails to achieve second-order convergence. We labeled the number of selected columns by the greedy algorithm in Figure~{\ref{fig:2DComGdyTau}} with black numbers and noticed that they do not show significant differences. Regarding the CPU time for the greedy algorithm under different time step settings, note that all tested scenarios concluded after a small number of iterations. Consequently, the CPU times were relatively low\footnote{The CPU time of the greedy algorithm used in each point-selecting procedure is approximately 0.20 seconds.} and provide limited value for comparative analysis. Reported runtime in MATLAB is obtained from an Intel-i7 processor. 

While in the greedy algorithm, we observe that SC-$2$ uses a default tolerance $\tau_r=\mecheps $ for the residual, which is too small to effectively stop the greedy iteration. This forces the selection of extra basis functions, inducing unnecessary ill-conditioning and instability. Therefore, choosing appropriate stopping criteria for the greedy algorithm is important to prevent these issues.
As a test, we ran the greedy algorithm with various tolerance values $0.1\leq \tau_r\leq 0.7$ to produce the error profiles in Figure \ref{fig:2DComGdyTau}. Compared to the default (black solid line), some runs with large $\tau_r$ terminate earlier (as expected) by selecting fewer columns (numbers in the figure). Using fewer columns resulted in 2-3 orders of magnitude accuracy improvement.
However, we noted that these settings' errors lack stability; some runs with only 2 columns are inaccurate.
If $\tau_r=0.1$ is too small to stop the iteration with SC-$2$, then SC-$1$ will stop the algorithm, yielding identical errors as the default. 

These observations indicate the greedy algorithm can stabilize solutions compared to methods without it. However, choosing appropriate stopping criteria is critical to achieving the desired accuracy.
\qed
\end{example}

The kernel-based collocation method using greedy algorithms could efficiently solve semi-discretized PDEs. However, when spatial and temporal discretization accuracies are mismatched, the resulting fully discrete systems may become ill-conditioned, impacting solution accuracy. As demonstrated in Example \ref{eg: 2D RD PDE}, using a smaller $\tau_r$ results in more selected columns to achieve higher spatial accuracy. When the order of spatial discretization significantly surpasses that of temporal discretization, a substantial mismatch in the scales of numerical representation between the two arises. This mismatch can lead to numerical instability.

In the rest of this section,  we develop two stopping strategies for the greedy algorithm to improve condition number and accuracy when solving fully-discrete systems. We propose to:
\begin{itemize}
\item \textbf{[SC-$1'$]} run a residual search within the newly added columns and output the trial sub-subspace that minimizes the least-squares residual, when the greedy algorithm stops due to a large condition number $\kappa( \A(\m',\n' ) ) > 1/\tau_{\kappa}$ with $\tau_\kappa$, and
\item \textbf{[SC-$2'$]}   apply a backtracking process to reduce selected column set accuracy down to a target residual $\tau_{r}' $ when the greedy algorithm stops due to a small residual $\|  \mathbf{r}_{\text{primal}} \|_\infty < \tau_{r} $.
\end{itemize}
These techniques enable the greedy algorithm to balance accuracy, stability, and efficiency based on the conditioning properties of the resulting discrete systems. We implement new stopping criteria by modifying lines 8 to 16 in algorithm \ref{alg:col sel algo}:

\begin{algorithm}
\caption{Column subspace selection algorithm with new stopping criteria}\label{alg:col sel algo SC'}
\begin{algorithmic}[1]
\State \textbf{Inputs:} $\A(\onetom, \oneton) \in \R^{m\times n}$,  $\mathbf{b} \in \R^{m}$ $\tau_{\kappa}$, $\tau_{r}$
\State \textbf{Initialize:} 
\State$\m \leftarrow $ row index with the maximum primal residual  
\State$\n \leftarrow $ column index with the maximum dual residual
\While{$|\n| < n$}
    \State Solve primal-dual subsystems (\ref{eq p res}) and (\ref{eq d res}) by QR factorization of $\A(\m, \n)$ 
    \State Calculate $\mathbf{r}_{\text{primal}}$ and $\mathbf{r}_{\text{dual}}$ 
    \If{$\mathbf{r}_{\text{primal}} < \tau_{r}$}   \hfill   $\triangleright$ \textbf{[SC-$2'$]} 
    \State Run backtracking process, $\n' \leftarrow$ columns set that makes the residual closest to  \State the target residual $\tau_r'$, break
    \EndIf
    \State Select rows $\m'$ and columns $\n'$ in batches, $\m \leftarrow \m'$, $\n \leftarrow \n'$
    \State Estimate the condition number $\kappa( \A(\m',\n' ) )$  
    \If{$\kappa(\A(\m, \n)) > \tau_{\kappa}$}       \hfill   $\triangleright$ \textbf{[SC-$1'$]}
    \State Run residual search, $\n'' \leftarrow$ columns set corresponding to the residual minimum 
    \State value 
    \State $\n \leftarrow \n'$
    \State Break
    \EndIf
    
\EndWhile
\end{algorithmic}
\end{algorithm}

\subsection{[SC-$1'$] Stop by large condition number}
\label{sec: Ill-conditioned situation}
\begin{figure}
\centering
  \subfloat[Using selected rows in computation]{\begin{overpic}[width=6.5cm,scale=.25]{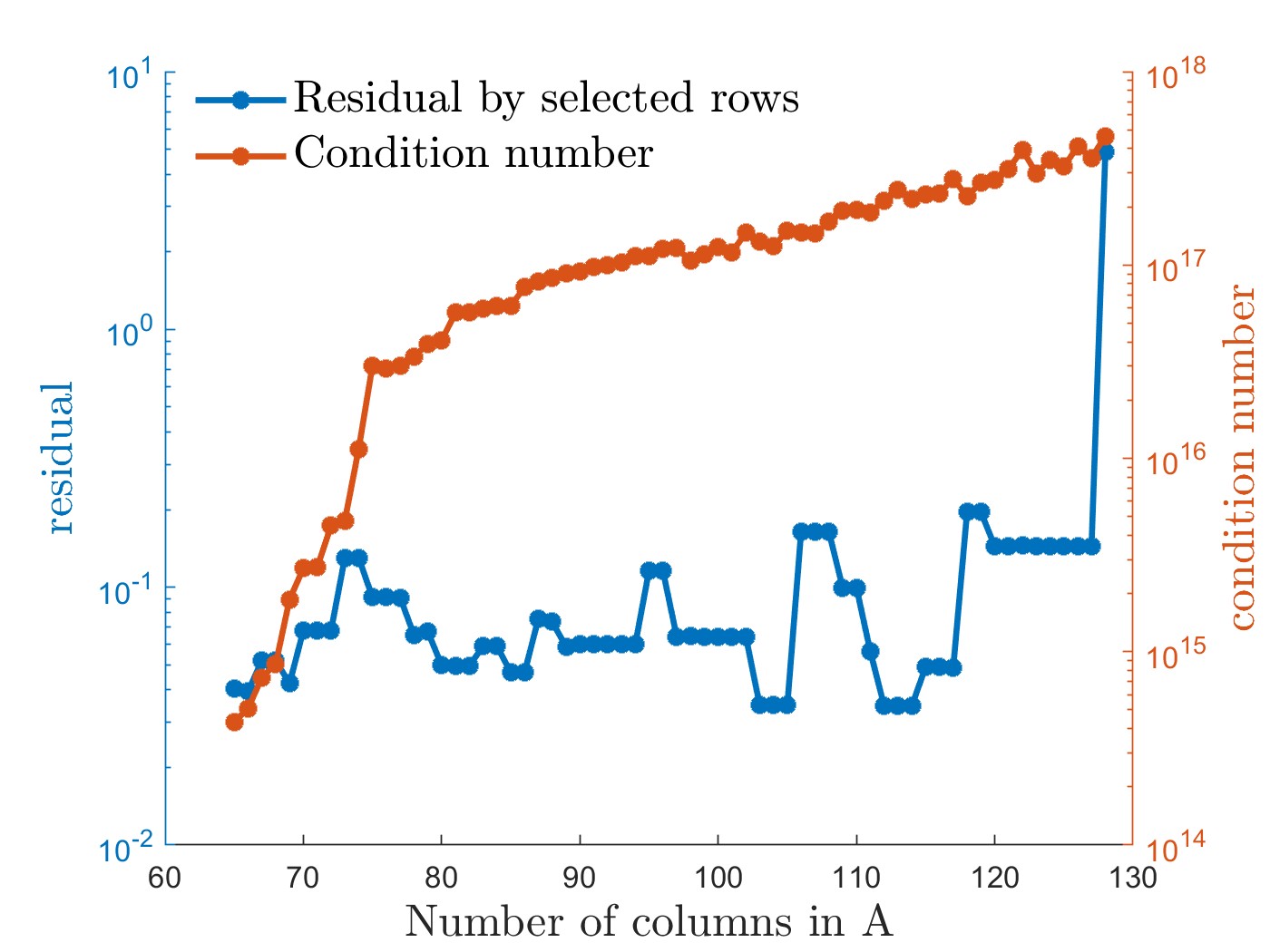}
    \multiput(15,40)(6,0){10}{\linethickness{0.2mm}\line(1,0){3}}
    \put(75,38){$1/\tau_{\kappa}$}
    \end{overpic}
    \label{fig:ResCondXselill}}
  \subfloat[Using all rows in computation]{\begin{overpic}[width=6.5cm,scale=.25]{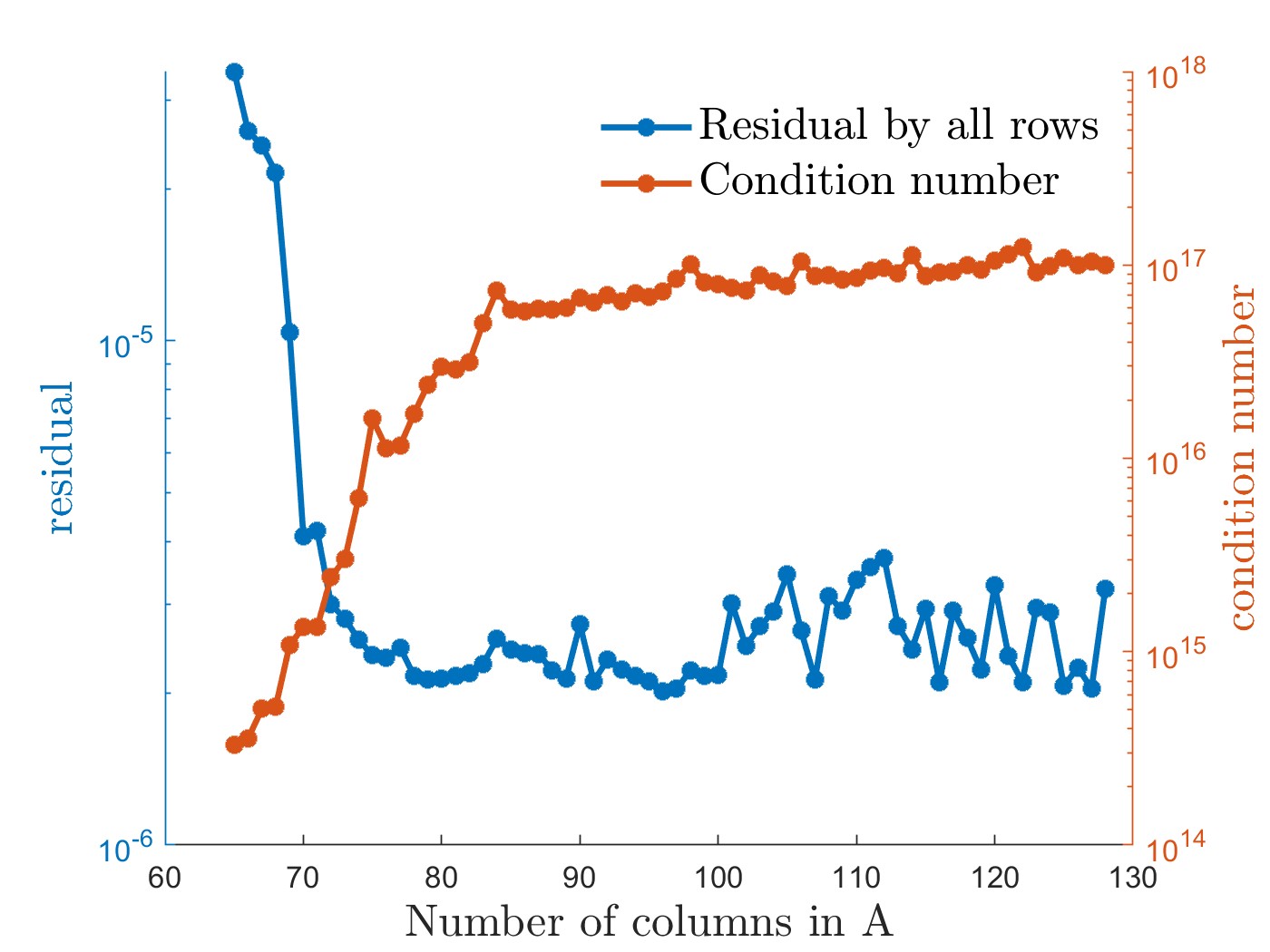}
    \end{overpic}
    \label{fig:ResCondXfullill}}
   \caption{Greedy algorithm stops by large condition number, 128 columns are selected in the column space -- A schematic demonstration of residual value and condition numbers of submatrices of $\A$ with expending columns $\n\to\n'$ and (a) selected rows $\m$ and (b) all rows $\onetom$ in the column space.  
   }\label{fig:Ill Res Cond}
\end{figure}

Figure \ref{fig:ResCondXselill} schematically illustrates the behavior of the residuals and the condition numbers of the subsystem as selected columns are incrementally added. This progression continues until the subsystem reaches ill-conditioning, at which point the original greedy method is stopped by SC-$1$.
The red curve shows the condition numbers of submatrices of $\A$ with rows $\m'$ against expanding columns from the index set $\n$ in the previous iteration to that of the current iteration $\n'$, denoted by $\n\to\n'$.
Its left endpoint satisfies $\kappa(\A(\m,\n))<1/\tau_\kappa$ and so the iteration continues; whereas its right endpoint satisfies $\kappa(\A(\m',\n'))>1/\tau_\kappa$ and so the greedy algorithm stops. In the original SC-$1$, a backtracking process traces back along the same curve. It seeks the point with $\kappa(\A(\m,\n))<1/\tau_\kappa$ to identify the final number of column selection $\greedyk$, as discussed in Section 2.
The blue curve shows the $\ell^\infty$-norm residual vectors of all rows for least-squares subproblems similar to  \eqref{eq p res}  with the selected rows in $\m$ and expanding $\n\to\n'$ columns in $\A$.

We remark that the observed large residual is not the true approximation of the trial subspace.
Large error is expected at the remaining rows in $\onetom\setminus\m'$ whose information was left out from the least-squares subproblems. The greedy method is used for columns (or trial functions) selection and PDEs were solved with all rows (or collocation points).

For the new stopping criterion SC-$1'$, we stop the greedy iteration using the same mechanism:
\[
 \kappa(\A(\m,\n)) < 1/\tau_\kappa \text{\quad $and$ \quad}\kappa(\A(\m',\n')) > 1/\tau_\kappa.
\]
For small $\tau_\kappa$, we are dealing with ill-conditioned subsystems and should take a stability-first approach. The goal of the greedy algorithm is then to select column subspaces that are computationally stable.

Figure \ref{fig:ResCondXfullill} corresponds to Figure \ref{fig:ResCondXselill} but uses all rows in the submatrices. Note the different $y$-axis ranges in these figures. First, using more rows has a minimal effect on the condition numbers and the red curves are similar in both figures. The blue curve shows the $\ell^\infty$-norm error of using all rows and expanding the $\n$ to $\n'$ columns in $\A$ to least-squares approximate $\mathbf{b}$, i.e.,
\[
    \A(\onetom, \n\to\n') \boldsymbol{\lambda} = \mathbf{b}.
\]
This residual curve shows the true approximation power of trial subspaces.
We also observe that the trial subspace that yields the smallest approximation error in approximating $\mathbf{b}$  is not easy to detect from the curve of the condition number.
We propose searching this residual curve instead and having the greedy algorithm return $\greedyk$ columns at which the residual is minimum and the corresponding column index set $\n''$.

The computational overhead is the difference in cost of QR factoring the greedy selected $|\m'| \times |\n''|$ submatrix and the $|\m'| \times |\n'|$ that allows us to compute the whole residual curve in Figure~\ref{fig:ResCondXfullill}. The estimate is $\calo(|\m'| \times (|\n'|^2-|\n''|^2) )$.

\subsection{[SC-$2'$] Stop by small residual}
\label{sec: Well-conditioned large-scale problems}

\begin{figure}
\centering
  \subfloat[Using selected rows in computation]{\begin{overpic}[width=6.5cm,scale=.25]{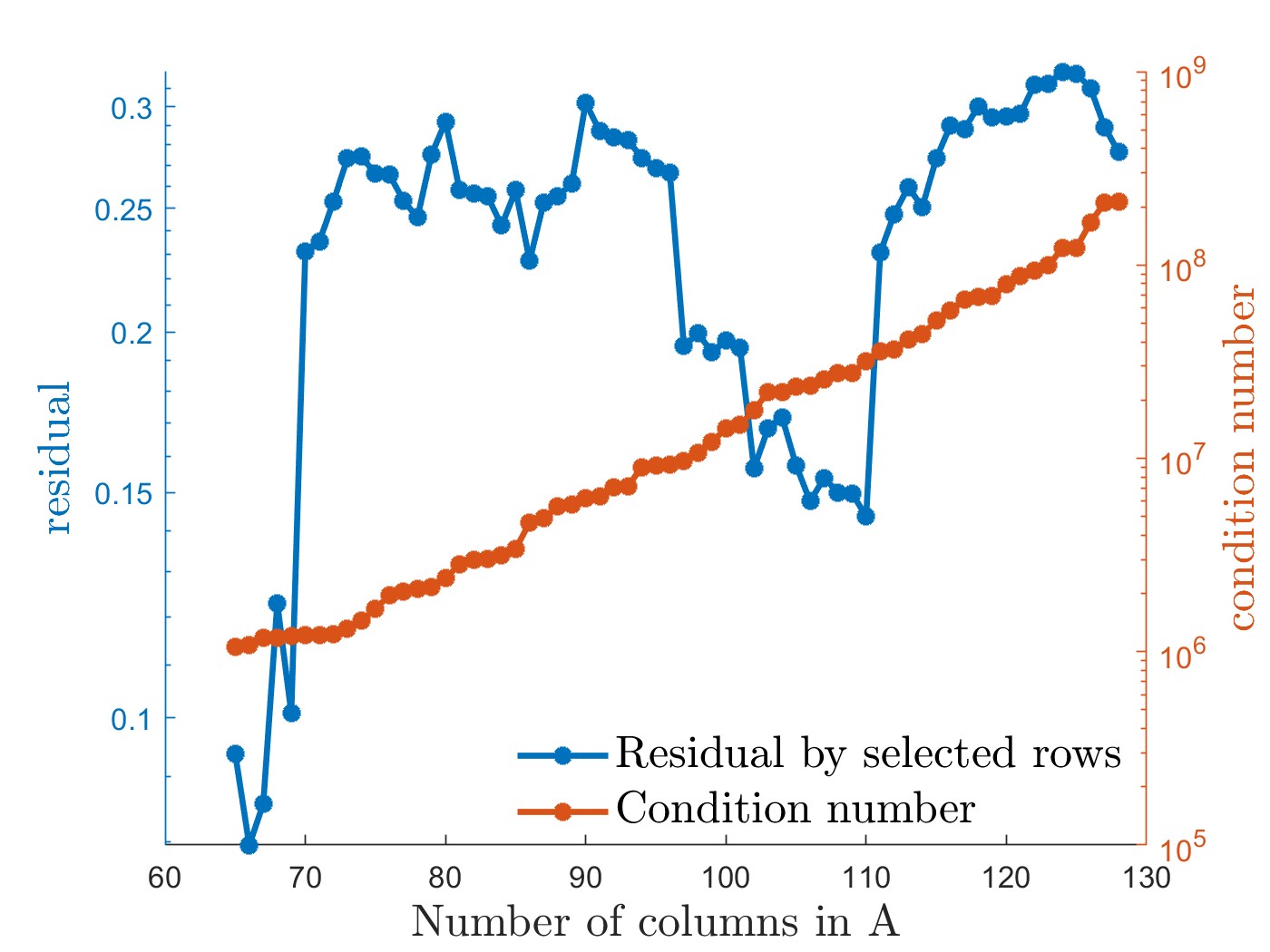}
    \multiput(20,36)(6,0){10}{\linethickness{0.2mm}\line(1,0){3}}
    \put(79,35){$\tau_{r}$}
    \end{overpic}
    \label{fig:ResCondXselwell}}
  \subfloat[Using all rows in computation]{\begin{overpic}[width=6.5cm,scale=.25]{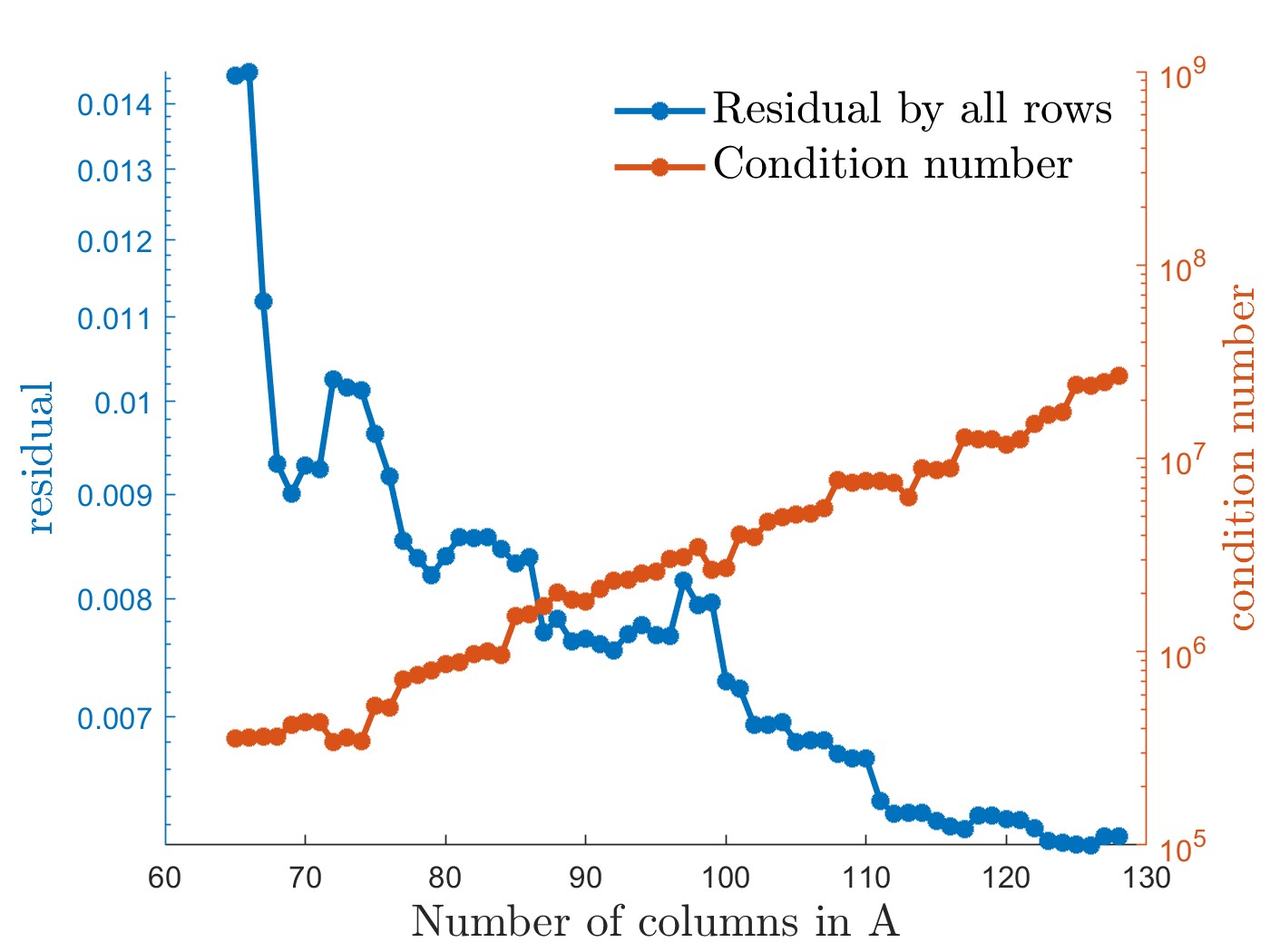}
    \end{overpic}
    \label{fig:ResCondXfullwell}}
   \caption{Greedy algorithm stops by small residual, 128 columns are selected in the column space --A schematic demonstration of residual value and condition numbers of submatrices of $\A$ with expending columns $\n\to\n'$ and (a) selected rows $\m$, and (b) all rows $\onetom$ in the column space. 
         }
\label{fig:well Res Cond}
\end{figure}

This is an aspect the original greedy algorithm largely ignores. The greedy algorithm only runs on residual values using the selected rows and columns, see equations \eqref{eq p res}--\eqref{eq d res}.
Setting $\tau_r=\varepsilon_{mech}$ means the greedy algorithm almost never stops because of this small residual criteria. Even in the case when the discretized matrix system is well-conditioned,
the greedy algorithm in the default setting will select the whole column space (but without a huge overhead in the block version) to yield the best possible accuracy.
Figure \ref{fig:well Res Cond} were the counterpart of Figure \ref{fig:Ill Res Cond} but the greedy algorithm was stopped by a small residual for some given (large) tolerance $\tau_r$. In Figure~\ref{fig:well Res Cond}, the condition numbers are all well below the dangerous zone $1/\varepsilon_{mech}$. In Figure~\ref{fig:ResCondXselwell}, we once again observed that the residual values obtained using the selected rows are large and do not reflect the true approximation power of the trial subspace.
The residual using all rows in Figure~\ref{fig:ResCondXfullwell} is not a quantity that exists within the greedy algorithm. Yet it could be way smaller than the tolerance value $\tau_r$.

Suppose the greedy algorithm selected column index set $\n'$ after stopping by small residual (based on selected rows indexed by $\m'$), $\|  \mathbf{r}_{\text{primal},\m'} \|_\infty < \tau_{r}$. The original greedy algorithm does not have any postprocessing.
In the new SC-$2'$, we propose to run the stopping criteria using residual from all rows when
\[
     \|  \mathbf{r}_{\text{primal},\onetom} \|_\infty < \tau_{r} .
\]
That is, we use the blue curve in Figure~\ref{fig:ResCondXfullwell} instead of that in Figure~\ref{fig:ResCondXselwell} to stop the greedy iteration.
Then we apply a  backtracking process, similar to the one presented in Section \ref{sec greedy review},  to select the first $\greedyk$ indices in $\n'$ to generate the largest subset $\n'' \subset \n'$ so that the final column subspace selection yields a residual for the  $\A(\onetom, \n'')$ subsystem is just below another tolerance $\tau_r'$.
If the residual using all rows and $\n'$ is above $\tau_r'$, we simply set $\n''=\n'$ and return all selected columns.
Like in SC-$1'$, the QR factorization of $\A(\onetom, \n')$ allows us to evaluate any data point on the residual curve. Since we double the number of selected columns in each (block-version) greedy iteration, the overhead of SC-$2'$ is the extra QR factorization to intermediate submatrix $\A(\onetom,1\to\n)$ in the previous iteration that costs $\calo(m\greedyk^2)$.

In the context of meshfree time-stepping schemes, a highly accurate spatial scheme is wasteful in terms of computational cost when the overall error is dominated by temporal discretization.
It makes sense to greedily select just enough spatial approximation power to match the temporal error in order to reduce spatial accuracy and compensate for computational efficiency.
The tolerance $\tau_r'$ should be based on the finite difference scheme and time step $\dt$ used in temporal discretization, which yields the semi-discrete systems. For the second-order scheme in this paper, we propose:
\begin{equation}\label{eq tau}
  \tau_\kappa=\calo(\mecheps/\dt ),\quad\tau_r=\calo(\dt ), \quad\text{and}\quad\tau_r'=\calo(\dt ^2).
\end{equation}
In the rest of this paper, we use 1 as the Big-Oh constant, yielding satisfactory results in 2D and 3D test problems.

\begin{example}\label{eg: 2D RD PDE new}\textbf{(New stopping criteria in action)}
\begin{figure}
\centering
  \subfloat[$\epsilon = 1$]{\includegraphics[width=4.5cm]
 {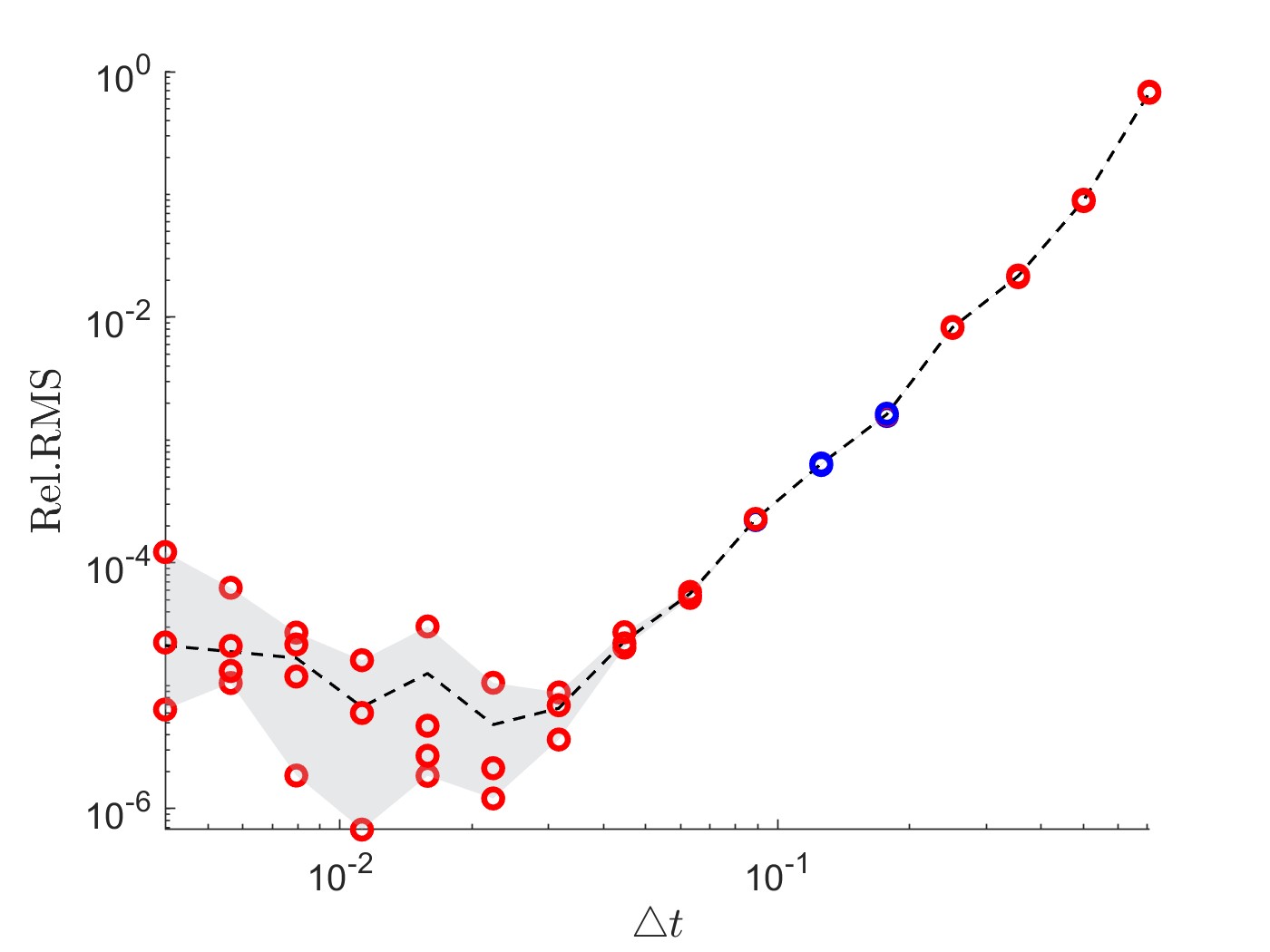}\label{fig:linerr2DMS61}}
  \subfloat[$\epsilon = 3$]{\includegraphics[width=4.5cm]{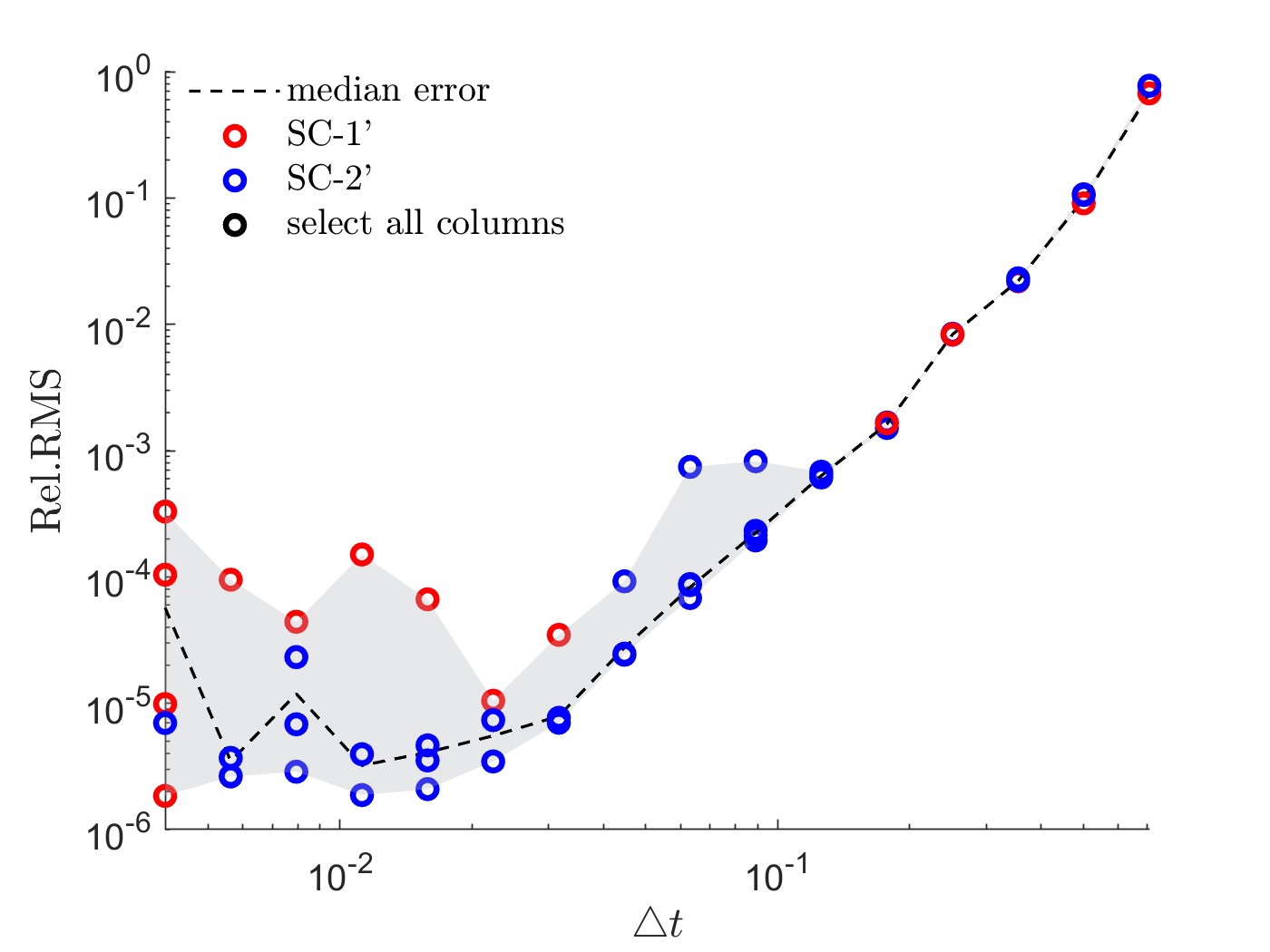}\label{fig:linerr2DMS63}}
  \subfloat[$\epsilon = 6$]{\includegraphics[width=4.5cm]{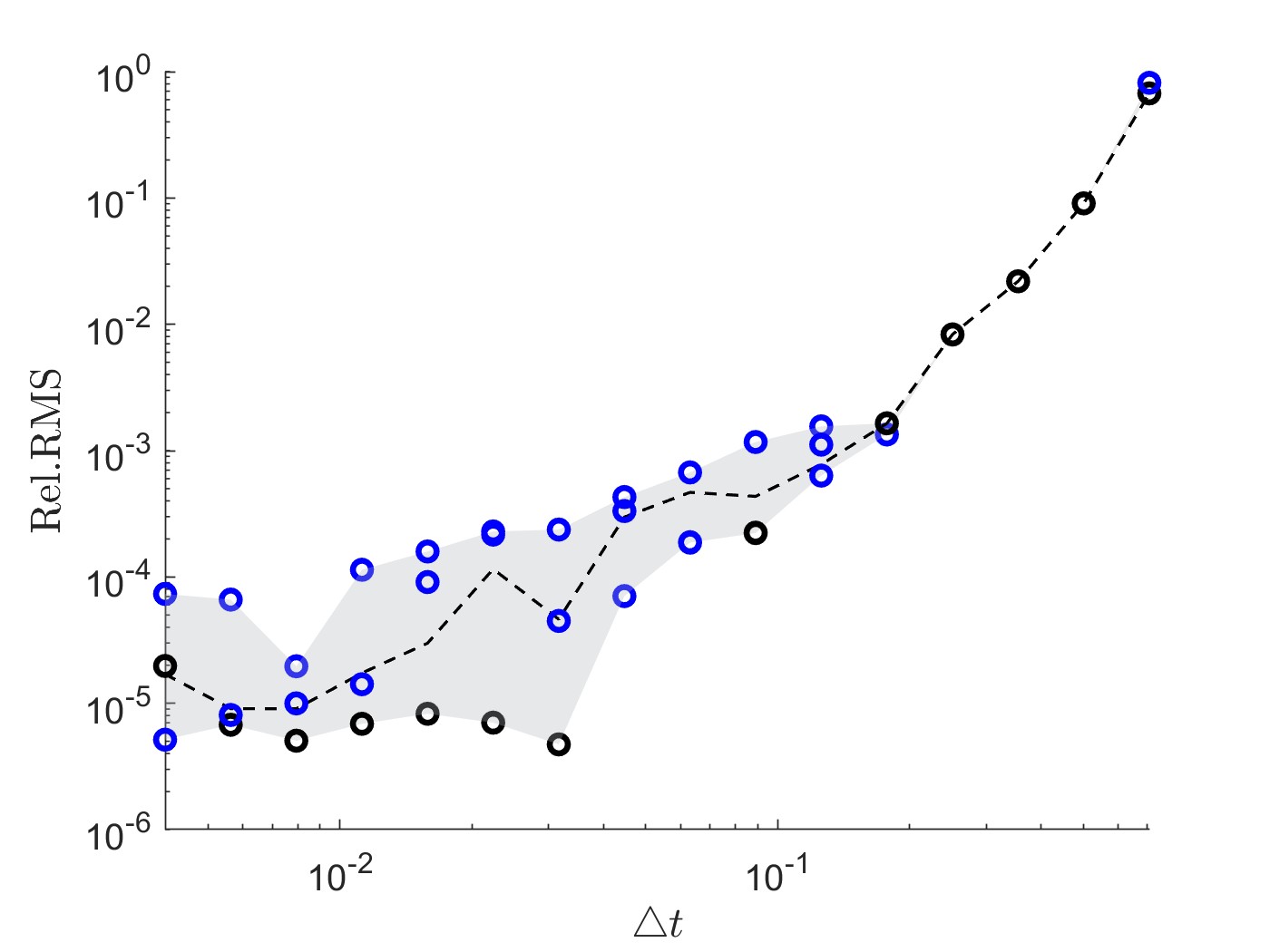}\label{fig:linerr2DMS66}}
   \caption{
 Example~\ref{eg: 2D RD PDE new}, error profiles for solving a 2D heat equation using a meshfree time-stepping method with greedy algorithm and the MS kernel with $\mu=6$.
  Different values of $\varepsilon$ and $\dt $ were used, and the shaded areas in each plot show the error range for $n\in[500,1000]$ data points, with the median as a dashed line.
The stopping criteria that terminate the greedy algorithm are shown as black circles for all columns being selected, \red{red circles} for SC-$1'$, and \blue{blue circles} for SC-$2'$.
         }
\label{fig:linerr2DMS}
\end{figure}
We utilize SC-$1'$ and SC-$2'$ into the greedy algorithm and solve the heat equation in Example \ref{eg: 2D RD PDE} again. For temporal discretization, we still use the CN scheme. For spatial discretization, we use the MS kernel in dimension $d=2$ with smoothness order $\mu=6$
\begin{equation}\label{eq ms6}
  \Phi_\mu(x-y) := \|x-y\|_2^{\mu-d/2} \calk_{\mu-d/2}(\|x-y\|_2)
  \qquad\text{for  $x,y\in\R^2$}
\end{equation}
and scale input argument of  this translation-invariant kernel  by $\Phi(r)\leftarrow \Phi(\epsilon r)$ with shape parameters $\epsilon=1$, 3, and 6.
We take $n=500,550,\ldots,1000$ Halton points to form the sets of trial centers and collocation points.

Instead of one curve per tested $n$, we show the range of relative root mean errors in Figure (\ref{fig:linerr2DMS}) along with the median of all tested $n$ against different $\dt$. Note that the $\dt$ value determines all greedy tolerance by \eqref{eq tau}. The red circle for SC-$1'$ and blue circle for SC-$2'$ in the figures indicate the active stopping criteria that terminate the greedy algorithm in that particular run. If the greedy algorithm selects all columns without activating any stopping criteria (usually when $n$ is small and the system matrix $\A$ is well-conditioned), it is labelled by black circles in the figure.
Our experimental results show that the new stopping criteria generate more stable solutions with a wide range of meshfree setups, as evidenced by the fluctuation range of the error in Figure \ref{fig:linerr2DMS}.
In the case of $\epsilon=1$,  the  MS kernel quickly leads to an ill-conditioned matrix in the fully-discretized system as $n$ increases. We observe in Figure~\ref{fig:linerr2DMS61} that most runs were stopped SC-$1'$ by large conditioned numbers. With everything else fixed, increasing $\epsilon$ yields a system matrix with a smaller condition number. We choose to show the results for $\epsilon=3$ since it clearly demonstrates errors of similar magnitude when greedy is terminated by either stopping criteria.
The errors using all columns are generally smaller since SC-$2'$ is designed for the sake of computational cost. Selecting all columns yields the highest possible approximation power in the full discretization and the lowest error, though at a higher computational cost.

We set up a similar test problem in 3D with the following heat solution:
\[
u^{*}([x_1,x_2,x_3]^T, t) = \sin(2\pi x_{1}) \sin(\pi x_{2})\sin(\pi x_{3})  \exp{(-2\pi^2t)} + (1 - \exp{(-2\pi^2t)})/2.
\]
Using the MS kernel \eqref{eq ms6} with $d=3$ and $\epsilon=3$ centered at $n=5000$ to $8000$ Halton points in $[0,1]^3$ as trial centers and collocation points, we solve the heat equation (as in the 2D case) with various time step $\dt$ (hence, various tolerance for the greedy algorithm). The resulting relative root mean squared errors are shown in Figure~\ref{fig:3d heat}a.
Similar to Figure~\ref{fig:linerr2DMS63}, we see that the blue circle for SC-$2'$ is more likely to yield smaller errors for any $\dt$. Figure~\ref{fig:3d heat}b shows the error function for $n = 7000$ and $\dt = 0.01$ with the greedy algorithm.

\begin{figure}
  \centering
  \begin{overpic}[width=0.45\textwidth,trim=0 0 0 0, clip=true,tics=10]{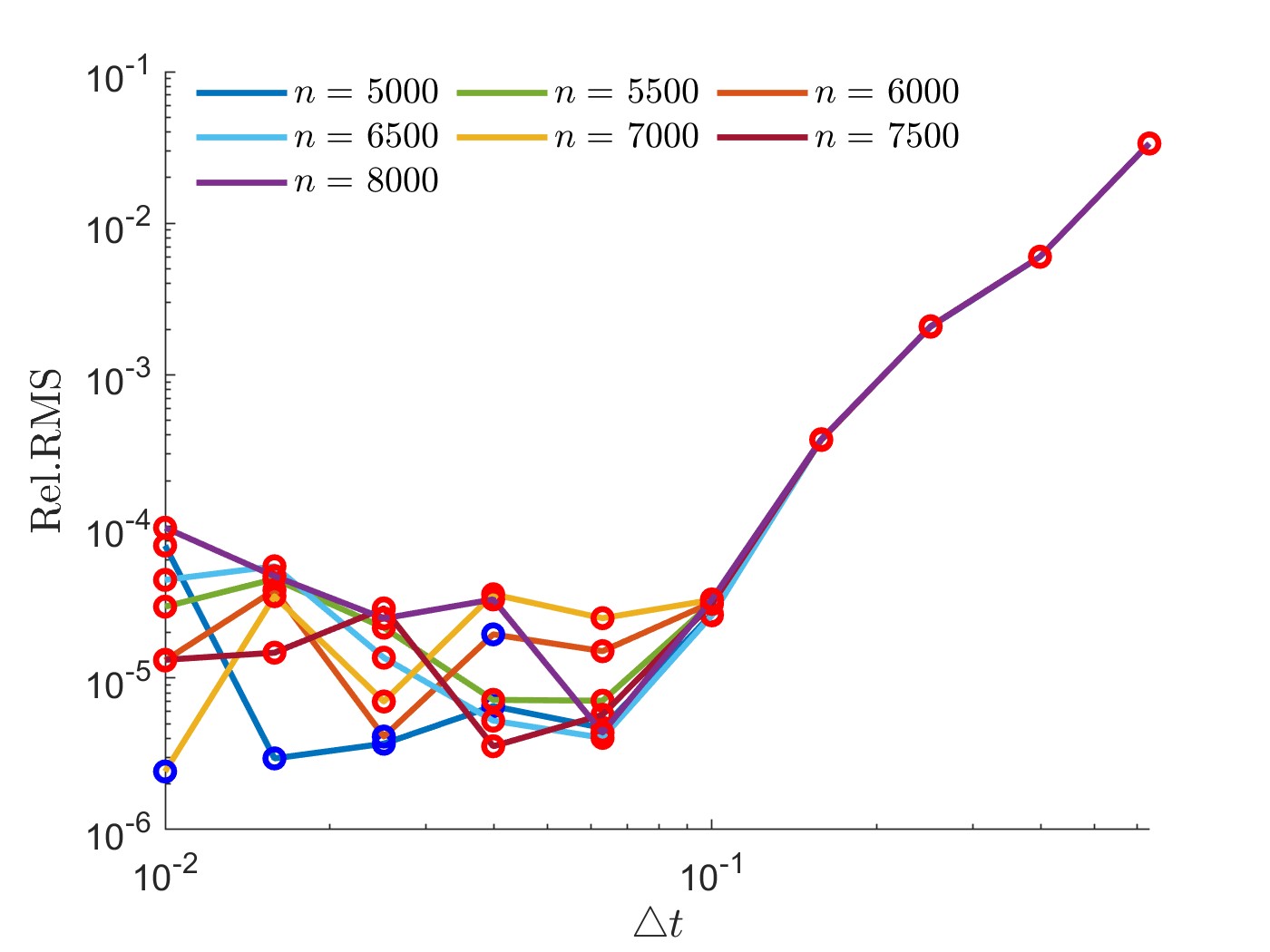}
    \put(45,-3){\scriptsize (a)}
\end{overpic}
\begin{overpic}[width=0.45\textwidth,trim=50 50 50 50, clip=true,tics=10]{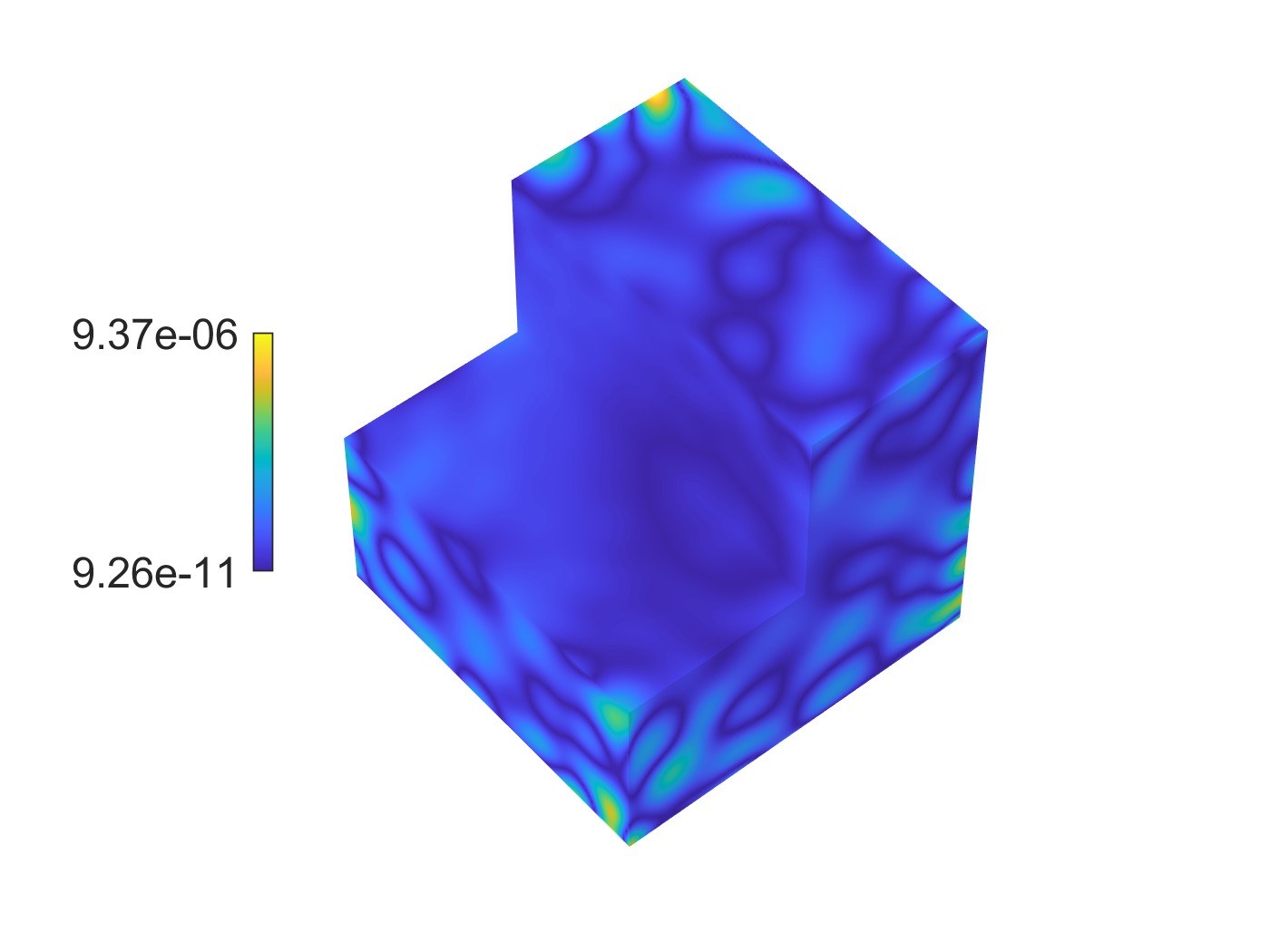}
    \put(45,-3){\scriptsize (b)}
\end{overpic}
\caption{Example~\ref{eg: 2D RD PDE new}, (a) error profiles for the solution of a 3D heat equation using a meshfree time-stepping method with the new greedy stopping criteria for different numbers of data points $n$. (b) error function for $n = 7000$, $\dt=0.01$.
   }\label{fig:3d heat}
\end{figure}
\qed
\end{example}
The previous examples have provided motivation and insights into the effectiveness of the new stopping criteria across different setups. However, a more thorough investigation through extensive numerical experiments is needed to verify the robustness of our proposed stopping criteria and tolerance selection strategies.

\section{Coupled bulk-surface pattern formations}\label{sec num exmp}
Semi-linear coupled bulk-surface reaction-diffusion equations are crucial in modeling phenomena where processes at the surface significantly influence or are influenced by the dynamics occurring in the bulk. Such models find applications in various fields including biology, chemistry, materials science, and environmental engineering. These equations consist of reaction-diffusion terms in the bulk and on the surface, with coupling typically occurring through boundary conditions linking the fluxes or concentrations between the bulk and the surface. For example, in developmental biology, it is crucial to understand how polarized states arise and are maintained, characterized by uneven distributions of chemical substances such as proteins and lipids \cite{Symmetrybreaking2014, Rtz2011TuringII}. Another application of the bulk-surface reaction-diffusion equation is in developing solid-state hydrogen storage materials using metal hydrides, where slow kinetics of hydrogen uptake and release involve adsorption and dissociation of hydrogen, its diffusion through the solid, and a phase transition in the metal-hydride from lower to higher hydrogen content \cite{2022IJMMM..29...32L}. This section explores the application of these equations in developmental biology, specifically in the area of pattern formation.

Let $\Omega\subset\R^d$ be some bounded (bulk) domain and $\cals=\partial \Omega$ be its smooth boundary with well-defined normal vector field $\nn:\cals\to\R^d$. We verify the new greedy stopping criteria with the following semi-linear coupled bulk-surface reaction-diffusion equations in both 2D and 3D for pattern formation:
\begin{equation}\label{eq: CBS}
\begin{split}
&
\begin{dcases}
\frac{\partial u}{\partial t} = D_{u}\nabla^{2}u +  f_{1}(u, v) ,  \\
\frac{\partial v}{\partial t} = D_{v} \nabla^{2}v + f_{2}(u, v) ,   \\
\end{dcases} %
\hskip 11.3em\relax
\textrm{in } \Omega \times [0, T]\\
& \begin{dcases}
\frac{\partial w}{\partial t} = D_{w}\nabla_{\cals}^{2}w +  f_{1}(w, s) - h_{1}(u, w),
 \\
\frac{\partial s}{\partial t} = D_{s}\nabla_{\cals}^{2}s +  f_{2}(w, s) - h_{2}(v, s),
\end{dcases}
\hskip 6.0em\relax
\textrm{on } \cals \times [0, T]
\end{split}
\end{equation}
with coupling boundary conditions
\begin{equation}\label{eq: bc CBS}
\begin{dcases}
D_{u}\nabla u \boldsymbol\cdot \nn = h_{1}(u, w),  \\
D_{v}\nabla v \boldsymbol\cdot \nn = h_{2}(v, w).
\end{dcases}\hskip 13em\relax
\textrm{on } \cals \times [0, T].
\end{equation}
The surface equations in \eqref{eq: CBS} contain differential operators on surfaces, that can be defined via a projection to the tangent space of $\cals$
\[
P : = I_{d} - \nn\nn^T,
\]
where $I_{d}$ denotes the $d \times d$ identity matrix. The surface gradient and the Laplace-Beltrami operator, a.k.a. surface Laplacian, are defined respectively by
\[
    \nabla_{\cals}= P\nabla
    \quad\text{and}\quad
     \nabla_{\cals}^{2} = \nabla_{\cals} \boldsymbol\cdot \nabla_{\cals} = (P \nabla) \boldsymbol\cdot (P \nabla),
\]
using the standard gradient operator $\nabla$ for functions defined in $\R^d$.

\begin{table}
\caption{Parameters for Turing's pattern formation used in coupled bulk-surface reaction-diffusion equations \eqref{eq: CBS}--\eqref{eq: bc CBS}.}\label{ta: parameters for CBS}
\begin{center}
\begin{tabular}{|c||c | c| c| c| c| c| c |c| c |c| c |c|}
 \hline

& $a$  & $b$ & $\alpha_{1}$ & $\alpha_{2}$ & $\beta_{1}$ & $\beta_{2}$ &  $D_{v}$ & $D_{s}$ & $q$ & $\gamma$  \\ \hline
Fig.~\ref{fig: CBS2D w/o g demo}--\ref{fig: CBS2D w g} 
& \multirow{ 2}{*}{\vdots} & \multirow{ 2}{*}{\vdots} &  \multirow{ 2}{*}{\vdots} & \multirow{ 2}{*}{\vdots} & \multirow{ 2}{*}{\vdots} & \multirow{ 2}{*}{\vdots} &  2 & 2 & 1/12 & 30
\\
Fig.~\ref{fig: CBS2D w g} 
&  &  &   &  &  &  &  1 & 1 & 1/12 & 30
\\
Fig.~\ref{fig: CBS2D w g stripe} 
& \multirow{ 2}{*}{1/10} & \multirow{ 2}{*}{9/10} & \multirow{ 2}{*}{5/12}  & \multirow{ 2}{*}{5} & \multirow{ 2}{*}{5/12} & \multirow{ 2}{*}{5} &  {5} & {5} & {1/10} & {500}
\\
Fig.~\ref{fig: CBS3D torus eps=1}--\ref{fig: CBS3D torus eps=4 5}
&  &  &   &  &  & &  {3} & {3} & {1/12} & {40}
\\
Fig.~\ref{fig: CBS3D cycl}
& \multirow{ 2}{*}{\vdots} & \multirow{ 2}{*}{\vdots} & \multirow{ 2}{*}{\vdots}  & \multirow{ 2}{*}{\vdots} & \multirow{ 2}{*}{\vdots} & \multirow{ 2}{*}{\vdots} &  {6} & {6} & {1/12} & {30}
\\
Fig.~\ref{fig: CBS3D ellip}
&  &  &   &  &  & &  {3} & {3} & {1/12} & {30}
\\
\hline
\end{tabular}
\end{center}
\end{table}

The source functions in \eqref{eq: CBS}--\eqref{eq: bc CBS} are given by the nonlinear reaction kinetics and the coupling of the internal bulk dynamics in $\Omega$ to the surface dynamics on $\cals$. The interaction between unknown functions $u,v:\Omega\to\R$ and $w,s:\cals\to\R$ yields Turing patterns in the bulk and on the surface.
The nonlinear reaction kinetics are:
\begin{eqnarray*}
f_{1}(u, v) = \gamma(a - u + u^2v), & \qquad& f_{2}(u, v) = \gamma(b - u^{2}v).
\end{eqnarray*}
The coupling of the internal dynamics to the surface dynamics is:
\begin{eqnarray*}
h_{1}(u, w)= \alpha_{1}w - \beta_{1}u, & \qquad & h_{2}(v, s) = \alpha_{2}s - \beta_{2}v,
\end{eqnarray*}
with parameter values listed in Table \ref{ta: parameters for CBS} with   $D_u = qD_v$, $D_w = qD_s$.
As shown in \cite{BSMz}, the coupled bulk-surface PDEs \eqref{eq: CBS} has a unique homogeneous equilibrium state when
\begin{equation}
    \label{eq: steady state CB}
    f_{1}(w, s) - h_{1}(u, w) = 0 \quad \textrm{and} \quad f_{2}(w, s) - h_{2}(v, s) = 0.
\end{equation}
For any $\gamma$, solving \eqref{eq: steady state CB} yields corresponding (constant) function values at equilibrium as
\begin{equation}\label{eq bspf ic}
  \Big(a+b,\, \frac{b}{(a+b)^2},\,a+b, \,\frac{b}{(a+b)^2}\Big) = (u_0,\, v_0,\, w_0,\, s_0).
\end{equation}
We set the initial conditions $(u_0,\, v_0,\, w_0,\, s_0)$ for the coupled bulk surface PDEs \eqref{eq: CBS} to be the homogeneous equilibrium state with $a,b$ as shown in Table \ref{ta: parameters for CBS}. Small approximation errors in these initial conditions then serve as sources of perturbations. According to Turing instability theory \cite{Turing1}, in the absence of diffusion a uniform steady state remains stable. However, it becomes unstable to small spatial perturbations. In our meshfree context, minor disturbances from numerical approximation error can cause an otherwise stable homogeneous equilibrium state to become unstable, resulting in the formation of patterns.

\subsection{Discretizing PDEs for the greedy algorithm}
We use the SBDF2 scheme in equation \eqref{eq: sbdf2} for temporal discretization, where all nonlinear reaction terms in equation \eqref{eq: CBS} were evaluated at the previous time step. Let $U^k,V^k,W^k,S^k$ represent the semi-discretized solutions as in Section \ref{sec select tol}. We also use an explicit scheme for bulk-surface dynamics.
The semi-discretized equation for the bulk function $u$ is given by:
\begin{equation}\label{eq BSPF bulk}
\frac{1}{2\dt _{k}}\Big(3\u^{k} -4\u^{k-1} + \u^{k-2}\Big) = 2f_1(U^{k-1},V^{k-1}) - f_1(U^{k-2},V^{k-2}) + D_u\nabla^{2}\u^{k}
\quad\text{in $\Omega$,}
\end{equation}
subject to the boundary condition
\begin{equation}\label{eq BSPF bc}
   D_{u}\nabla U^k \boldsymbol\cdot \nn = h_1( U^{k-1},W^{k-1})
   \quad\text{on $\cals$,}
\end{equation}
which is decoupled from the other solutions.
The semi-discretized equation for the surface function $w$ is
\begin{eqnarray}
  \frac{1}{2\dt _{k}}\Big(3W^{k} -4W^{k-1} + W^{k-2}\Big) = 2f_1(W^{k-1},S^{k-1}) - 2 h_1(W^{k-1},S^{k-1})
  \nonumber    \\
    -f_1(W^{k-2},V^{k-2}) +  h_1(W^{k-2},V^{k-2})
+ D_w\nabla^{2}_\cals W^{k}
\quad\text{on $\cals$.}     \label{eq BSPF surf}
\end{eqnarray}
We obtain the semi-discretized systems for $V^k$ and $S^k$ similarly to the equation for $U^k$ and $W^k$. In the end, this yields four sequences of linear elliptic PDEs that incorporate solution history only through the right-hand source functions.

For spatial discretization, we construct two trial spaces to approximate the solutions in the bulk and on the surface, respectively, using the MS kernel in \eqref{eq ms6} that reproduces the Sobolev space $H^\mu(\R^d)$. Let $\mu_\Omega \geq \max(2,d/2)$ be the smoothness order used in the bulk kernel $\Phi_{\mu_\Omega}$.
When we restrict the kernel $\Phi_{\mu_\cals + 1/2}$ on the surface $\cals$ for any $\mu_\cals \geq \max(2,(d-1)/2)$, the restricted kernel $\Psi_{\mu_\cals} := \Phi_{\mu_\cals + 1/2}\big|\cals : \cals \times \cals$ is a Sobolev space $H^{\mu_\cals -1/2} (\cals)$ reproducing kernel with smoothness $\mu_\cals$, as shown in \cite{Fuselier+Wright-ScatDataInteEmbe:12}.
Based on convergence estimates for meshfree least-squares collocation for elliptic PDEs in the bulk \cite{Cheung+LingETAL-leaskerncollmeth:18} and on the surface \cite{Chen+Ling-Extrmeshcollmeth:20}, it is suggested \cite{CMBS} to choose the smoothness orders of the kernels such that
\begin{equation}
    \label{eq: smoothness order}
    \begin{dcases}
    \text{[SO-$1$] }\qquad\mu_{\cals} \geqslant \mu_\Omega + d/2 - 2,
    & \mu_\Omega \geq (9+d)/2,
    \qquad\text{or}
     \\
    \text{[SO-$2$] }\qquad\mu_\Omega \geqslant \mu_{\cals} -d/2,
    & \mu_\cals \geq 3+d,
    \end{dcases}
\end{equation}
in order to balance error estimates in the bulk and on the surface. Condition SO-$1$ arises from the error estimate for meshfree least-square collocation methods solving elliptic PDEs in the bulk and the error of surface functions treated as perturbations. A second condition results from reversing the roles of the bulk and surface functions within the same framework.

To simplify notations, we use two sets of data points $\Xi_\Omega\subset\overline{\Omega}$ and $\Xi_\cals\subset\cals$ for the bulk and surface functions, respectively, i.e.,
\[
    U^k,V^k \in \calu_{\Xi_{\Omega}, \Omega} := \Span{\{\Phi_{\mu_\Omega}(\boldsymbol\cdot\,, \xx)\, |\, \xx \in  \Xi_\Omega\}},
\]
and
\[
    W^k,S^k\in\calu_{\Xi_{\cals}, \cals} := \Span{\{\Psi_{\mu_\cals}(\boldsymbol\cdot\,, \xx)\, |\, \xx \in  \Xi_\cals\}}
    = \Span{\{\Phi_{\mu_\cals + 1/2}(\boldsymbol\cdot\,, \xx)\, |\, \xx \in  \Xi_\cals\}}.
\]
For the upcoming 2D and 3D simulations, we simply use $\mu_\Omega = 6$ and $\mu_\cals =5.5$ which satisfy the first condition SO-$1$ in \eqref{eq: smoothness order}. Consequently, we will use the MS kernel $\Phi_6$ in constructing both the bulk and surface collocation systems.

We use sets of collocation points $Z_\Omega \subset \Omega\cup\cals$ and $Z_\cals\subset \cals$ for bulk and surface functions respectively to set up the meshfree collocation matrix system, as discussed in Section~\ref{sec select tol}.
Specifically, for bulk functions $U^k, V^k$, we can obtain matrix equations similar to \eqref{eq: alam = b} except with different constant factors and a new boundary part in the matrix, namely $\Phi(Z_\Omega\cap\cals,\Xi_\Omega) \leftarrow [\nabla\Phi_{\mu_\Omega}\boldsymbol\cdot\nn](Z_\Omega\cap\cals,\Xi_\Omega)$.

We assume prior analytical knowledge of the normal vector field $\nn$ and use the convergent Kansa-type analytical projection method in \cite{Chen+Ling-Extrmeshcollmeth:20} to approximate the surface Laplacian operator. This extrinsic approach allows collocations directly on surfaces and can fully discretize the surface equations.
The elliptic PDE in \eqref{eq BSPF surf} for semi-discretized surface function $W^k$ leads to the following meshfree collocation matrix when fully-discretized:
\[
  A_w(\dt_k) =\Big[ \big( 3 \Phi_{\mu_\cals + 1/2} - 2\dt _k [D_w\nabla^2_\cals\Phi_{\mu_\cals + 1/2}]\big)(Z_\cals, \Xi_\cals) \Big],
\]
depending on the step size.
The meshfree collocation matrix for $V^k$ looks almost identical but with a different diffusion coefficient $D_v$.
If we use an equal time step $\dt$, then all four meshfree collocation matrices (denoted by $\A_U$, $\A_V$, $\A_W$, and $\A_S$) remain fixed for all time, which is a cost-effective approach that we will take in the following calculations.
Now we have a sequence of matrix equations in the form of \eqref{eq Ax=b} for each unknown functions sequence $U^k$, $V^k$, $W^k$ and $S^k$ for $k=1,\ldots, K$.

Because our initial conditions \eqref{eq bspf ic} correspond to the homogeneous equilibrium state, analytically applying SBDF1 at the first time step yields constant functions $U^1, V^1, W^1, S^1$. This implies that, in the second time step, the right-hand side vector of the SBDF2 matrix system is a constant vector. In this case, we can run the greedy algorithm to select columns for each of the four unknown functions without needing to spend time computing the actual right-hand side vector. It is worth noting that the original greedy algorithm also runs with an all-one right-hand vector as part of the patch test procedure in the absence of a specified input vector.
Next, we obtain solutions in the second time step by solving the SBDF2 matrix system with all rows but only the greedily selected columns. This is where the perturbation comes into play. Solving the reduced system (or the full system in cases without the greedy algorithm) introduces small errors that serve as perturbations, allowing patterns to form according to the Turing instability mechanism described earlier.

In summary, here is the main numerical setup described above:
We use the MS $\Phi_6$ kernel in \eqref{eq ms6} with smoothness order $\mu_\Omega=6$ and $\mu_\cals=5.5$.
The tolerances of the greedy algorithm are set according to \eqref{eq tau} with a Big-Oh constant of 1.
All examples use an exactly determined setup with identical sets of trial centers $\Xi_{\boldsymbol\cdot}$ and collocation points $Z_{\boldsymbol\cdot}=\Xi_{\boldsymbol\cdot}$ in bulk domain and surface, respectively.
The number of data points $n=|\Xi|$ varies for different examples due to different domain volumes and will be reported individually.

\begin{example}\label{eg: CBSRD 2}\textbf{(Robustness verifications in 2D)}
Let $\Omega$ be the 2D unit ball. We solve the coupled bulk-surface pattern formation partial differential equations   \eqref{eq: CBS} in bulk $\Omega$ and surface $\cals=\partial\Omega$ subject to the boundary condition   \eqref{eq: bc CBS} and initial condition   \eqref{eq bspf ic}.
Under this simple setting, even in the absence of an analytic solution, we can use symmetry to gain insight into the interfacial patterns and identify pattern defects in the numerical solutions.
Such interfacial patterns would not form in isolated bulk or surface systems; patterns in the bulk and on the surface are expected to be highly correlated.

\begin{figure}
  \centering
  \begin{overpic}[width=0.27\textwidth,trim=100 50 100 50, clip=true,tics=10]{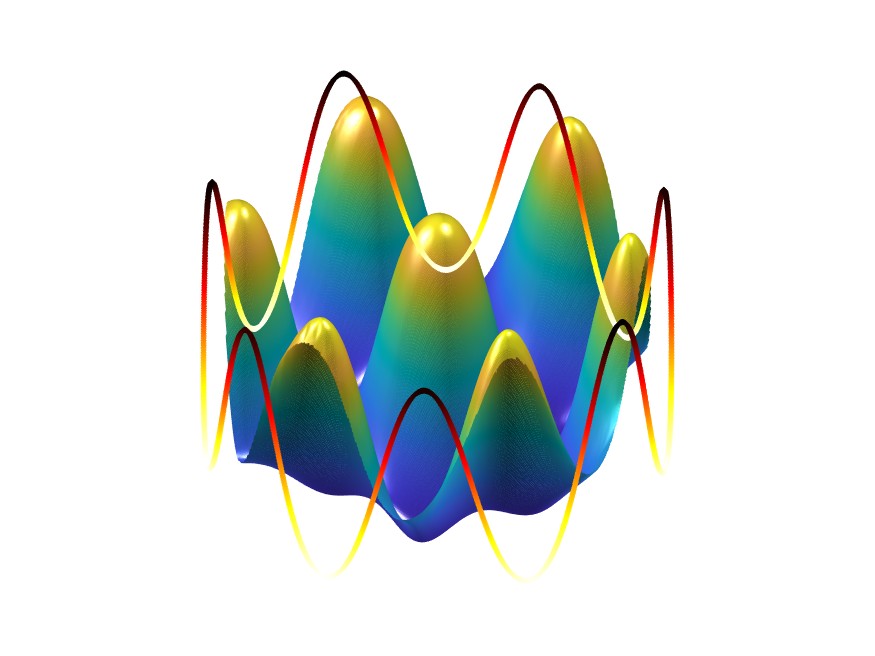}
    \put(-20,25){\rotatebox{90}{\textbf{Without}}}
    \put(-10,28){\rotatebox{90}{\textbf{Greedy}}}
    \put(10,-5){\scriptsize (a) $n=717+100$, $\dt=0.005$}
\end{overpic}
\begin{overpic}[width=0.27\textwidth,trim=100 0 100 0, clip=true,tics=10]{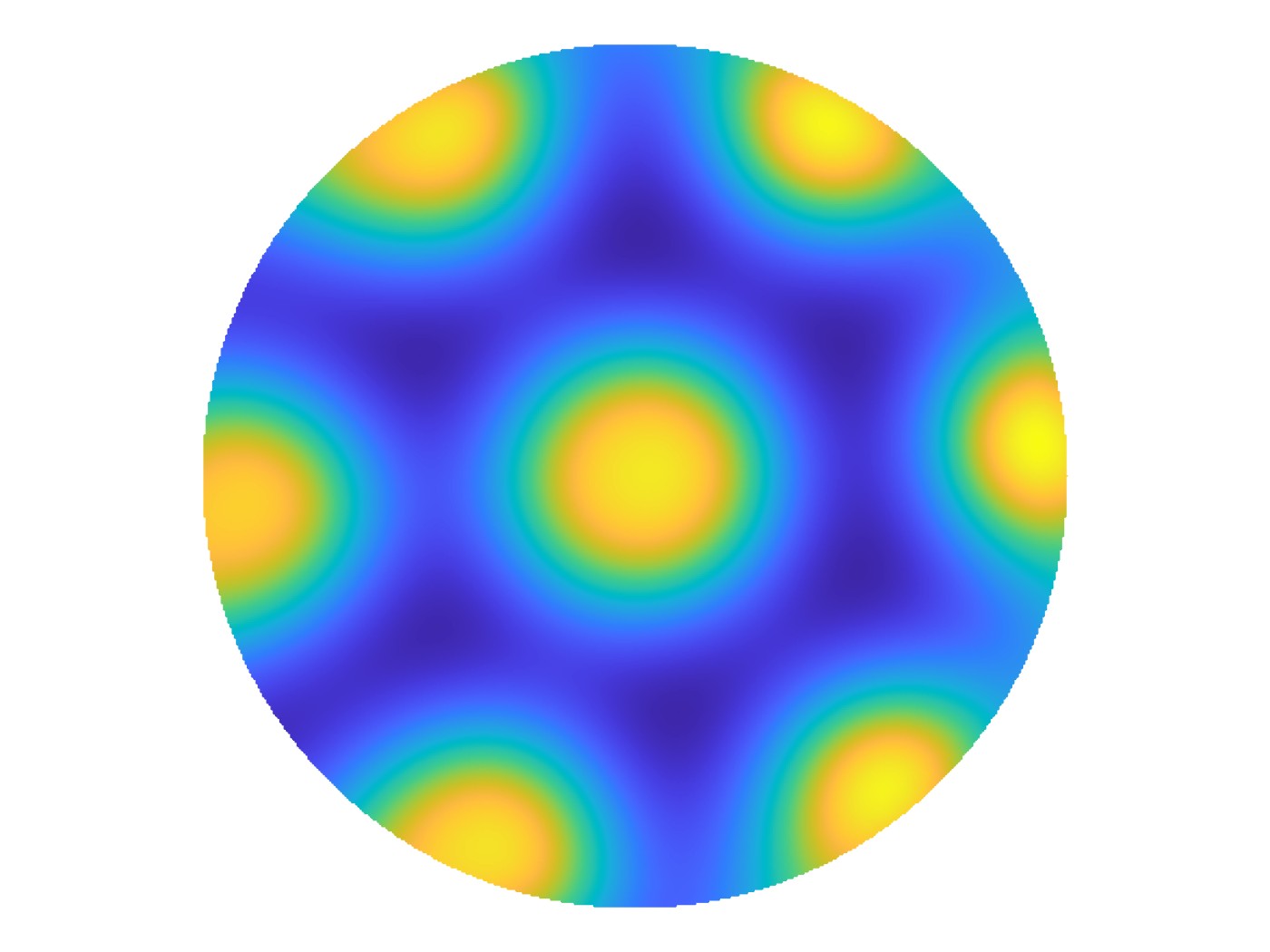}
    \put(20,-5){\scriptsize (b) Bulk solution $u$}
\end{overpic}
\begin{overpic}[width=0.27\textwidth,trim=100 0 100 0, clip=true,tics=10]{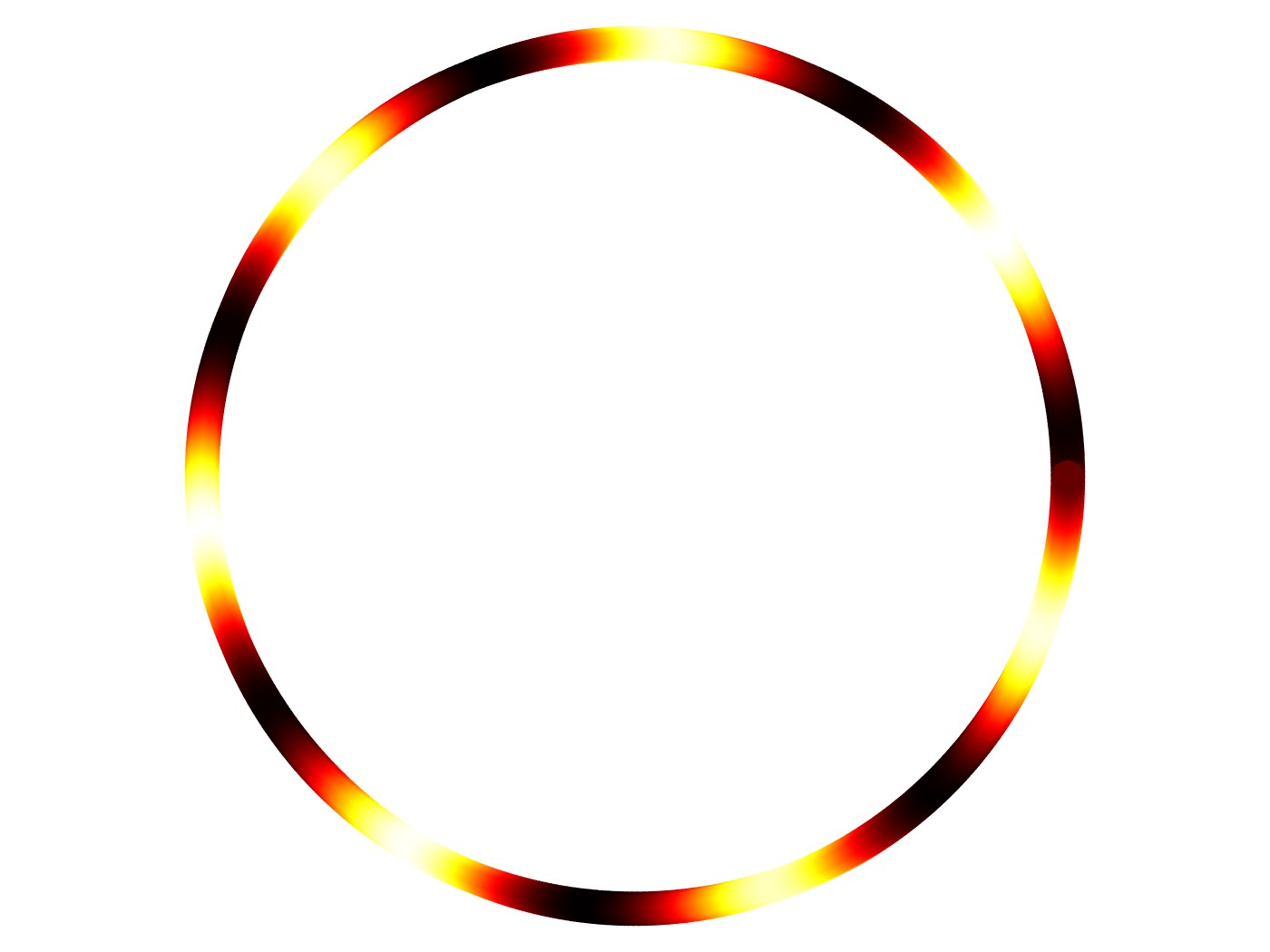}
    \put(10,-5){\scriptsize (c) Surface solution $w$}
\end{overpic}
\caption{Example~\ref{eg: CBSRD 2}: Spot pattern formation by solving a 2D coupled bulk-surface reaction-diffusion equation, using the meshfree time stepping method, but without using the greedy algorithm. The meshfree method uses $n = n_{\Omega} + n_{\cals} = 717 + 100$ data points and a time step of $\dt=0.005$. (a) Bulk and surface solutions together are shown in 3D, (b) pattern formed in the bulk solution $u$, and (c) pattern formed in the surface solution $w$.
   }\label{fig: CBS2D w/o g demo}
\end{figure}

Using Turing's pattern formation parameters in \cite{CBSCM} for spots pattern formation, as shown in Table~\ref{ta: parameters for CBS}, we solved the 2D coupled bulk-surface reaction-diffusion equations to study pattern formation, using a meshfree time stepping method. Firstly, without the greedy algorithm, Figure~\ref{fig: CBS2D w/o g demo} shows the meshfree solutions $u$ in the bulk and $w$ on the surface obtained by employing $n = n_{\Omega} + n_{\cals} = 717 + 100$ data points and a time step of $\dt=0.005$. Here $n_{\Omega}$ and $n_{\cals}$ denote the number of basis functions in $\Omega$ and on $\cals$. To ensure that the solutions converge to steady states, we solve the reaction-diffusion until a sufficiently long time, say, $T=1000$. 
Figure~\ref{fig: CBS2D w/o g demo}a shows the bulk and surface solutions together in 3D, revealing their complex spatial interactions.
Figures~\ref{fig: CBS2D w/o g demo}b and \ref{fig: CBS2D w/o g demo}c show bird's eye views of the bulk solution $u$ and surface solution $w$, respectively. The bulk solution exhibits a spotted pattern while the surface solution develops a wavelike pattern around the circular domain.

We implement the greedy algorithm to determine if the patterns can be generated with fewer basis functions. To ensure numerical convergence, we vary the time steps and the number of data points. Figure \ref{fig: CBS2D w and w/o g} displays overlaid bulk and surface solutions of $\triangle t = 0.01$, $\triangle t = 0.005$, and $n = 717 + 100$, $n = 1031 + 120$, with and without the greedy algorithm. Upon initial inspection, all four meshfree parameter sets, irrespective of whether the greedy algorithm was used, generate similar bulk patterns, characterized by 7 spots in the bulk and either 6 or 7 peaks on the surface. We also examine the interfacial patterns. In Figure \ref{fig: CBS2D w surf phase}, we plot all eight (non-unique) surface patterns $w$, shifted so that the first peak coincides. The numerical solutions for the surface exhibit less consistency than those in the bulk. 
Notably, the greedy surface solutions with small $n$, say $n = 717+100$ exhibit a distinct characteristic from the other solutions, displaying only a six-peak wave-like pattern. We plot these special solutions in red, while all other solutions are plotted in blue. Figure \ref{fig: CBS2D w and w/o g}e and \ref{fig: CBS2D w and w/o g}g correspond to these solutions showing that the bulk and surface solutions are coupled together. We recorded the CPU times used to select the trial space for all variables using the greedy algorithm in Table \ref{ta: CBS2D eps and SC}. To reduce random variability, we counted the average CPU time of executing the greedy algorithm ten times and represented it by $t_{CPU}$ (seconds).

\begin{table}
\centering
\caption{The number of selected basis functions and the stopping criteria corresponding to the solutions corresponding greedy cases in Figure \ref{fig: CBS2D w and w/o g}.}\label{ta: CBS2D eps and SC}
\resizebox{\columnwidth}{!}{%
\begin{tabular}{|c||c | c| c| c|c|c|c|c|c|c|}
\hline
\multirow{ 2}{*}{} &  \multicolumn{5}{c|}{$n = 717+100$} & \multicolumn{5}{c|}{$n = 1031+120$}\\
\cline{2-11}
 & $n_{\Omega}'(717)$   & $n_{\cals}'(100)$    & SC($\Omega$)   & SC($\cals$) & $t_{CPU}$ & $n_{\Omega}'(1031)$   & $n_{\cals}'(120)$    & SC($\Omega$)   & SC($\cals$) & $t_{CPU}$ \\ \hline
$\triangle t = 0.01$  &  341 & 23   & SC-2$'$  & SC-2$'$  & 7.27 &  186 & 39   & SC-2$'$  & SC-2$'$ & 5.74\\
$\triangle t = 0.005$   &  451 & 23   & SC-2$'$  & SC-2$'$ &5.90 &  465 & 24   & SC-2$'$  & SC-2$'$ & 5.58\\
\hline
\end{tabular}%
}
\end{table}

In Table \ref{ta: CBS2D eps and SC}, we present the stopping criteria and the number of basis functions selected for the greedy solutions corresponding to Figure \ref{fig: CBS2D w and w/o g} and Figure \ref{fig: CBS2D w surf phase}, with varying values of $n_{\Omega}$ and $n_{\cals}$. For both functions $u$ and $w$, fewer basis functions are selected to produce consistent patterns. Using smaller $\triangle t$ results in more basis functions being selected. We remark that all kernels use $\epsilon=6$ and run to $T=1000$. Figures~\ref{fig: CBS2D w/o g demo}--\ref{fig: CBS2D w and w/o g} is kept at the same range for the sake of fair comparison of spot size.

\begin{figure}
  \centering
  \begin{overpic}[width=0.24\textwidth,trim=150 20 150 0, clip=true,tics=10]{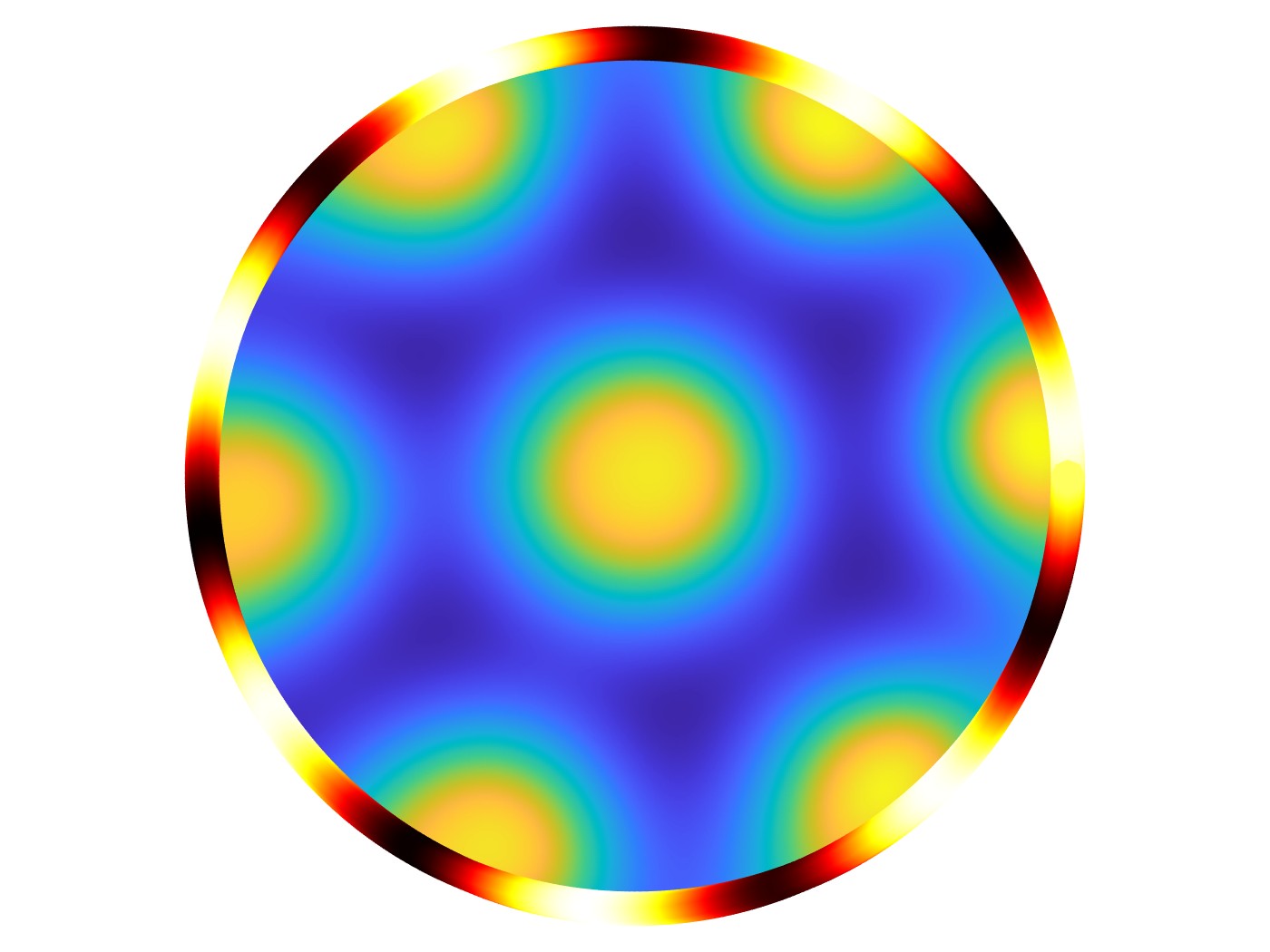}
    \put(-20,25){\rotatebox{90}{\textbf{Without}}}
    \put(-10,28){\rotatebox{90}{\textbf{Greedy}}}
    \put(20,95){$n=717+100$ }
     \put(80,110){$\dt=0.01$ }
    \put(17,-8){\scriptsize (a) 7 spots, 7 peaks}
\end{overpic}
\begin{overpic}[width=0.24\textwidth,trim=150 20 150 0, clip=true,tics=10]{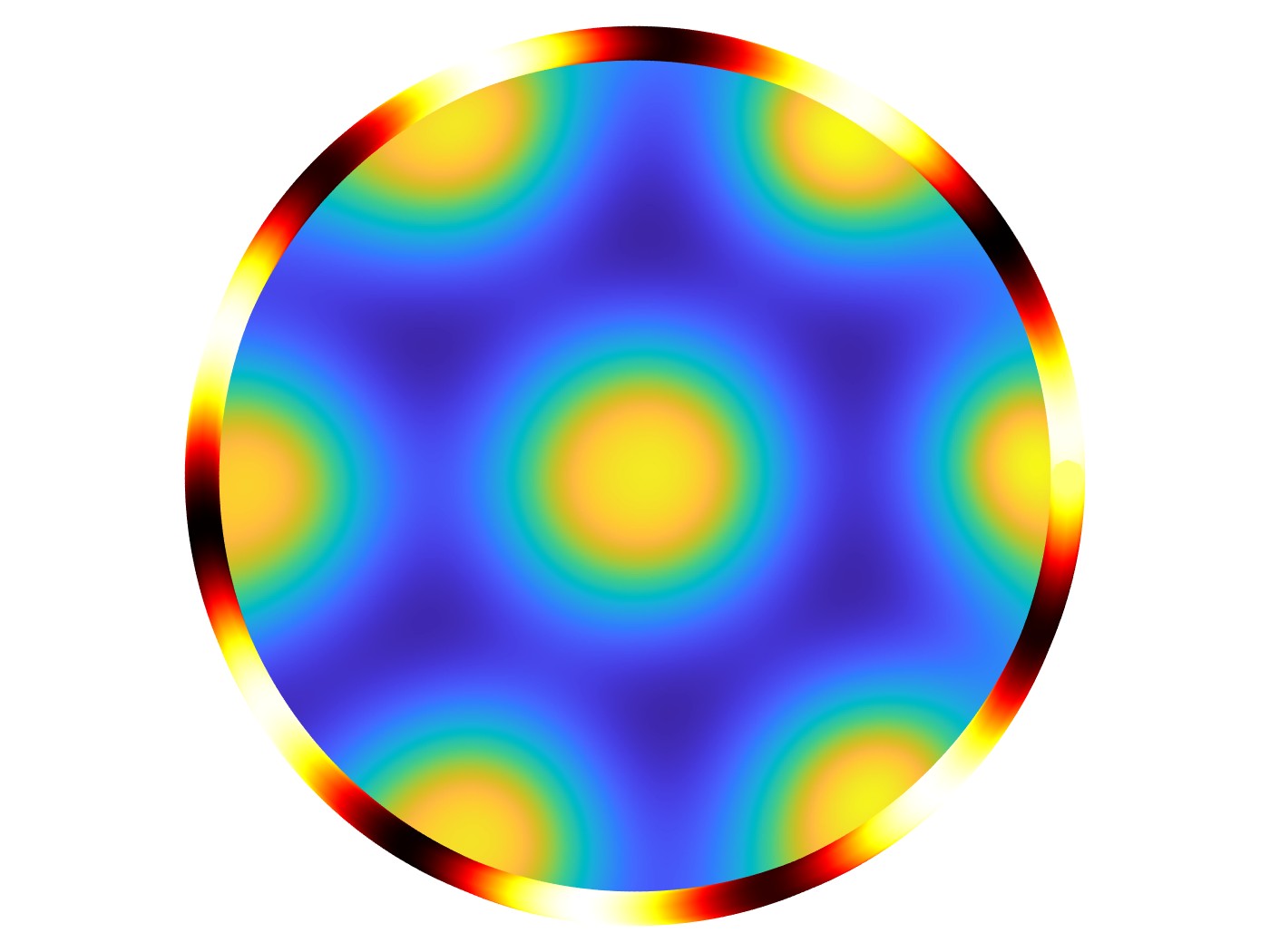}
    \put(20,95){$n=1031+120$ }
    \put(17,-8){\scriptsize (b) 7 spots, 7 peaks}
\end{overpic}
\begin{overpic}[width=0.24\textwidth,trim=150 20 150 0, clip=true,tics=10]{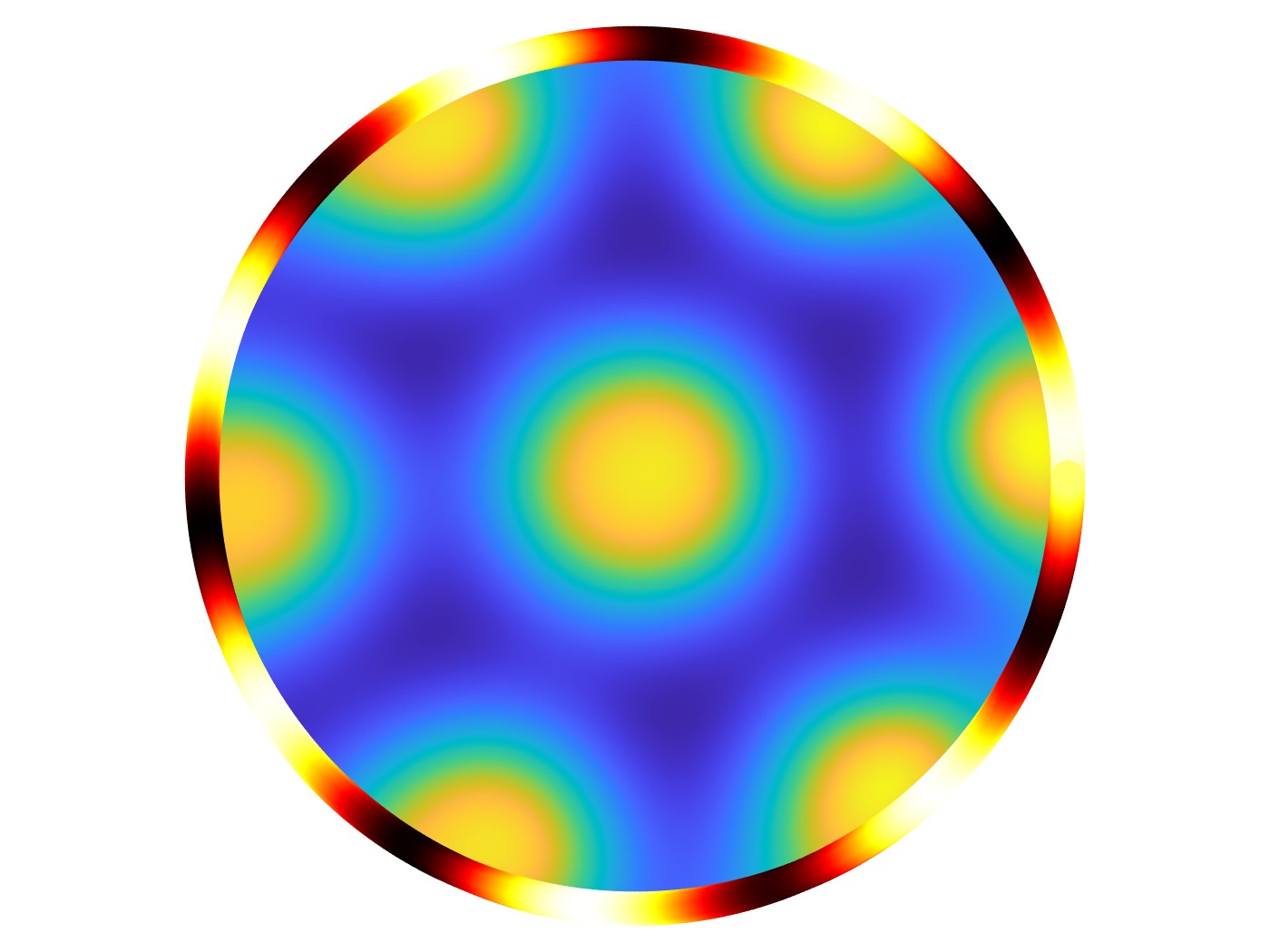}
    \put(20,95){$n=717+100$ }
     \put(78,110){$\dt=0.005$ }
    \put(17,-8){\scriptsize (c) 7 spots, 7 peaks}
\end{overpic}
\begin{overpic}[width=0.24\textwidth,trim=150 20 150 0, clip=true,tics=10]  {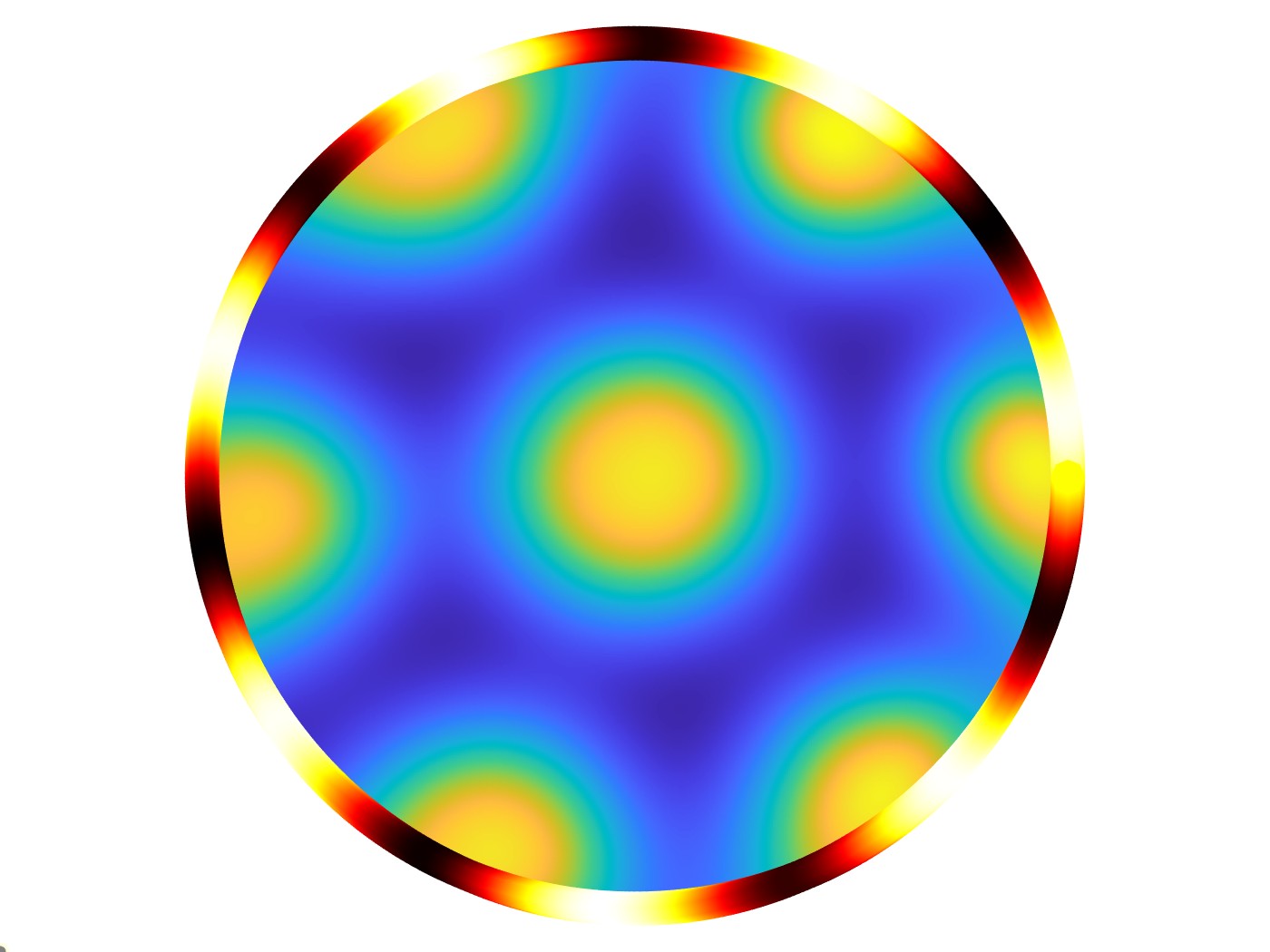}
    \put(20,95){$n=1031+120$} 
    \put(17,-8){\scriptsize (d) 7 spots, 7 peaks}
\end{overpic}
\\
\medskip
\medskip
\begin{overpic}[width=0.24\textwidth,trim=150 20 150 0, clip=true,tics=10]{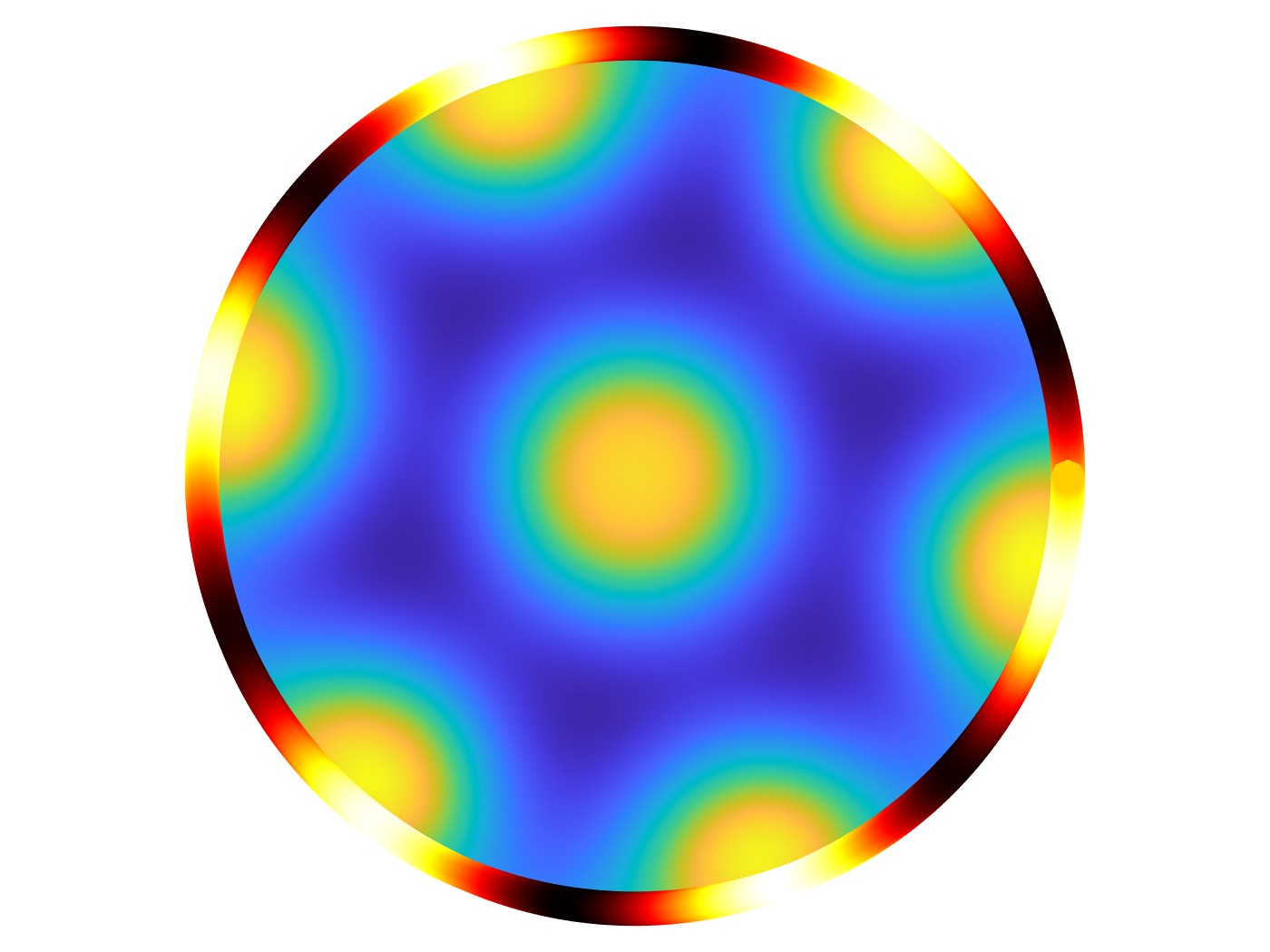}
     \put(-20,33){\rotatebox{90}{\textbf{With}}}
    \put(-10,28){\rotatebox{90}{\textbf{Greedy}}}
    \put(17,-8){\scriptsize (e) 7 spots, 6 peaks}
\end{overpic}
\begin{overpic}[width=0.24\textwidth,trim=150 20 150 0, clip=true,tics=10]{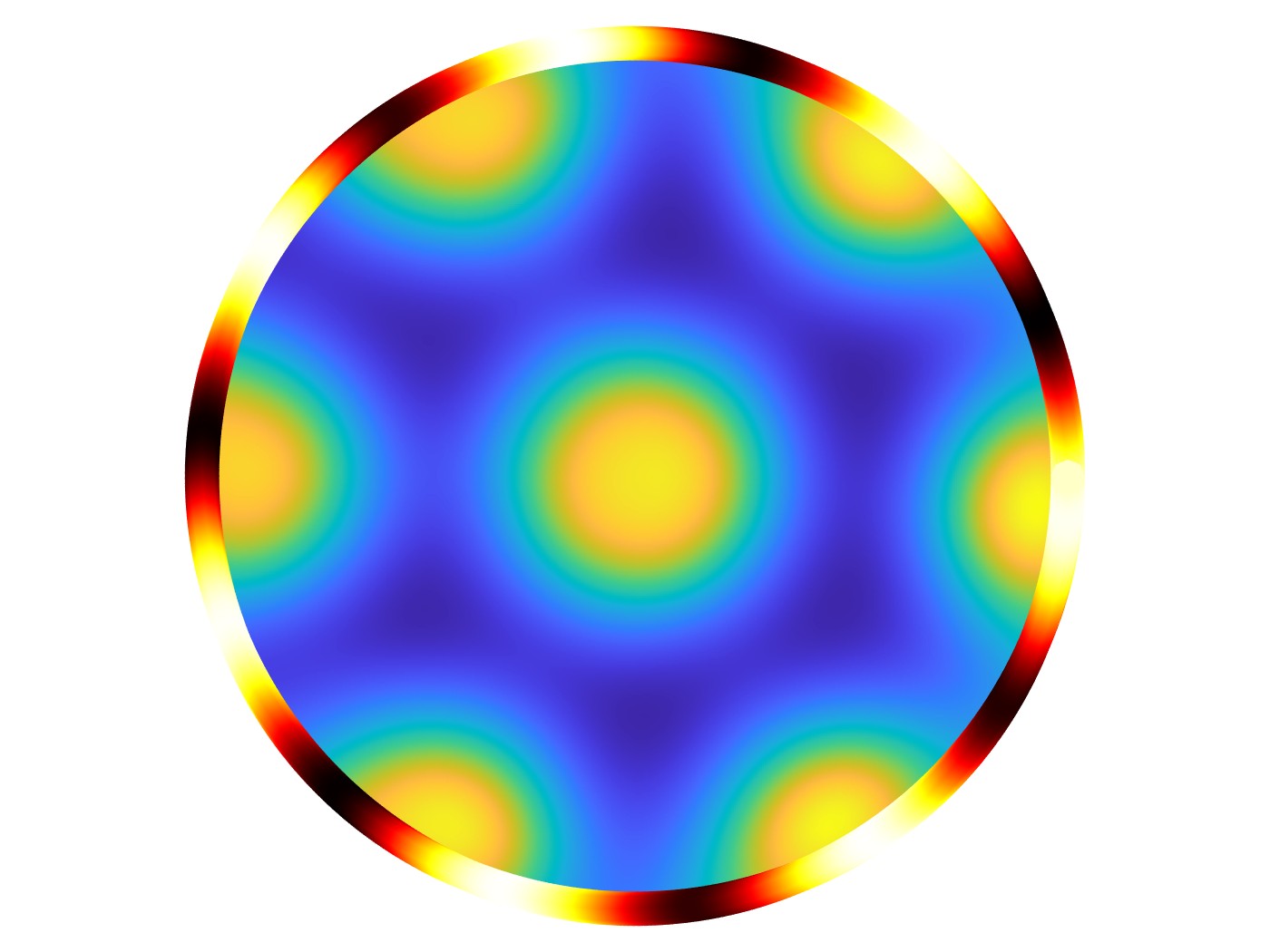}
    \put(17,-8){\scriptsize (f) 7 spots, 7 peaks}
\end{overpic}
\begin{overpic}[width=0.24\textwidth,trim=150 20 150 0, clip=true,tics=10]{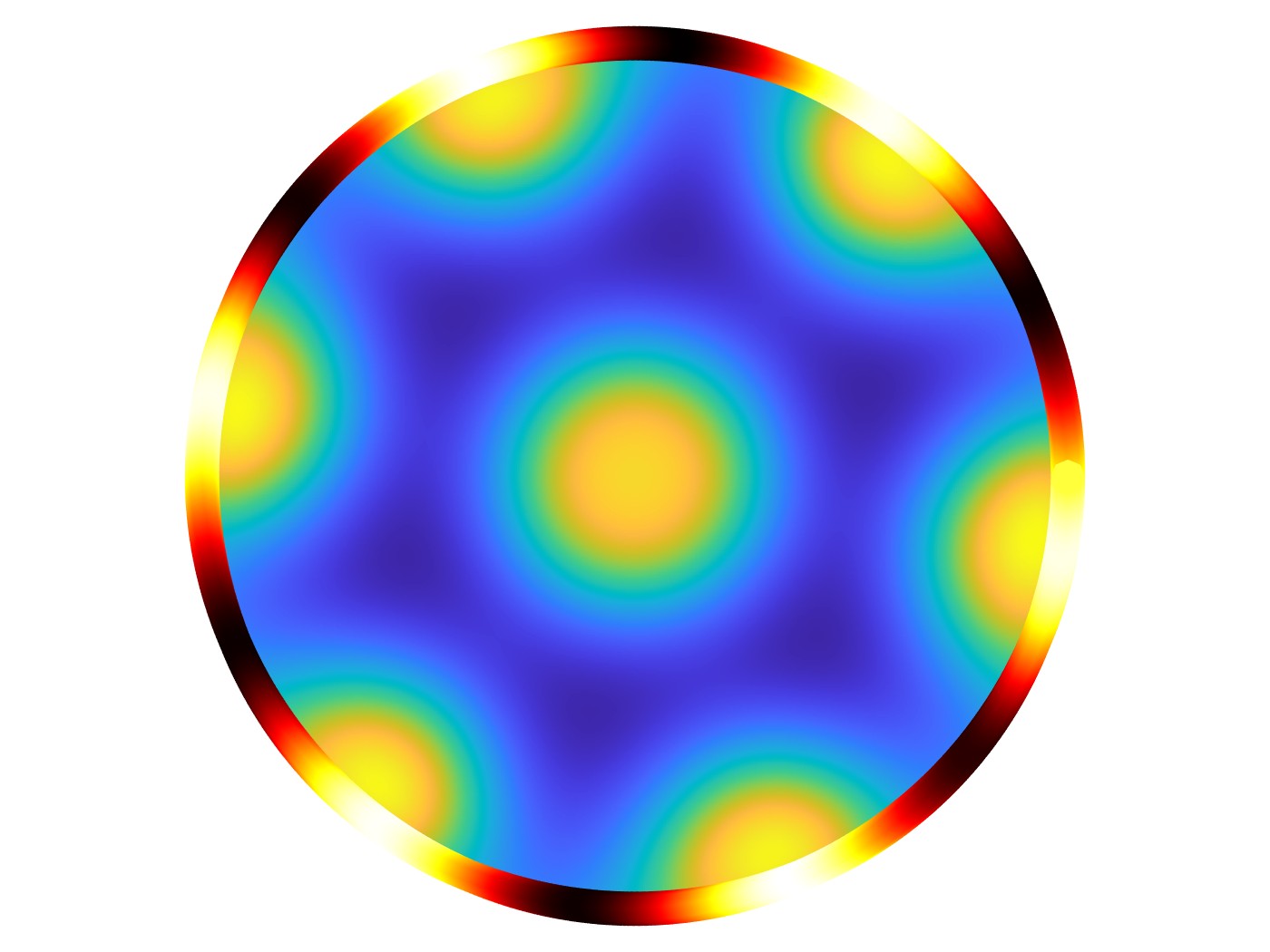}
    \put(17,-8){\scriptsize (g) 7 spots, 6 peaks}
\end{overpic}
\begin{overpic}[width=0.24\textwidth,trim=150 20 150 0, clip=true,tics=10]  {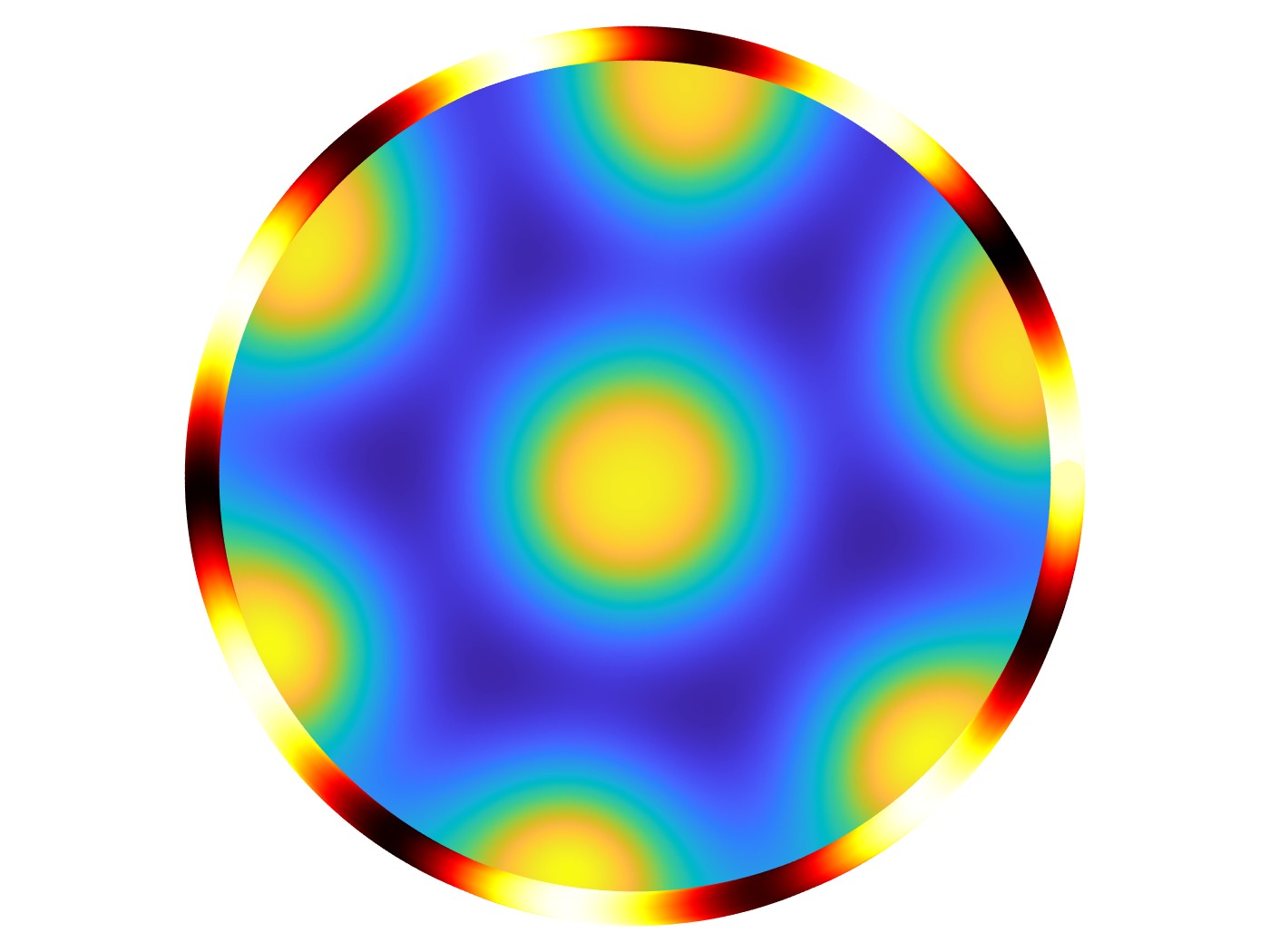}
    \put(17,-8){\scriptsize (h) 7 spots, 7 peaks}
\end{overpic}
\caption{Example~\ref{eg: CBSRD 2}: Bulk and surface solutions obtained by solving the 2D coupled bulk-surface reaction-diffusion equation with and without the greedy algorithm. The meshfree method uses $n=717+100$, $n = 1031+120$ data points and $\dt = 0.005$, $\dt=0.01$ time steps. The number of bases used and stopping criteria in the greedy cases are in Table \ref{ta: CBS2D eps and SC}.
   }\label{fig: CBS2D w and w/o g}
\end{figure}

\begin{figure}
  \centering
\begin{overpic}[width=0.95\textwidth,trim=20 0 20 0, clip=true,tics=10]{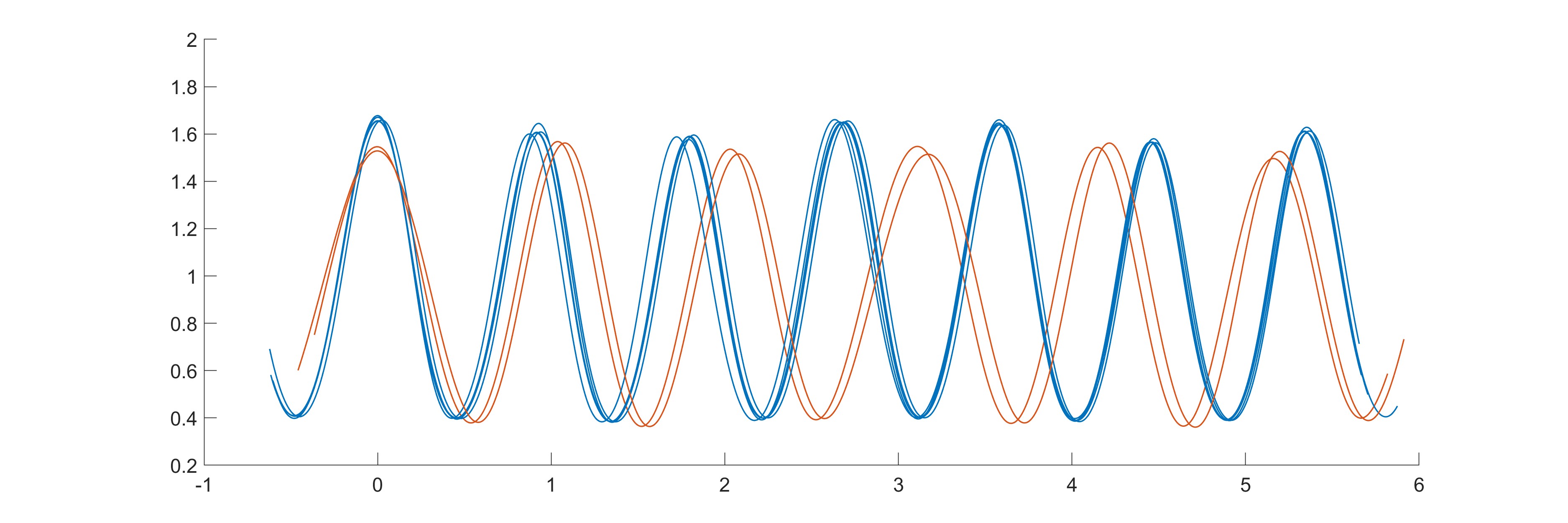}
    \put(7,12){\rotatebox{90}{\scriptsize$w(\theta+\triangle\theta)$}}
    \put(55,0){\scriptsize$\theta$ }
    \put(45,30){\linethickness{0.2mm}\color{frenchblue}\line(1,0){4.8}}%
    \put(66,29.7){\scriptsize Greedy on $n=717+100$}
    \put(60,30){\linethickness{0.2mm}\color{sinopia}\line(1,0){4.8}}%
    \put(51,29.7){\scriptsize all other}
\end{overpic}
\caption{Example~\ref{eg: CBSRD 2}: The surface solutions from Figure~\ref{fig: CBS2D w and w/o g} are individually shifted (with different $\triangle t$) to align the location of the first peak and plotted against $\theta$.
}\label{fig: CBS2D w surf phase}
\end{figure}
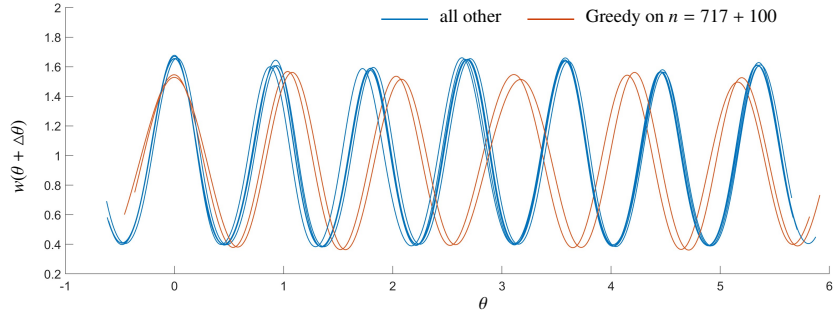

We present additional experiments that demonstrate how the greedy algorithm can rectify unstable numerical results that arise due to poor numerical settings. Figure \ref{fig: CBS2D w g} shows the solutions obtained with the greedy algorithm under different parameter settings. Specifically, building upon the scenarios depicted in Figure \ref{fig: CBS2D w and w/o g}, we take more data points wherein $n$ was set to $1612+150$ with $\triangle t$ being $0.01$ and $0.005$ in Figure \ref{fig: CBS2D w g}a-\ref{fig: CBS2D w g}b. Additionally, in Figure \ref{fig: CBS2D w g}c-\ref{fig: CBS2D w g}d, we varied the diffusion coefficients from $D_{v} = D_{s} = 2$ to $D_{v} = D_{s} = 1$ with $\triangle t = 0.005$, $n = 717+100$, and $n = 1031+120$. While there remains some uncertainty regarding whether the greedy algorithm accurately captures the distribution of the patterns, particularly in cases where $D_{v} = D_{s} = 1$, it is worth noting that the appearance of the patterns is indicative of the success of the greedy algorithm in comparison to the blow-up solution without the greedy algorithm under these settings.

\begin{figure}
  \centering
\begin{overpic}[width=0.24\textwidth,trim=150 20 150 0, clip=true,tics=10]  {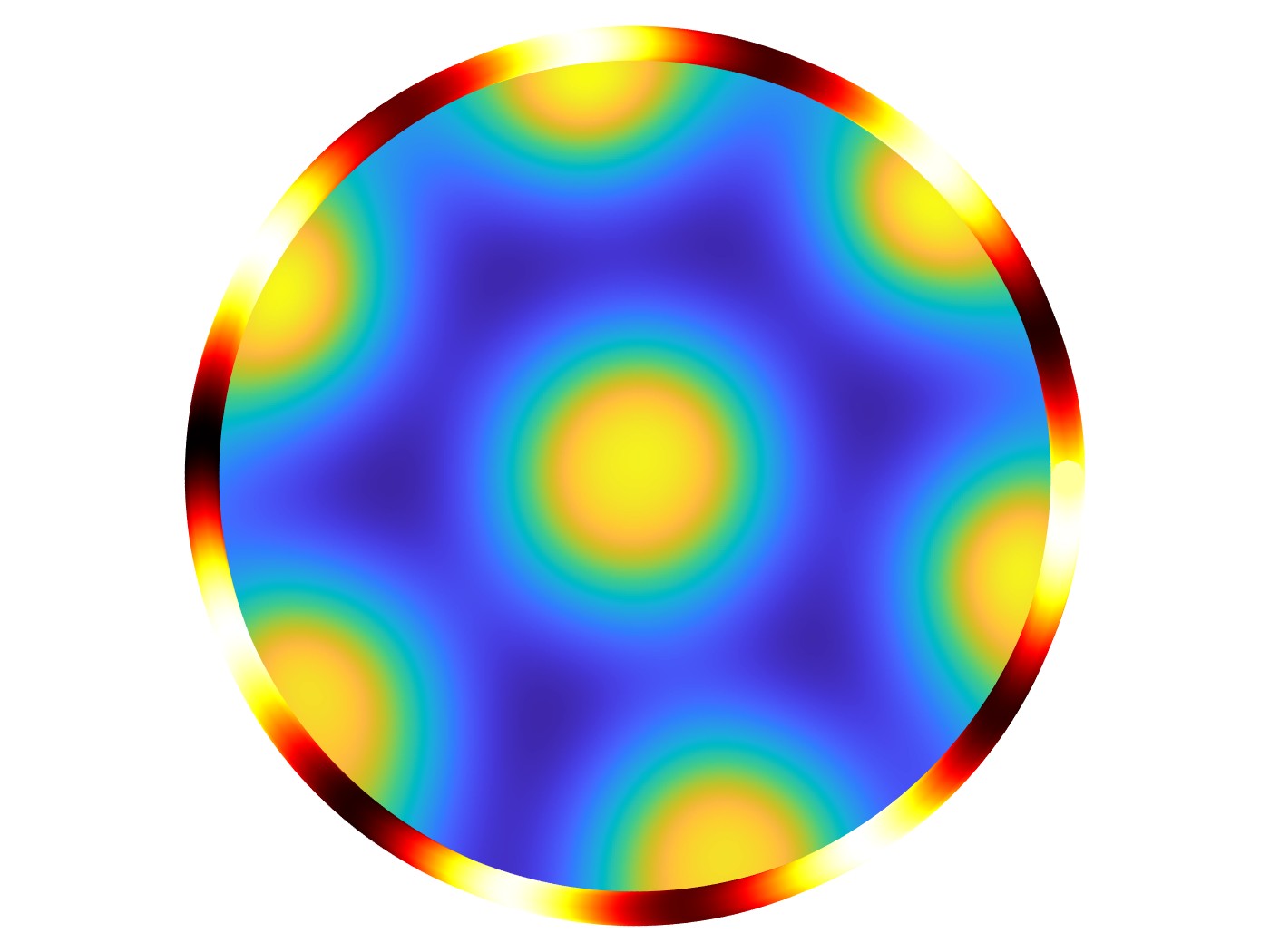}
     \put(-20,33){\rotatebox{90}{\textbf{With}}}
    \put(-10,28){\rotatebox{90}{\textbf{Greedy}}}
     \put(78,100){$D_{v} = D_{s} = 2$}
     \put(30, -8){\scriptsize (a) $\dt = 0.01$,}
     \put(30,-16){\scriptsize $n = 1612 + 150$} 
\end{overpic}
\begin{overpic}[width=0.24\textwidth,trim=150 20 150 0, clip=true,tics=10]{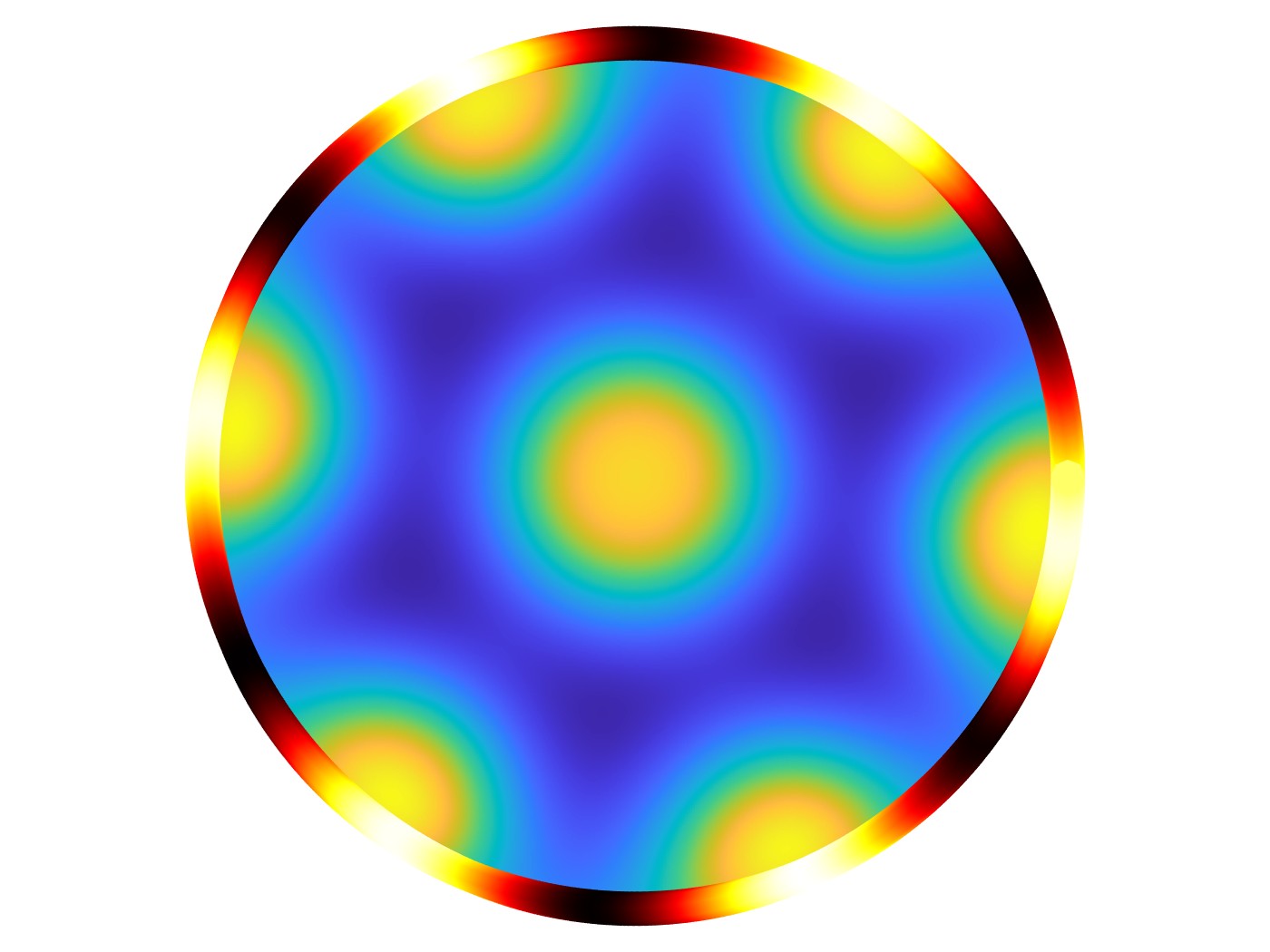}
    \put(30, -8){\scriptsize (b) $\dt = 0.005$,}
     \put(30,-16){\scriptsize $n = 1612 + 150$} 
\end{overpic}
\begin{overpic}[width=0.24\textwidth,trim=150 20 150 0, clip=true,tics=10]{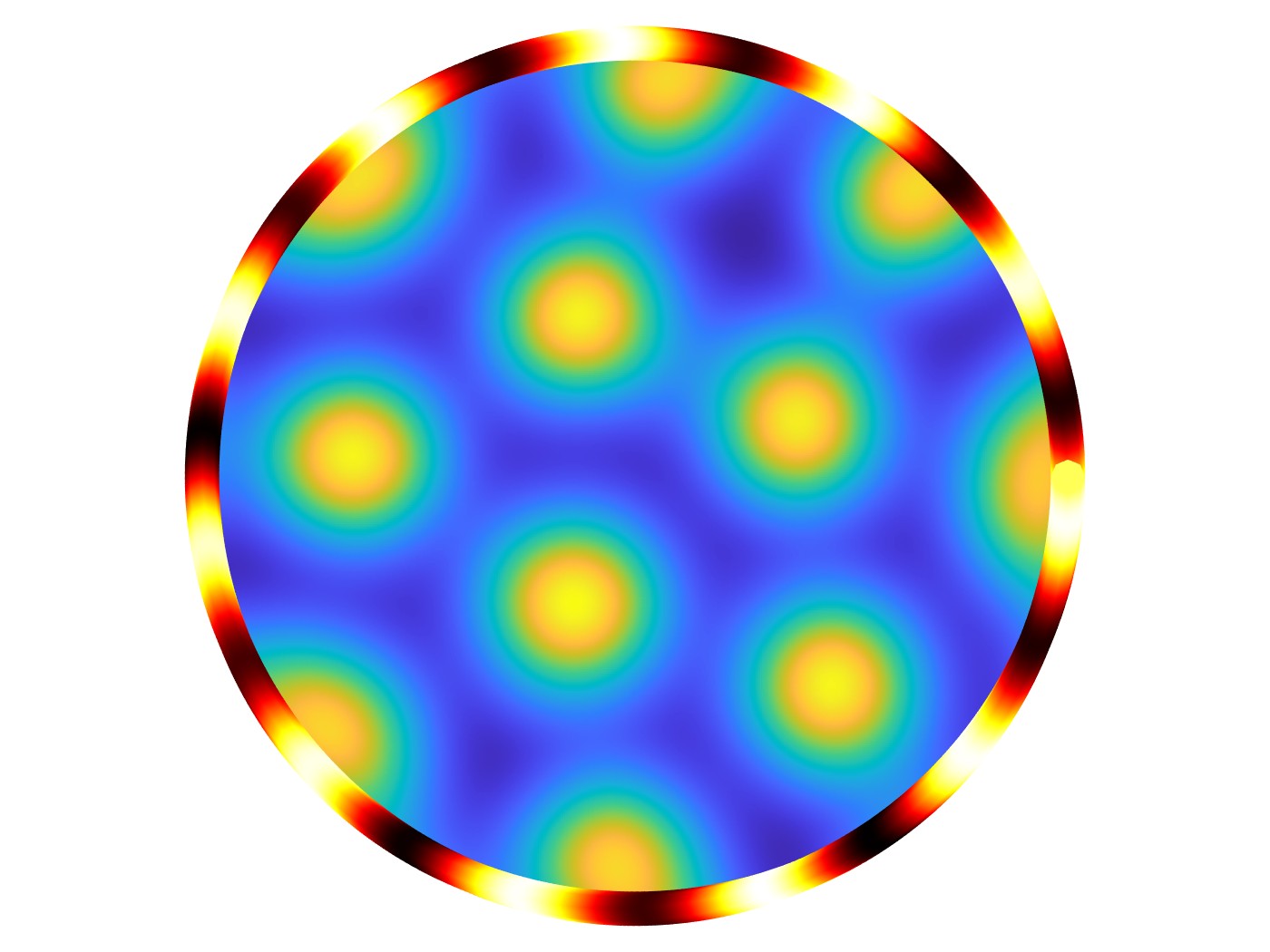}
     \put(80,100){$D_{v} = D_{s} = 1$ }
     \put(30, -8){\scriptsize (c) $\dt = 0.005$,}
     \put(35,-16){\scriptsize $n = 717+100$}
\end{overpic}
\begin{overpic}[width=0.24\textwidth,trim=150 20 150 0, clip=true,tics=10]{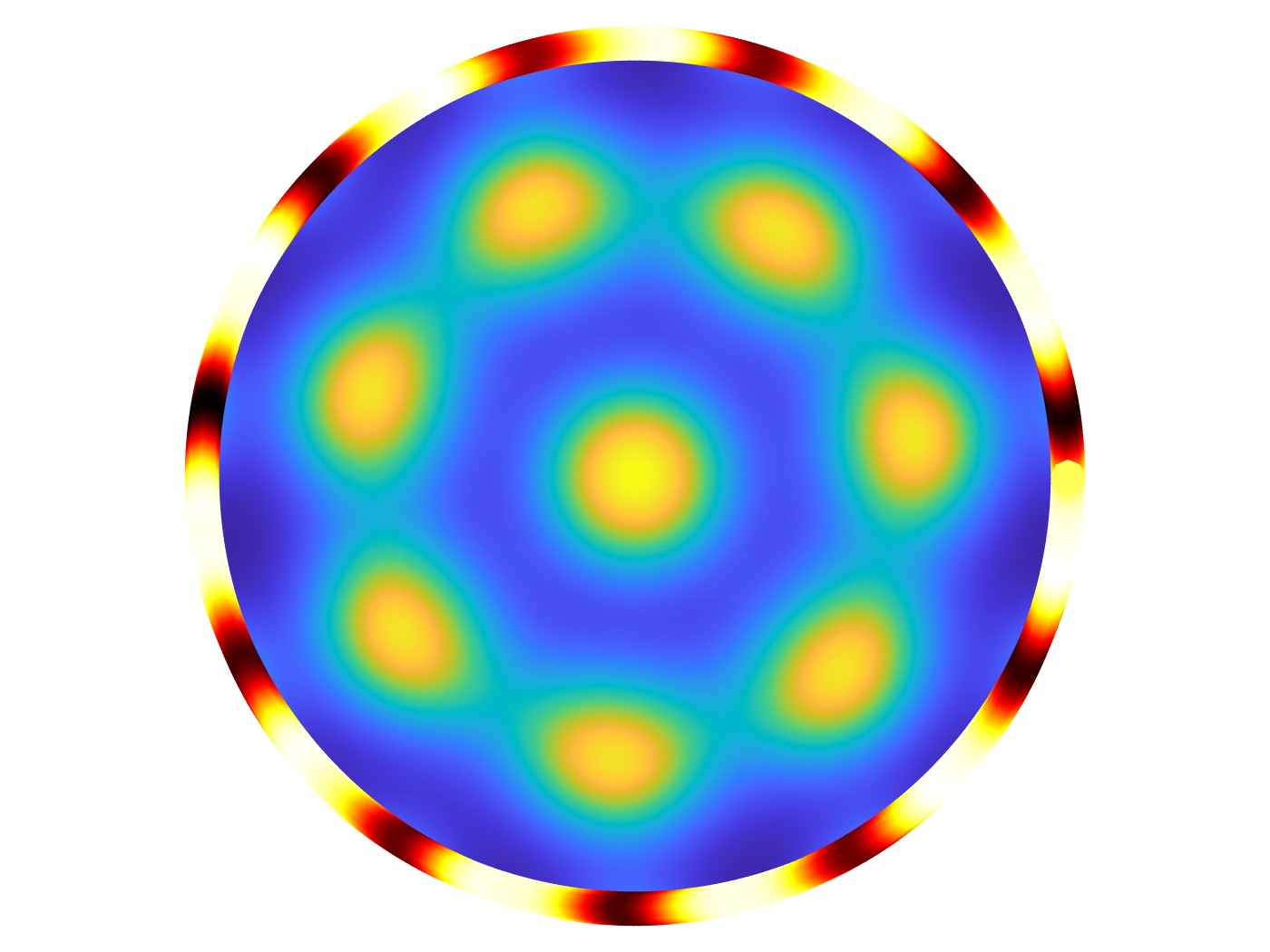}
     \put(30, -8){\scriptsize (d) $\dt = 0.005$,}
     \put(32,-16){\scriptsize $n = 1031+120$} 
\end{overpic}
\medskip
\caption{Example~\ref{eg: CBSRD 2}: Corresponding to the greedy cases in Figure \ref{fig:  CBS2D w and w/o g}, bulk and surface solutions by solving 2D coupled bulk-surface reaction-diffusion equation using $D_{v} = D_{s} = 2$ and $D_{v} = D_{s} = 1$.
   }\label{fig: CBS2D w g}
   \medskip
\end{figure}
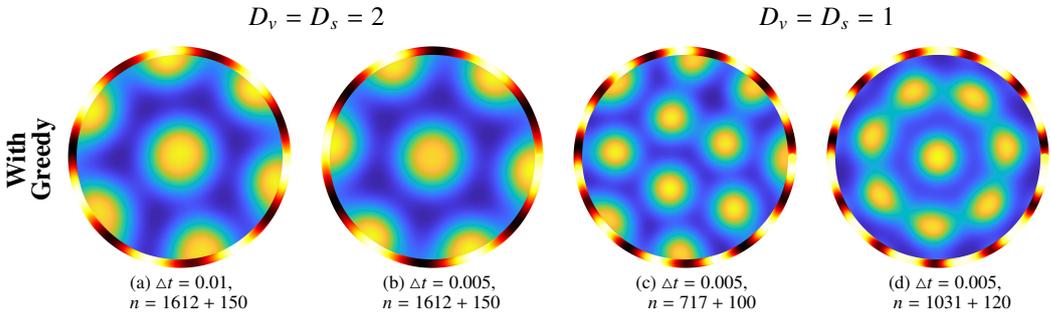

Before tackling 3D problems, we also attempt to solve the numerically more challenging stripe pattern formation problem. We do this by simply changing the parameters for Turing's pattern formation, as shown in Table~\ref{ta: parameters for CBS}. Since generating stripe patterns typically requires a higher level of spatial and temporal resolution than spot patterns due to their more complex and dynamic nature, we take a smaller time step as $\triangle t = 0.001$ and more data points as $n = 2869+200$ and $n = 4477+250$ in this case. The run time here is still taken as $T=1000$. Initially, we attempted to solve for stripe patterns under these settings without utilizing the greedy algorithm. However, as we had anticipated, the solution experienced a blow-up in the first few iterations. Then we employ the greedy algorithm, which selects fewer bases to obtain reasonable stripe patterns. Figure~\ref{fig: CBS2D w g stripe} showcases both solutions in the bulk and on the surface, which successfully capture the formation of stripe patterns in the bulk domain and exhibit very high-frequency wavelike surface solutions. The amplitudes of surface solutions appear uniform and consistent. In both patterns, the bulk solution at the boundary and the surface solution demonstrate high-frequency oscillations and are perfectly synchronized\footnote{Similar to the peak-to-trough matching seen in cases of 7 spots in the bulk in Figure~ \ref{fig: CBS2D w and w/o g}, there is a similar matching in these stripe patterns.}

\begin{figure}
  \centering
\begin{overpic}[width=0.23\textwidth,trim=150 20 150 20, clip=true,tics=10]
{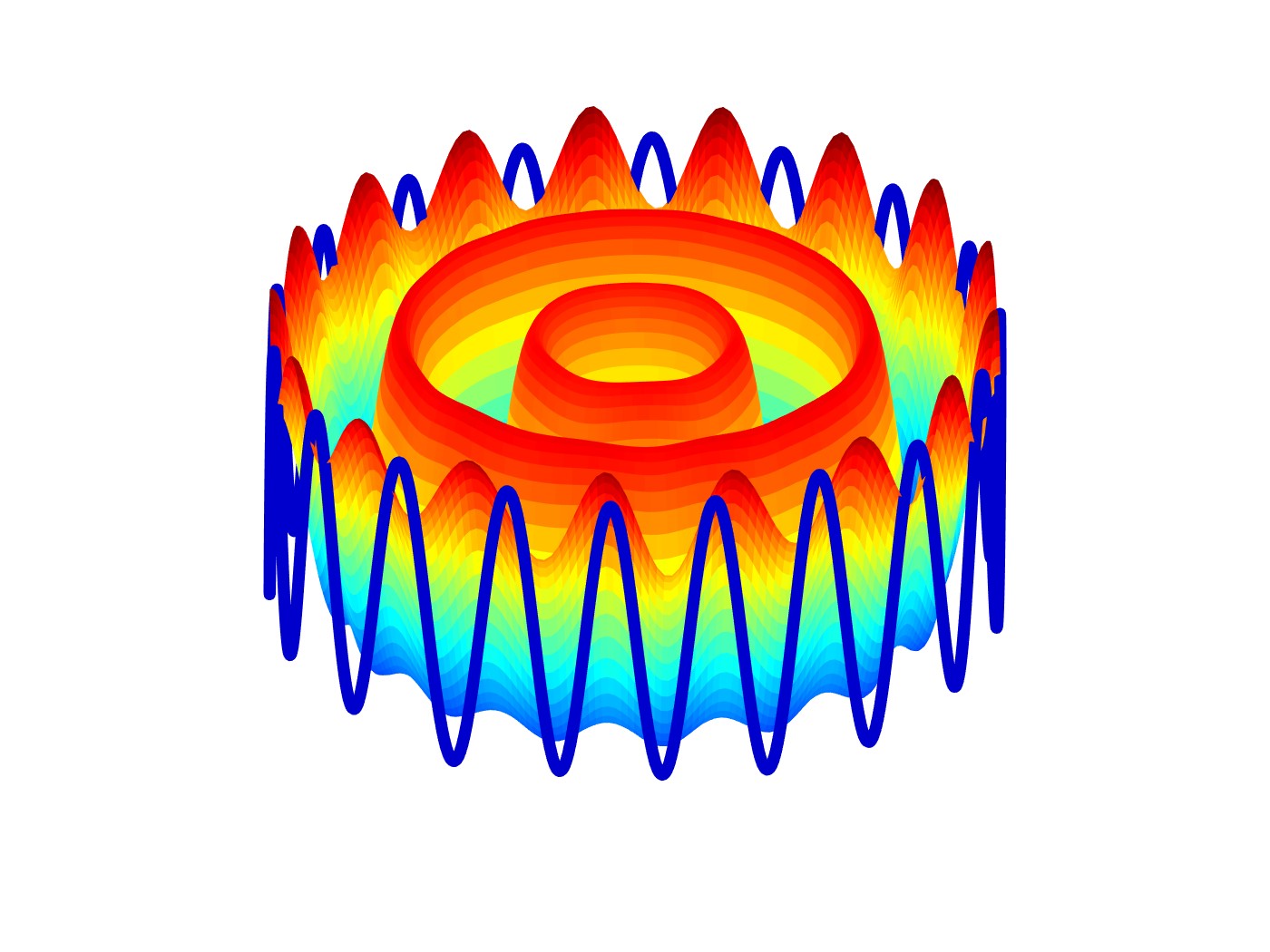}
    \put(65,95){$n=2869+200$} 
    \put(-20,33){\rotatebox{90}{\textbf{With}}}
    \put(-10,28){\rotatebox{90}{\textbf{Greedy}}}
\end{overpic}
\begin{overpic}[width=0.23\textwidth,trim=150 20 150 20, clip=true,tics=10]{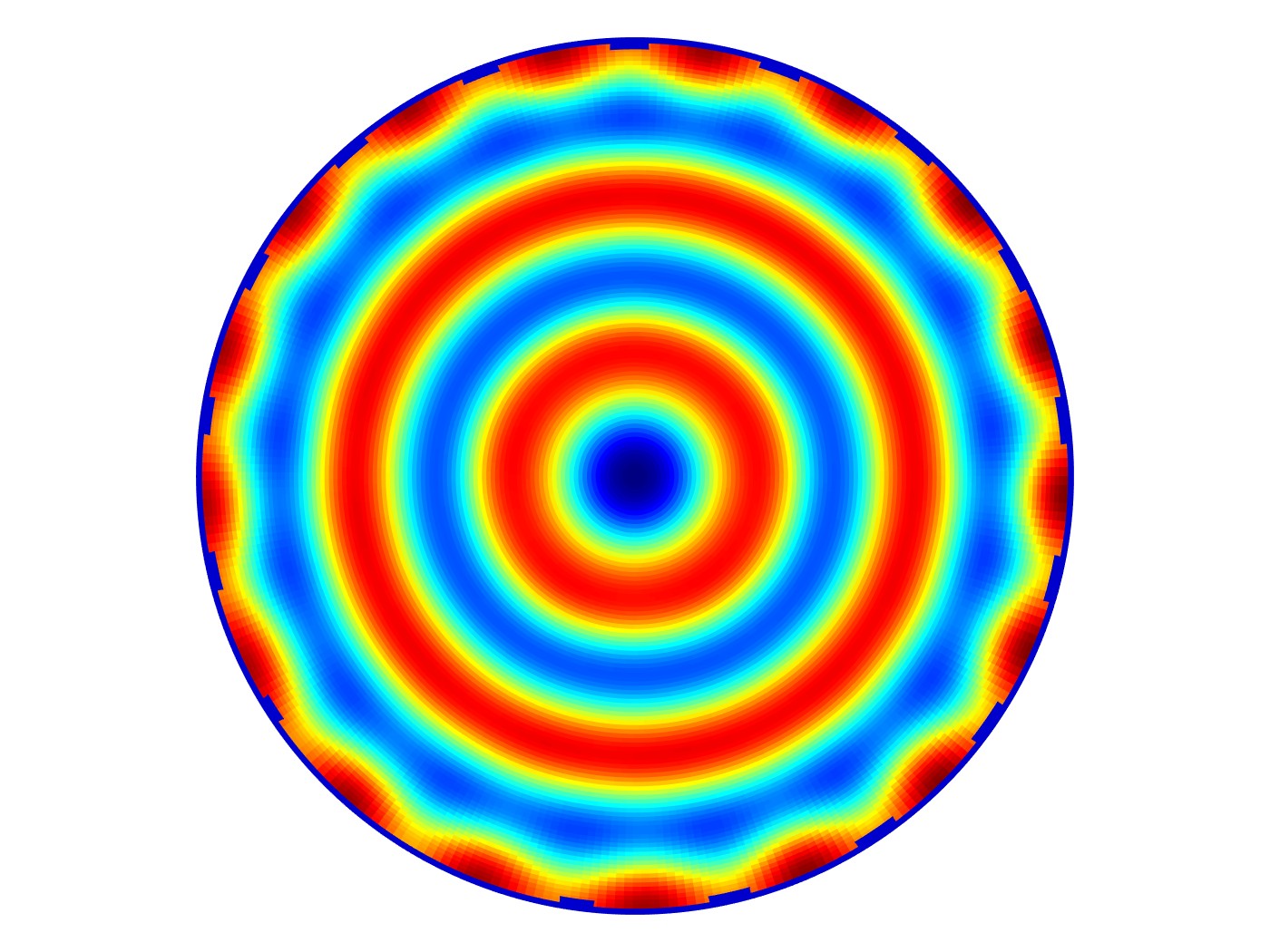}
\put(-23,10){\scriptsize 18 cycles}
\end{overpic}
\begin{overpic}[width=0.23\textwidth,trim=150 20 150 20, clip=true,tics=10]{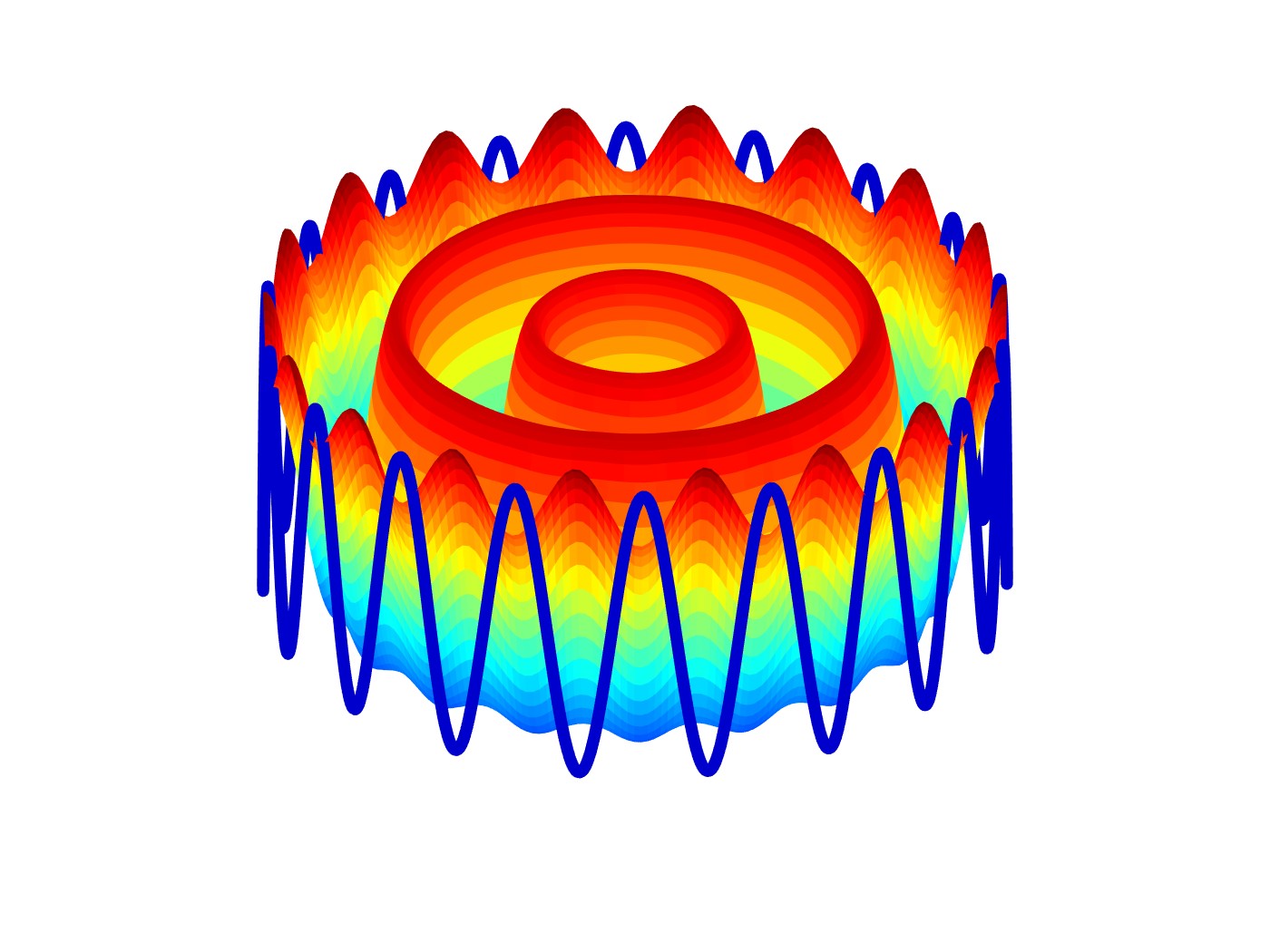}
    \put(65,95){$n=4477+250$} 
\end{overpic}
\begin{overpic}[width=0.23\textwidth,trim=150 20 150 20, clip=true,tics=10]{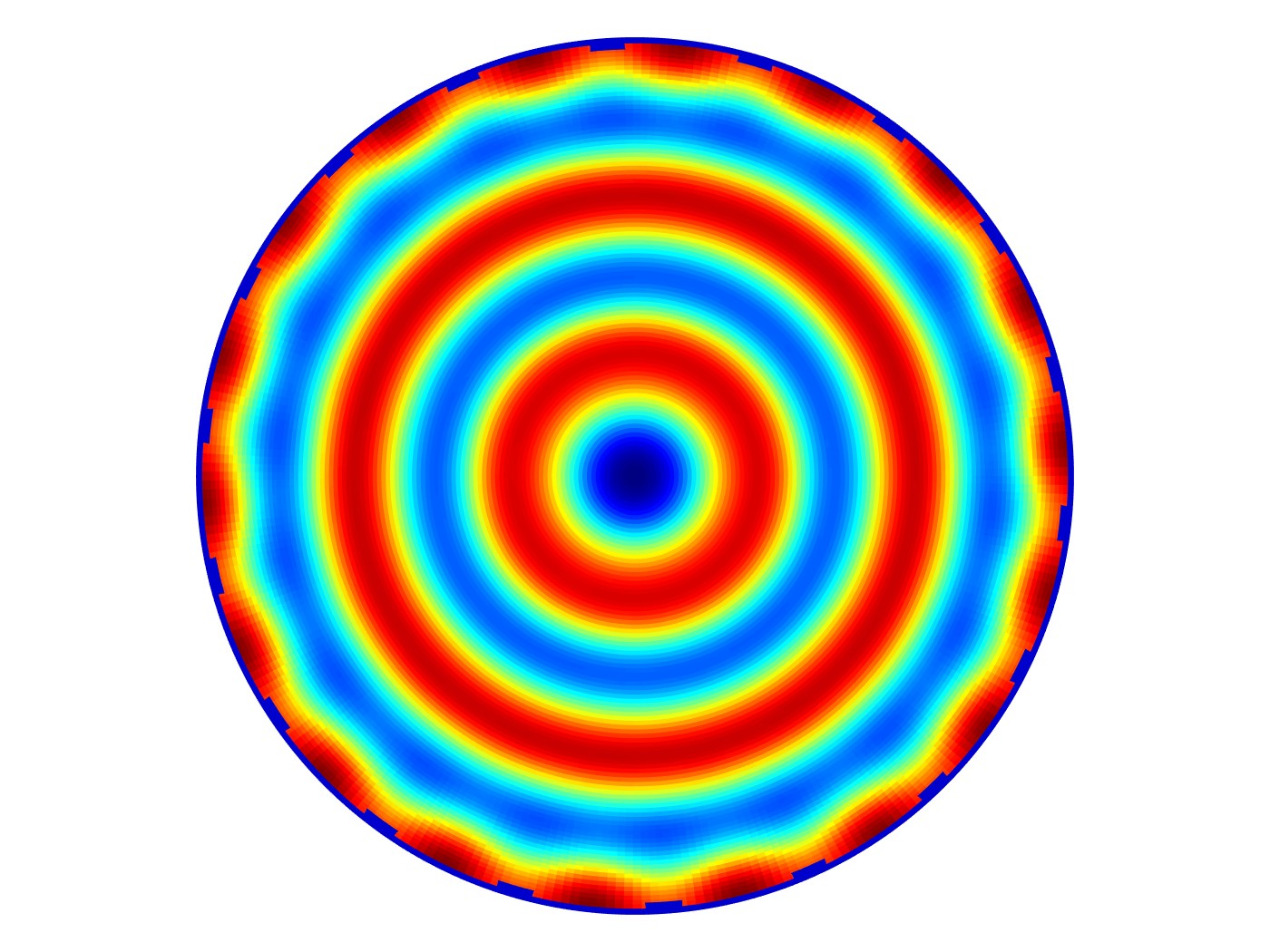}
    \put(-23,10){\scriptsize 18 cycles}
\end{overpic}
\caption{Example~\ref{eg: CBSRD 2}: Stripe pattern formation by solving a 2D coupled bulk-surface reaction-diffusion equation, using the meshfree time stepping method and the greedy algorithm. The meshfree method uses $n=2869+200$ and $4477+250$ data points and a time step of $\dt=0.001$.
For each test case, the bulk and surface patterns formed are shown together in 3D view and bird's eye view.
   }\label{fig: CBS2D w g stripe}
\end{figure}
\end{example}

\begin{example}\label{eg: CBSRD 3}\textbf{(Coupled Bulk-Surface Pattern formations in 3D)}

Moving forward with our simulation efforts, we now explore 3D coupled bulk-surface pattern formation.
Our aim is to test whether the greedy algorithm can produce patterns with relatively fewer basis functions in the cases when the solution is less oscillatory and smoother, i.e., spot formation. Compared to the 2D domain, generating patterns in 3D requires more accurate discretization due to the higher dimensionality. To simplify our approach, we focus on larger spot patterns in 3D domains instead of small patterns or slim stripes as seen in 2D cases. 

We implement the greedy algorithm approach as before without modification. Our initial experiments involve a torus with a simple geometrical structure, and we provide sufficient numerical solutions to verify the robustness of our approach. Subsequently, we will extend our experimentation to other 3D shapes, namely Dupin's cyclide and ellipsoid, to further validate our findings.

\begin{observation}[The greedy algorithm recovers from unsuccessful runs]

As in the 2D demonstration, we obtain bulk and surface solutions by solving the three-dimensional coupled bulk-surface reaction-diffusion equation was achieved by applying the meshfree time stepping method both with and without the greedy algorithm in Figure \ref{fig: CBS3D torus eps=1}. We still use the MS $\Phi_{6}$ kernel in 3D experiments. For the meshfree method, the value of $n$ was set to $2644+1430$, $\epsilon$ was fixed at 1 and the time step $\dt$ was allocated as 0.005 until the final time $T=1000$. We observe that in the absence of the greedy algorithm, faint and irregular patterns materialized, whereas the employment of the algorithm produced clearer, symmetric and well-separated spot configurations.

\end{observation}

\begin{figure}
  \centering
\begin{overpic}[width=0.24\textwidth,trim=130 150 130 150, clip=true,tics=10]{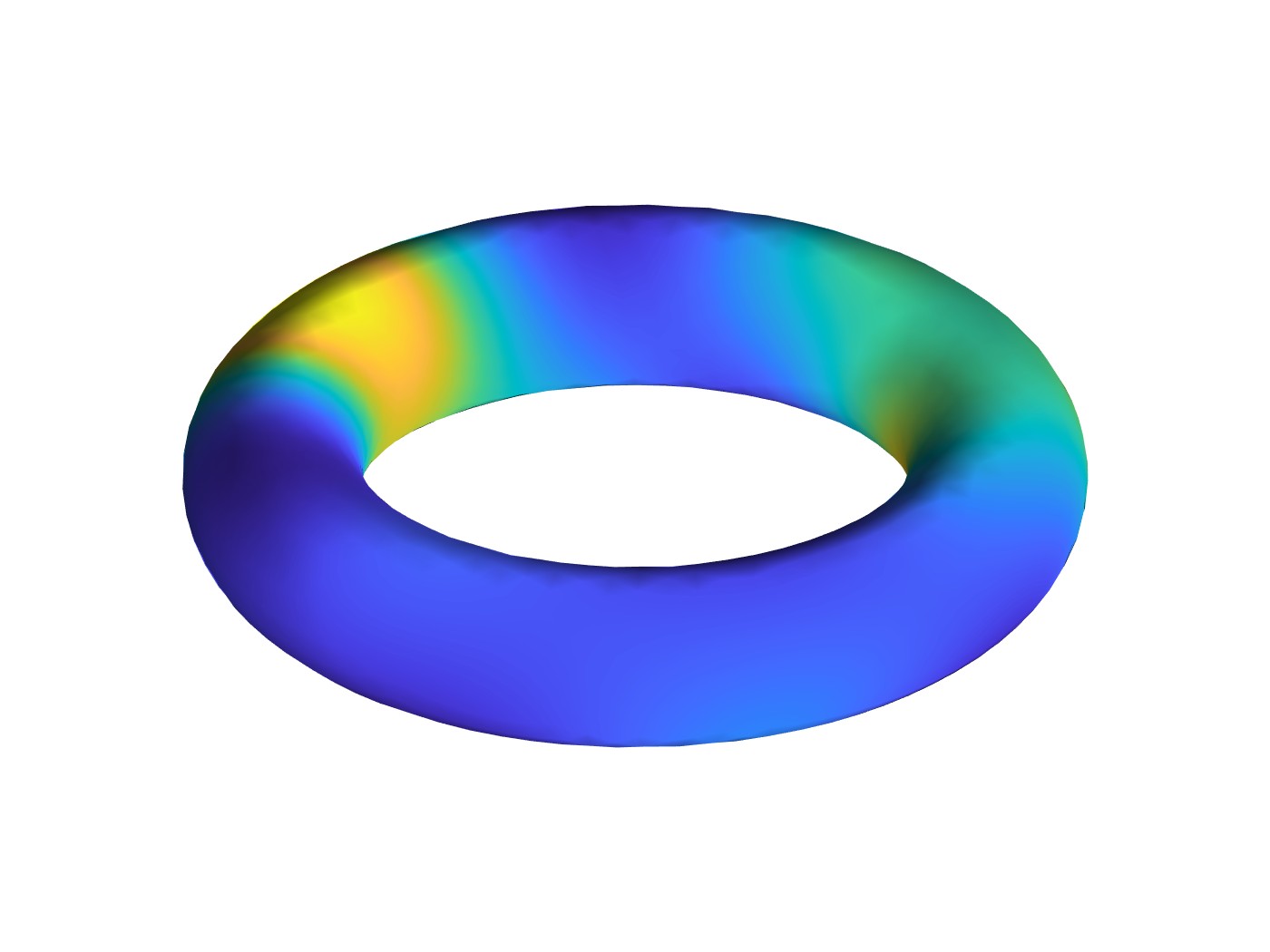}
    \put(45,70){\textbf{$u$}}
    \put(-20,15){\rotatebox{90}{\textbf{Without}}}
    \put(-10,18){\rotatebox{90}{\textbf{Greedy}}}
\end{overpic}
\begin{overpic}[width=0.24\textwidth,trim=130 150 130 150, clip=true,tics=10]{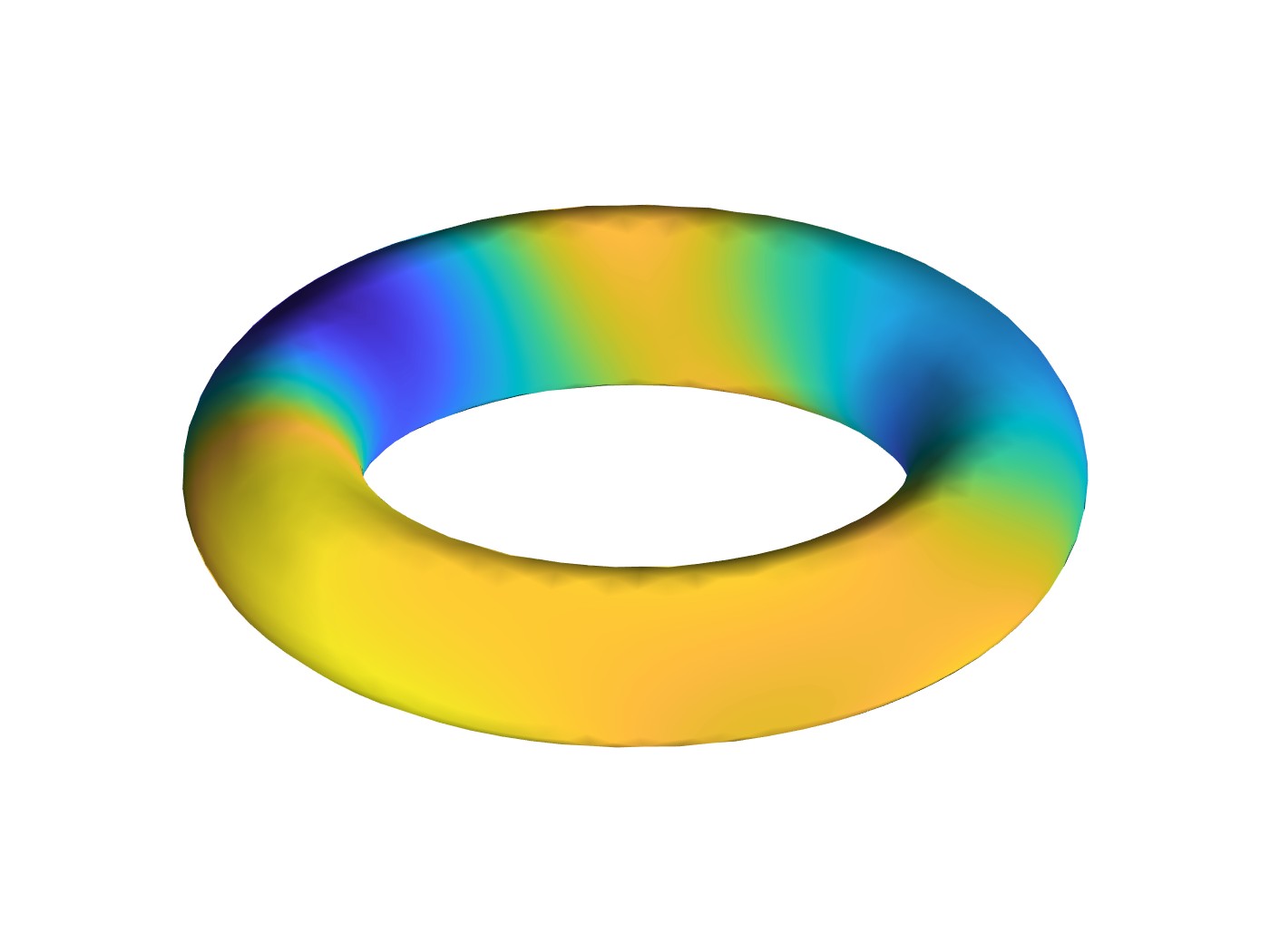}
    \put(45,70){\textbf{$v$}}
\end{overpic}
\begin{overpic}[width=0.24\textwidth,trim=130 150 130 150, clip=true,tics=10]{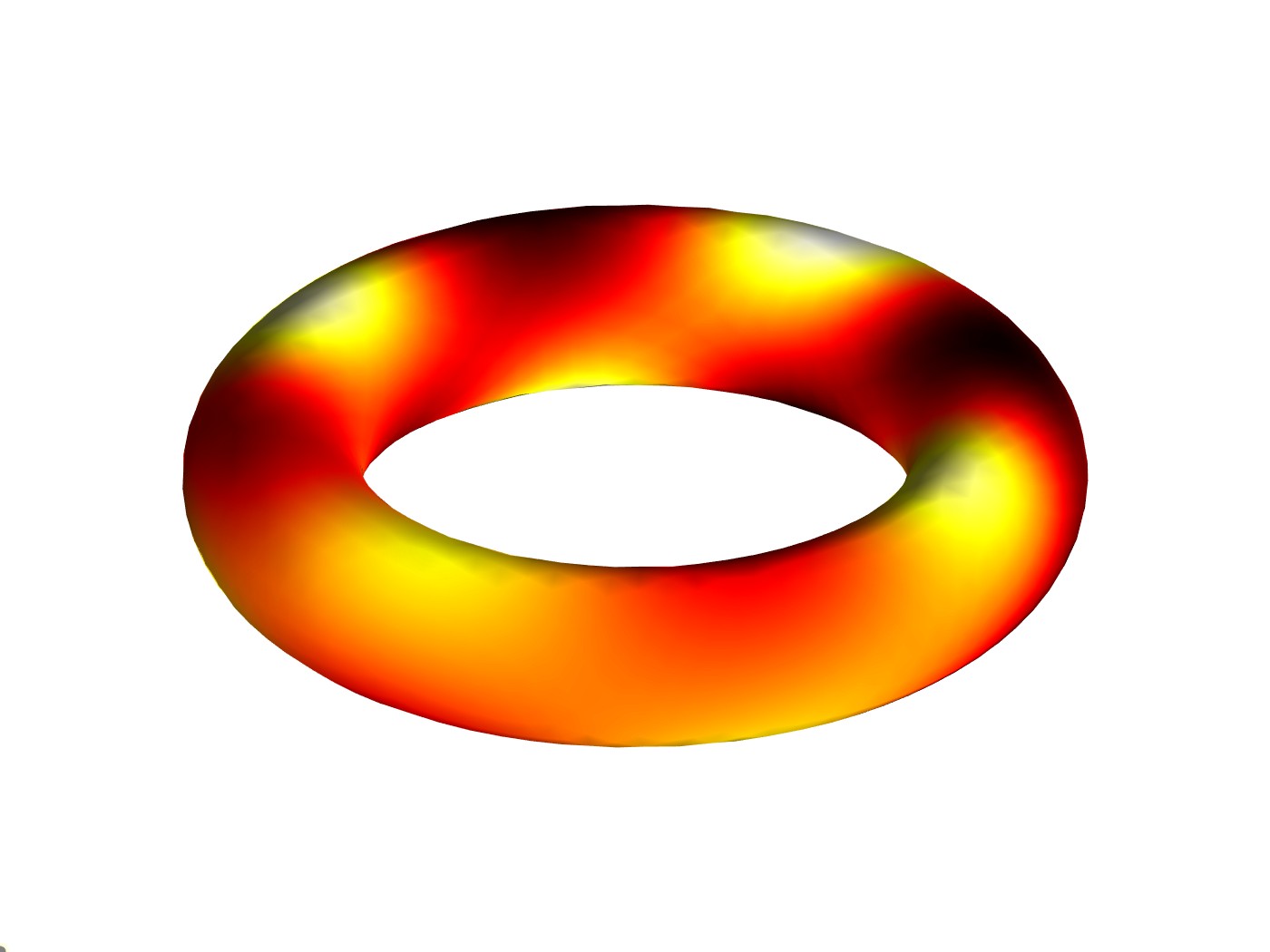}
    \put(45,70){\textbf{$w$}}
\end{overpic}
\begin{overpic}[width=0.24\textwidth,trim=130 150 130 150, clip=true,tics=10]{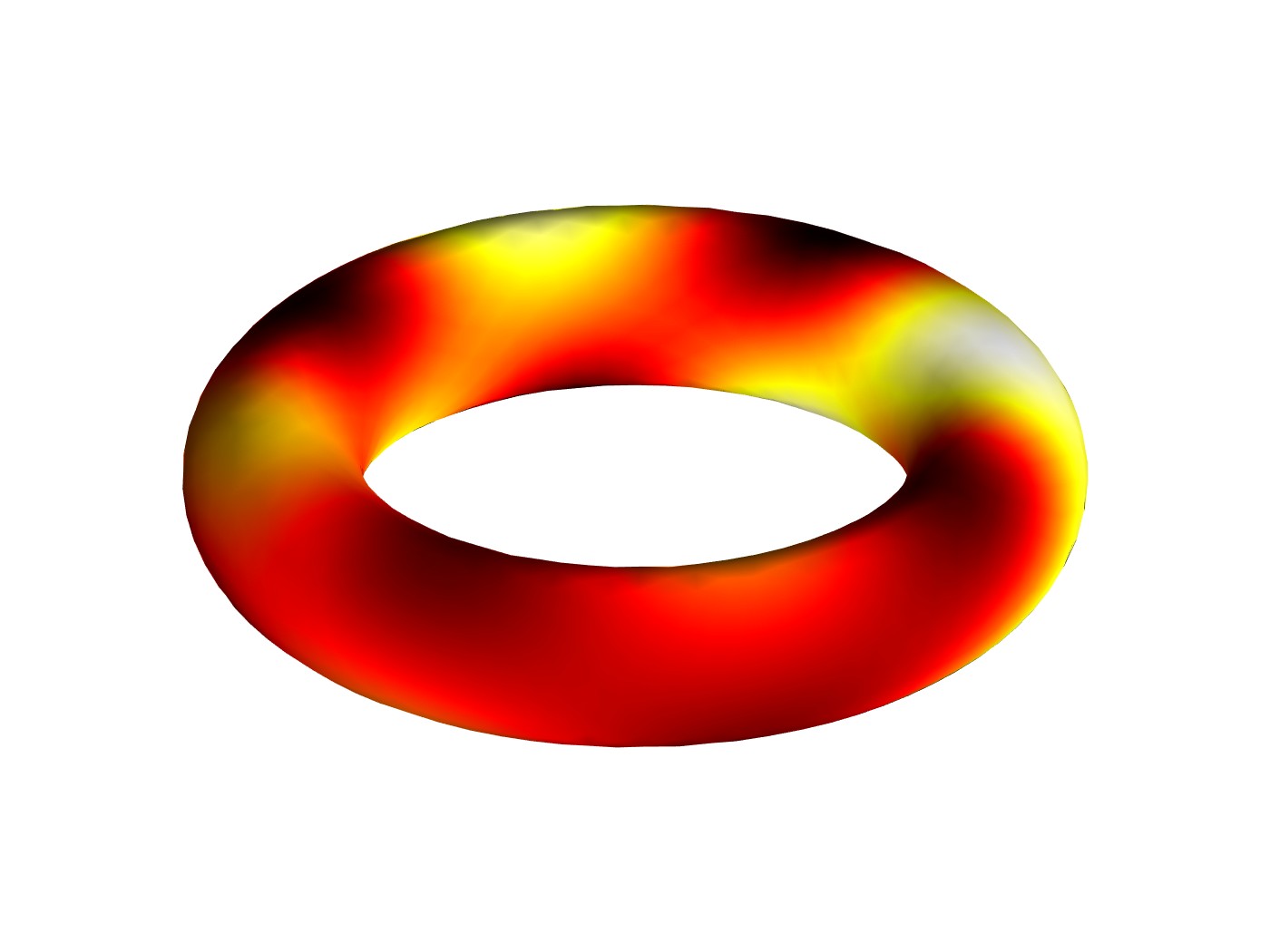}
    \put(45,70){\textbf{$s$}}
\end{overpic}
\\
\begin{overpic}[width=0.24\textwidth,trim=130 150 130 150, clip=true,tics=10]{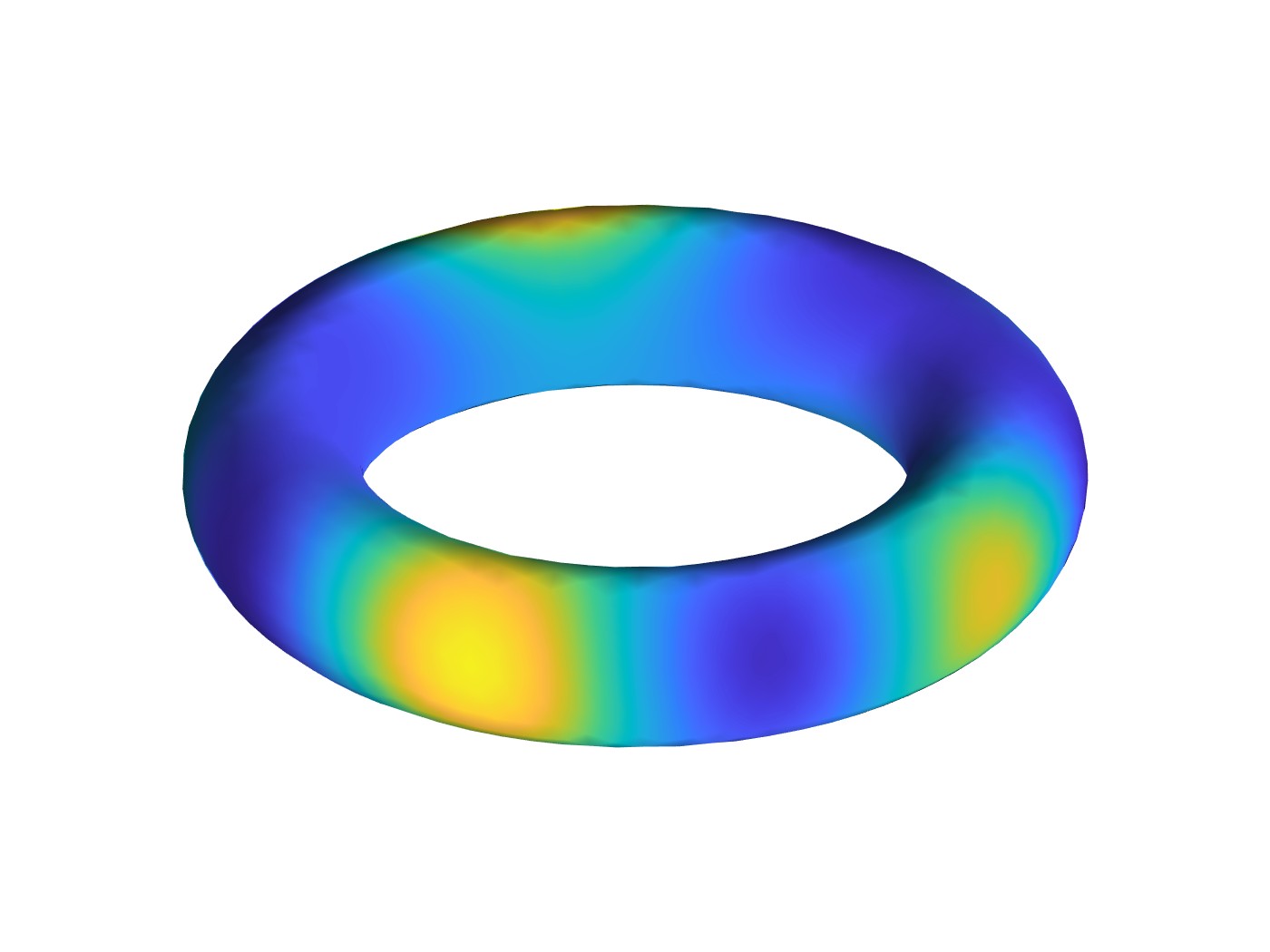}
    \put(-20,20){\rotatebox{90}{\textbf{With}}}
    \put(-10,16){\rotatebox{90}{\textbf{Greedy}}}
\end{overpic}
\begin{overpic}[width=0.24\textwidth,trim=130 150 130 150, clip=true,tics=10]{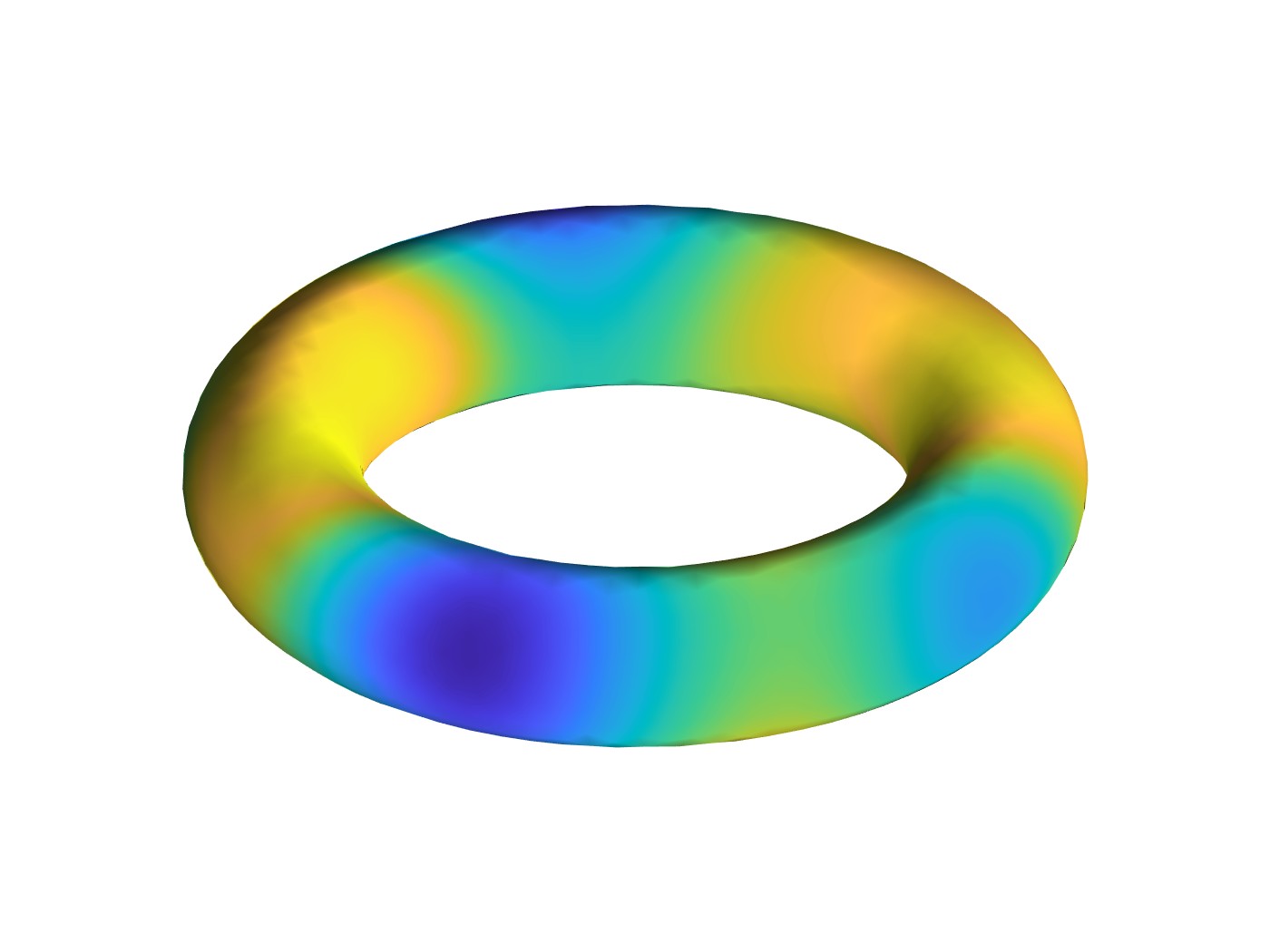}
\end{overpic}
\begin{overpic}[width=0.24\textwidth,trim=130 150 130 150, clip=true,tics=10]{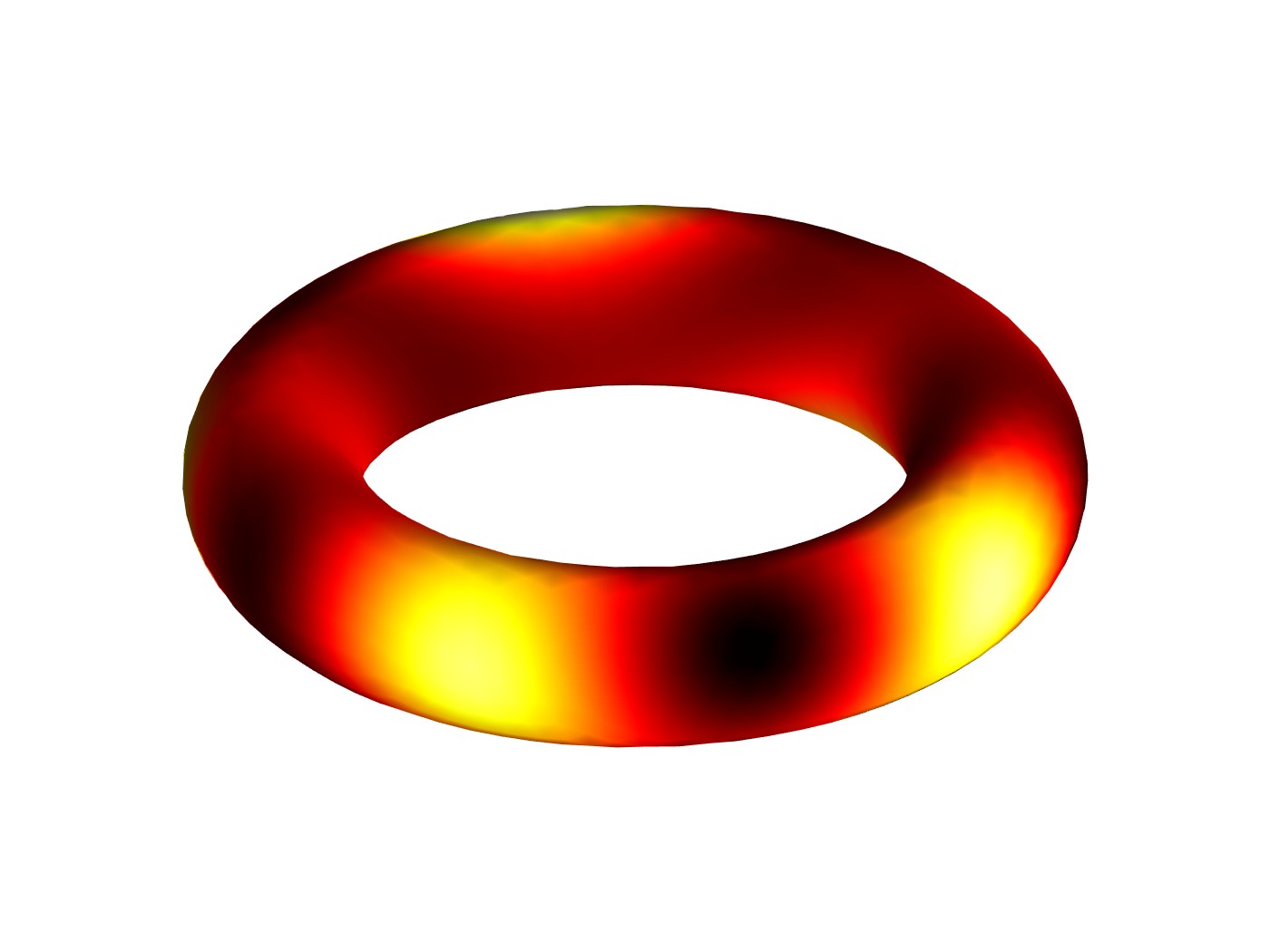}
\end{overpic}
\begin{overpic}[width=0.24\textwidth,trim=130 150 130 150, clip=true,tics=10]{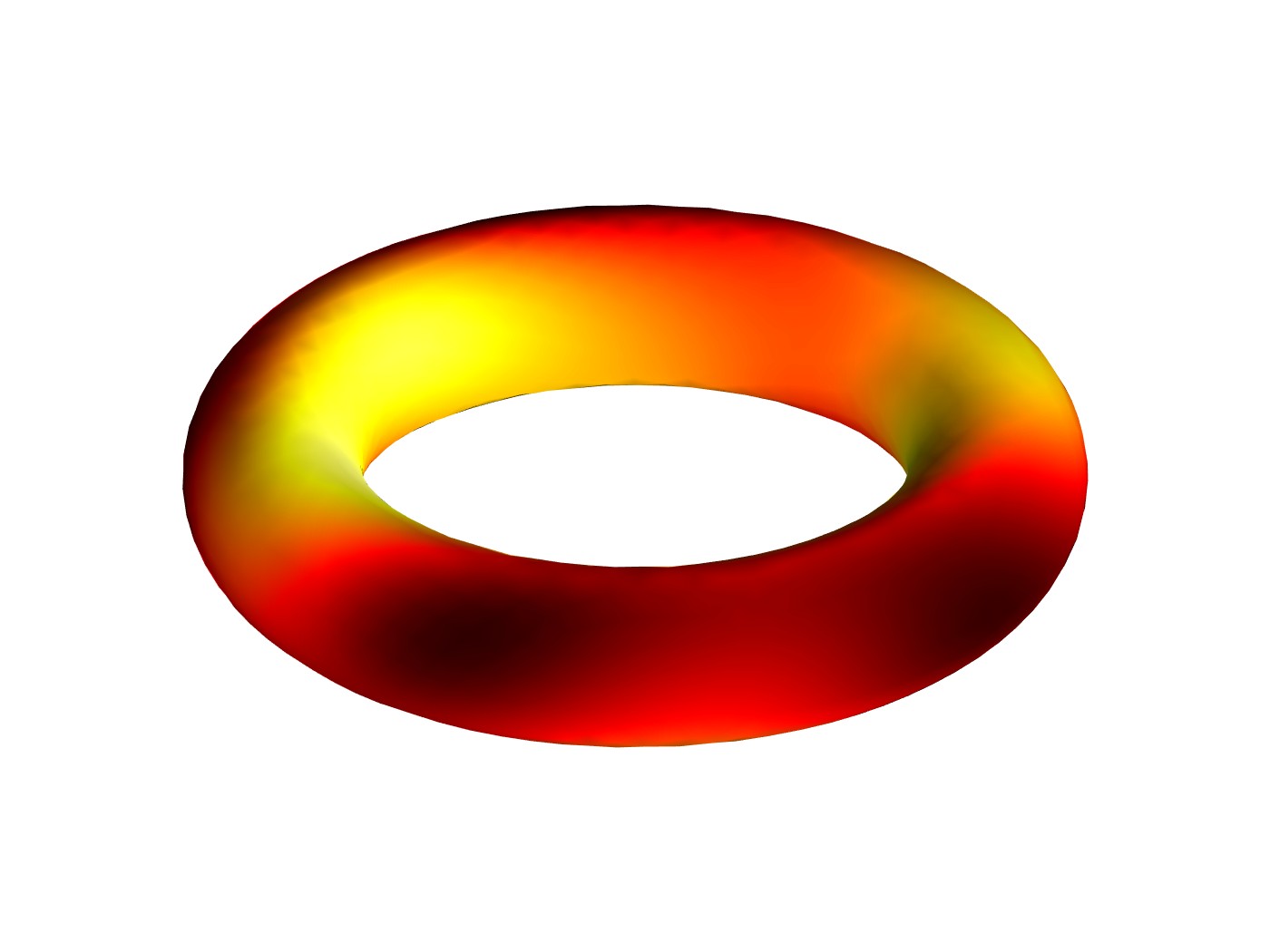}
\end{overpic}\\
\caption{Example~\ref{eg: CBSRD 3}: Bulk and surface solutions by solving the 3D coupled bulk-surface reaction-diffusion equation, using the meshfree time stepping method with and without the greedy algorithm in/on a torus. The meshfree method uses $n = 2644+1430$, $\epsilon = 1$ and time step $\triangle t = 0.005$. Without the greedy algorithm, faint and irregular patterns emerge, while using the algorithm produces clearer, symmetric and well-separated spot patterns.
   }\label{fig: CBS3D torus eps=1}
\end{figure}

\begin{observation}[The greedy algorithm improves efficiency in successful runs]

Next, we implement the experiments on the torus with some larger values of the shape parameter, say, $\epsilon = 4$ and $\epsilon=5$. We fix the number of basis functions as $n = 2644+1430$ and use time steps as $\triangle t = 0.01$ and $\dt = 0.005$. 
As seen in Figure \ref{fig: CBS3D torus eps=1}, the patterns of $u$ and $v$, as well as $w$ and $s$, are reversely correlated. Accordingly, we focus our discussion on the solutions of $u$ and $w$ for both the bulk and surface functions.

Figure \ref{fig: CBS3D torus eps=4 5} displays the bulk and surface function solutions in/on torus under different RBF shape parameters and time step settings. We observe that the patterns obtained with and without the greedy algorithm are nearly identical. Across all solutions, the spots tend to converge to the same pattern, with 4 spots in the bulk and 8 spots on the surface. This trend suggests that even with different numbers of basis functions, the greedy algorithm can facilitate pattern formation.

\begin{table}
\centering
\caption{The number of selected basis functions and the stopping criteria corresponding to the solutions corresponding greedy cases in Figure \ref{fig: CBS3D torus eps=4 5}.}\label{ta: CBS3D eps and SC}
\resizebox{\columnwidth}{!}{%
\begin{tabular}{|c||c | c| c| c|c|c|c|c|c|c|}
\hline
\multirow{ 2}{*}{} &  \multicolumn{5}{c|}{$\epsilon = 4$} & \multicolumn{5}{c|}{$\epsilon = 5$}\\
\cline{2-11}
 & $n_{\Omega}'$(2644)   & $n_{\cals}'$(1430)    & SC($\Omega$)   & SC($\cals$)  & $t_{CPU}$ &$n_{\Omega}'$(2644)   & $n_{\cals}'$(1430)    & SC($\Omega$)   & SC($\cals$) & $t_{CPU}$\\ \hline
 $\triangle t = 0.01$   &  753 & 116   & SC-2$'$  & SC-2$'$& 130.56 & 771 & 222   & SC-2$'$  & SC-2$'$ & 61.86\\
$\triangle t = 0.005$   &  776 & 125   & SC-2$'$  & SC-2$'$ & 115.96 &  2039 & 208   & SC-2$'$  & SC-2$'$ & 79.79\\
\hline
\end{tabular}%
}
\end{table}

Additionally, we list the stopping criteria and the number of selected basis functions out of $n = 2644+1430$ in Table \ref{ta: CBS3D eps and SC} provided by the greedy algorithm. Except for the situation when taking $\epsilon =5$ on the surface, taking $\triangle t = 0.005$ tends to select more basis than taking $\triangle t = 0.01$. Since using a different number of basis functions produces nearly identical patterns, it nearly has no effect on pattern formation compared to a solution without a greedy algorithm. Therefore, the larger time step $ \triangle t = 0.01$ with the greedy algorithm appears to be the most effective setting for generating patterns in this experiment.

\end{observation}
\begin{figure}
  \centering
\begin{overpic}[width=0.24\textwidth,trim=130 150 130 150, clip=true,tics=10]{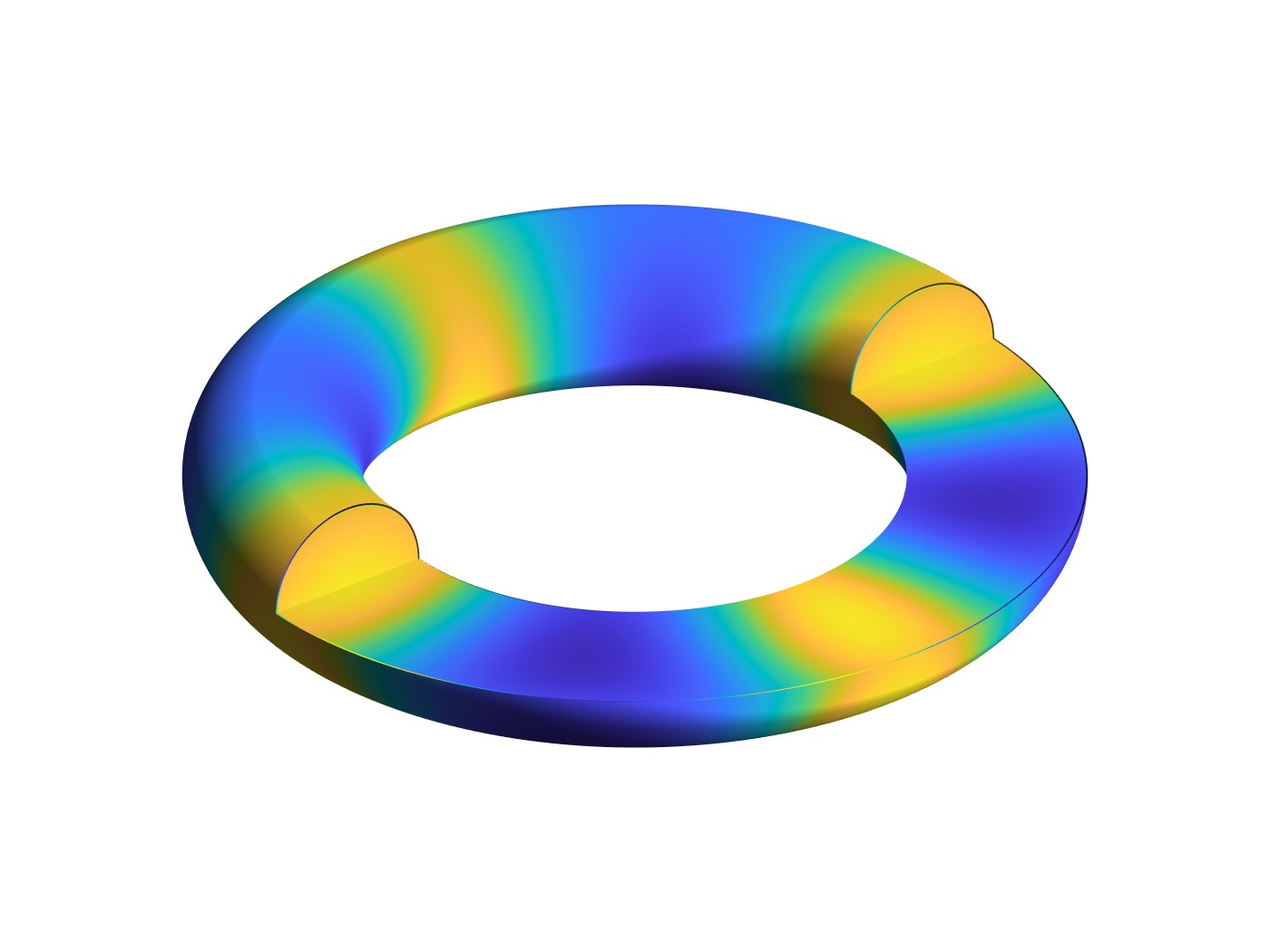}
    \put(75,2){$\yellow{{\raisebox{-.5ex}{\scalebox{2.5}{$\bullet$}}} }\!\times\!4$}
    \put(-25,-35){\rotatebox{90}{\textbf{Without Greedy}}}
    \put(88,75){$\epsilon = 4$}
    \put(45,67){$u$}
    \put(-8,10){\rotatebox{90}{$\dt=0.01$}}
\end{overpic}
\begin{overpic}[width=0.24\textwidth,trim=130 150 130 150, clip=true,tics=10]{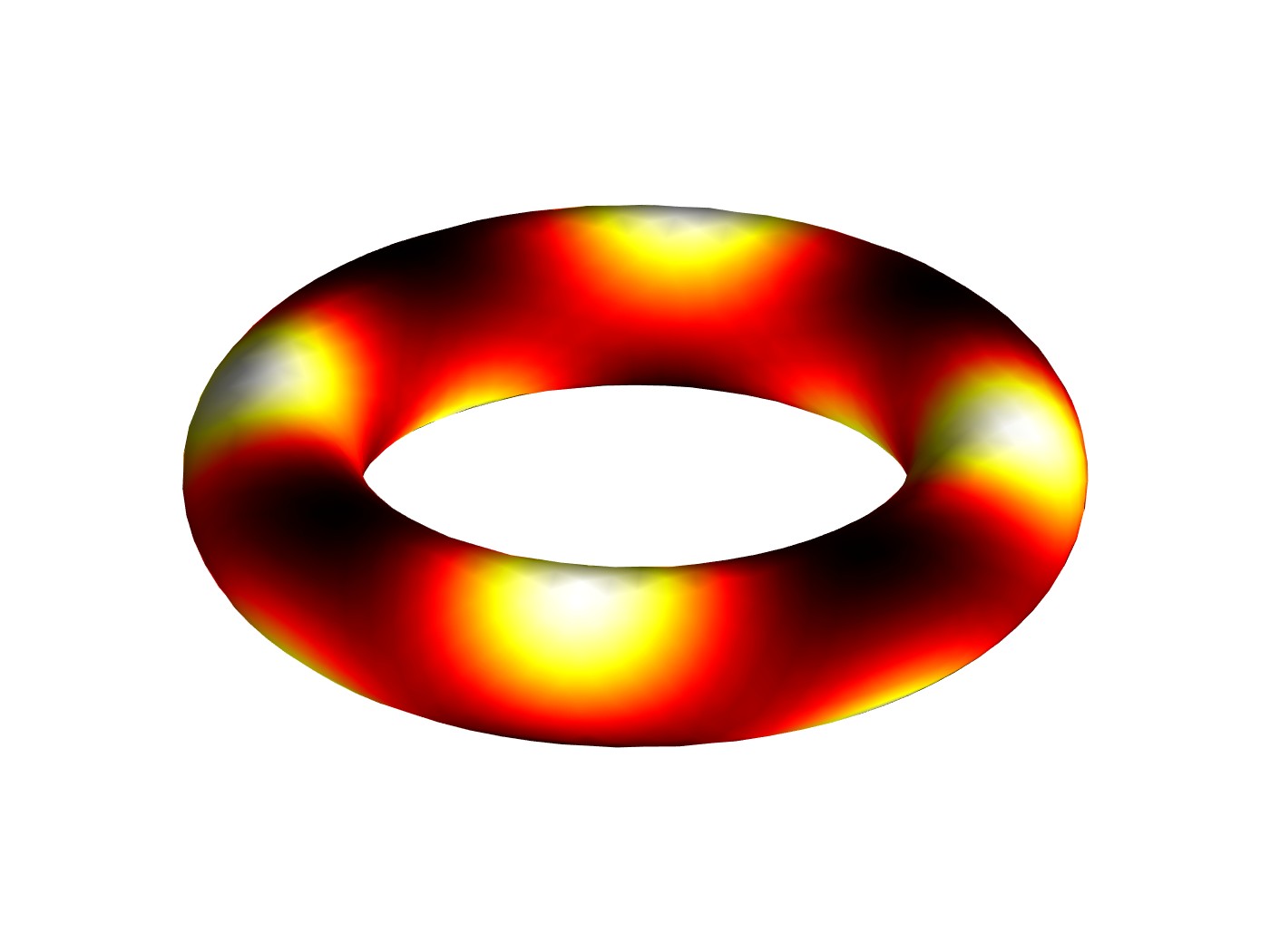}
\put(75,2){$\yellow{{\raisebox{-.5ex}{\scalebox{2.5}{$\bullet$}}} }\!\times\!8$}
\put(45,67){$w$}
\end{overpic}
\begin{overpic}[width=0.24\textwidth,trim=130 150 130 150, clip=true,tics=10]{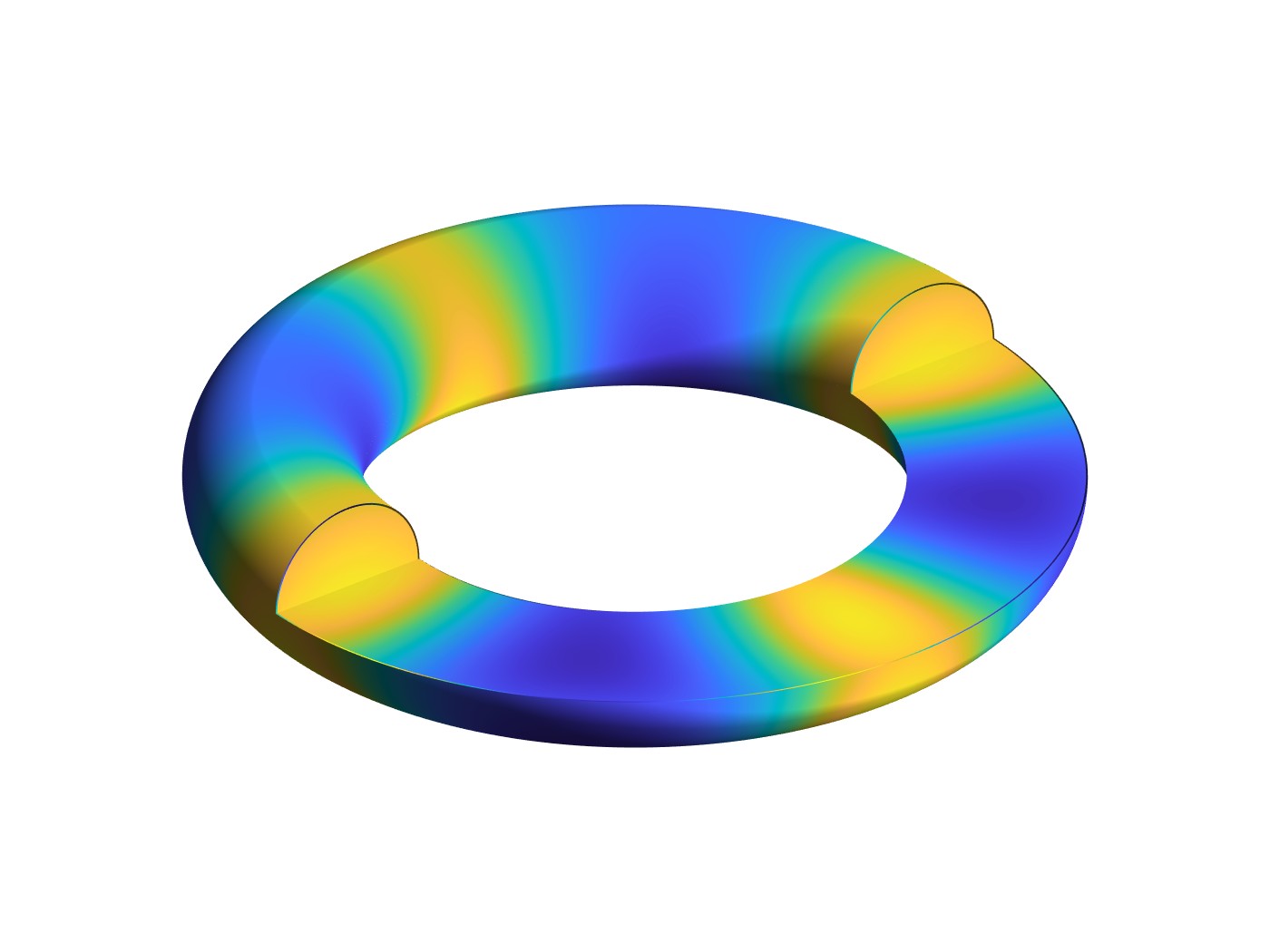}
\put(75,2){$\yellow{{\raisebox{-.5ex}{\scalebox{2.5}{$\bullet$}}} }\!\times\!4$}
\put(88,75){$\epsilon = 5$}
\put(45,67){$u$}
\end{overpic}
\begin{overpic}[width=0.24\textwidth,trim=130 150 130 150, clip=true,tics=10]{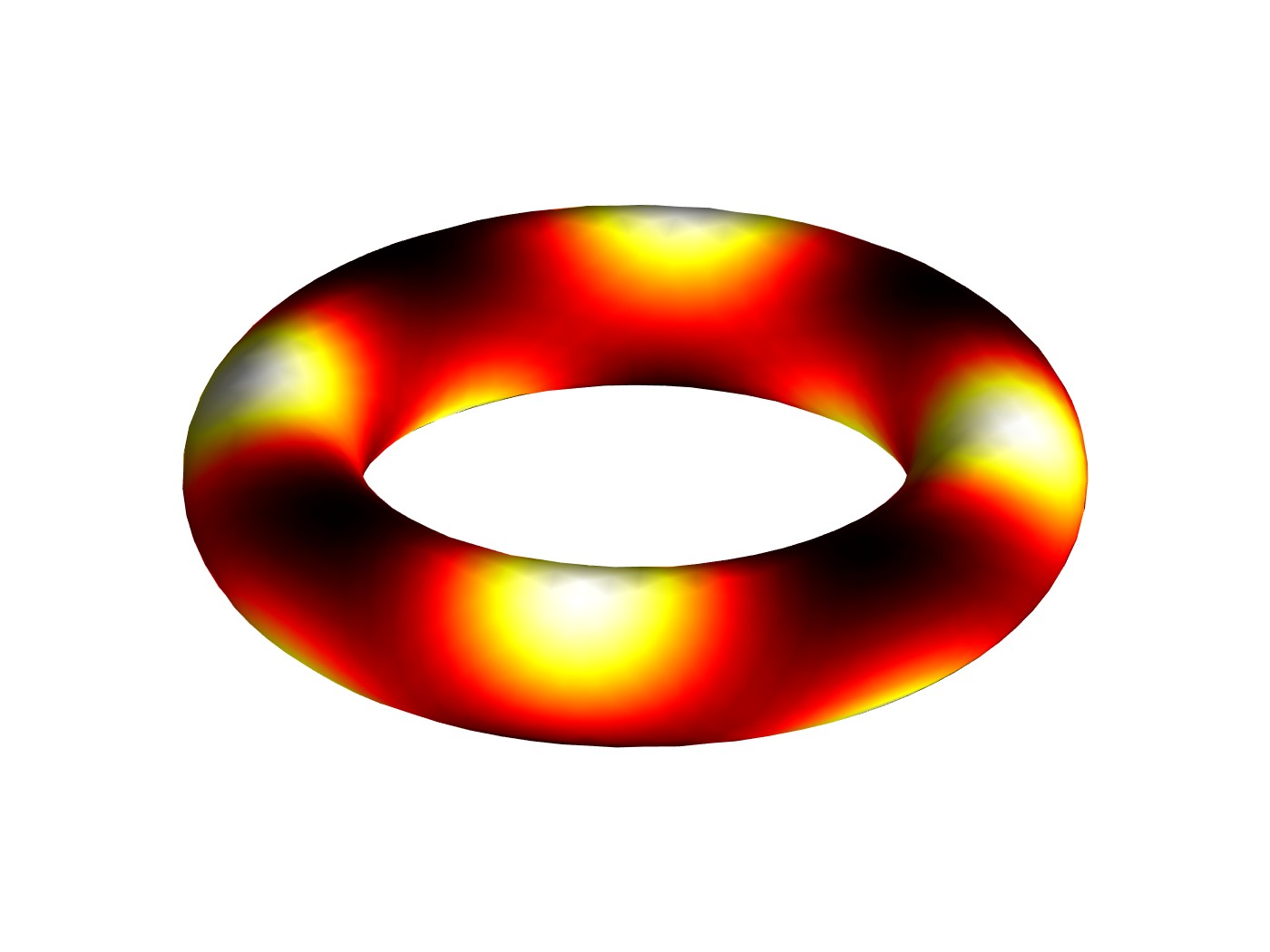}
\put(75,2){$\yellow{{\raisebox{-.5ex}{\scalebox{2.5}{$\bullet$}}} }\!\times\!8$}
\put(45,67){ $w$}
\end{overpic}
\\
  \begin{overpic}[width=0.24\textwidth,trim=130 150 130 150, clip=true,tics=10]{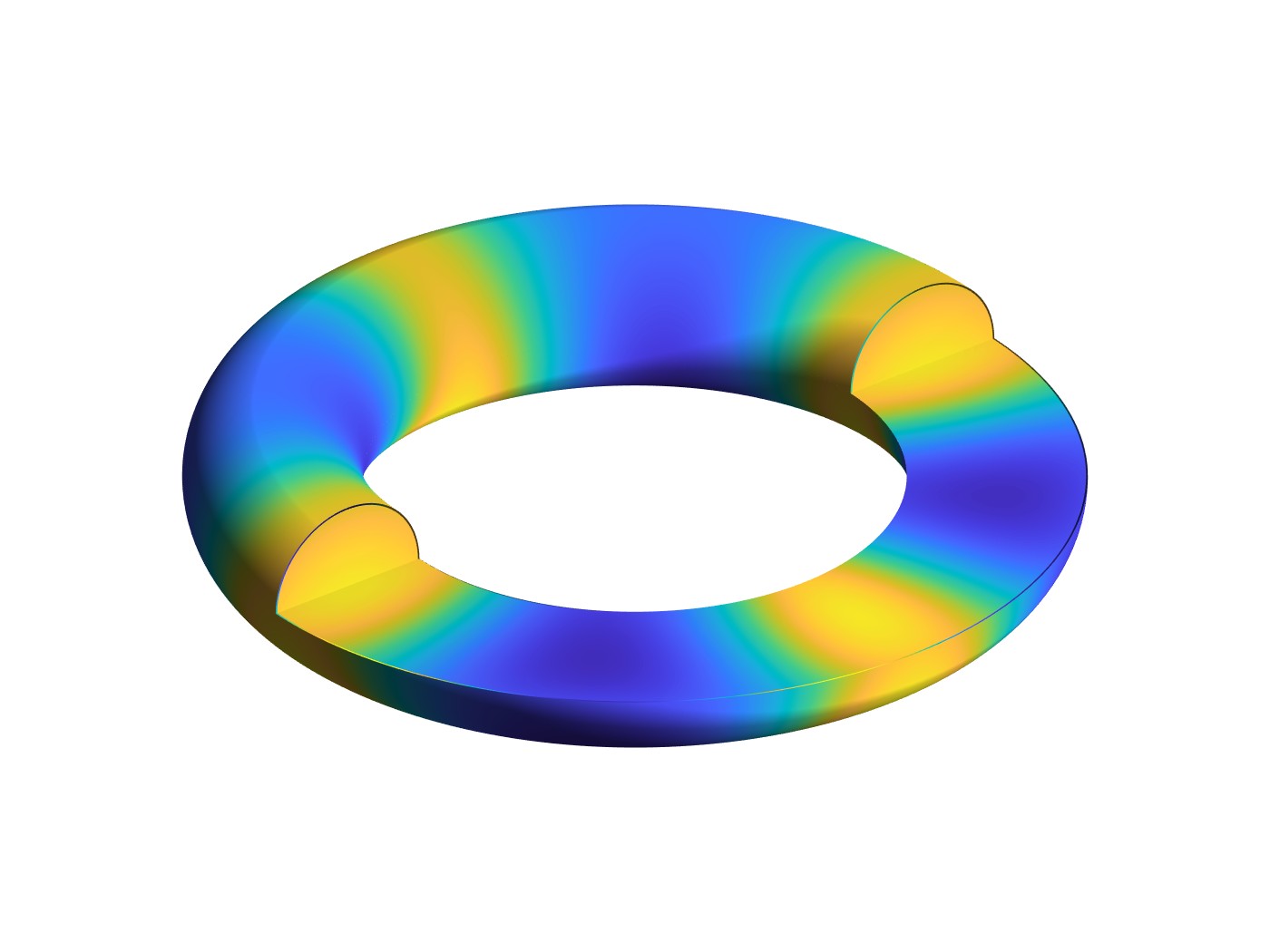}
    \put(75,2){$\yellow{{\raisebox{-.5ex}{\scalebox{2.5}{$\bullet$}}} }\!\times\!4$}
    \put(-8,10){\rotatebox{90}{$\dt=0.005$}}
\end{overpic}
\begin{overpic}[width=0.24\textwidth,trim=130 150 130 150, clip=true,tics=10]{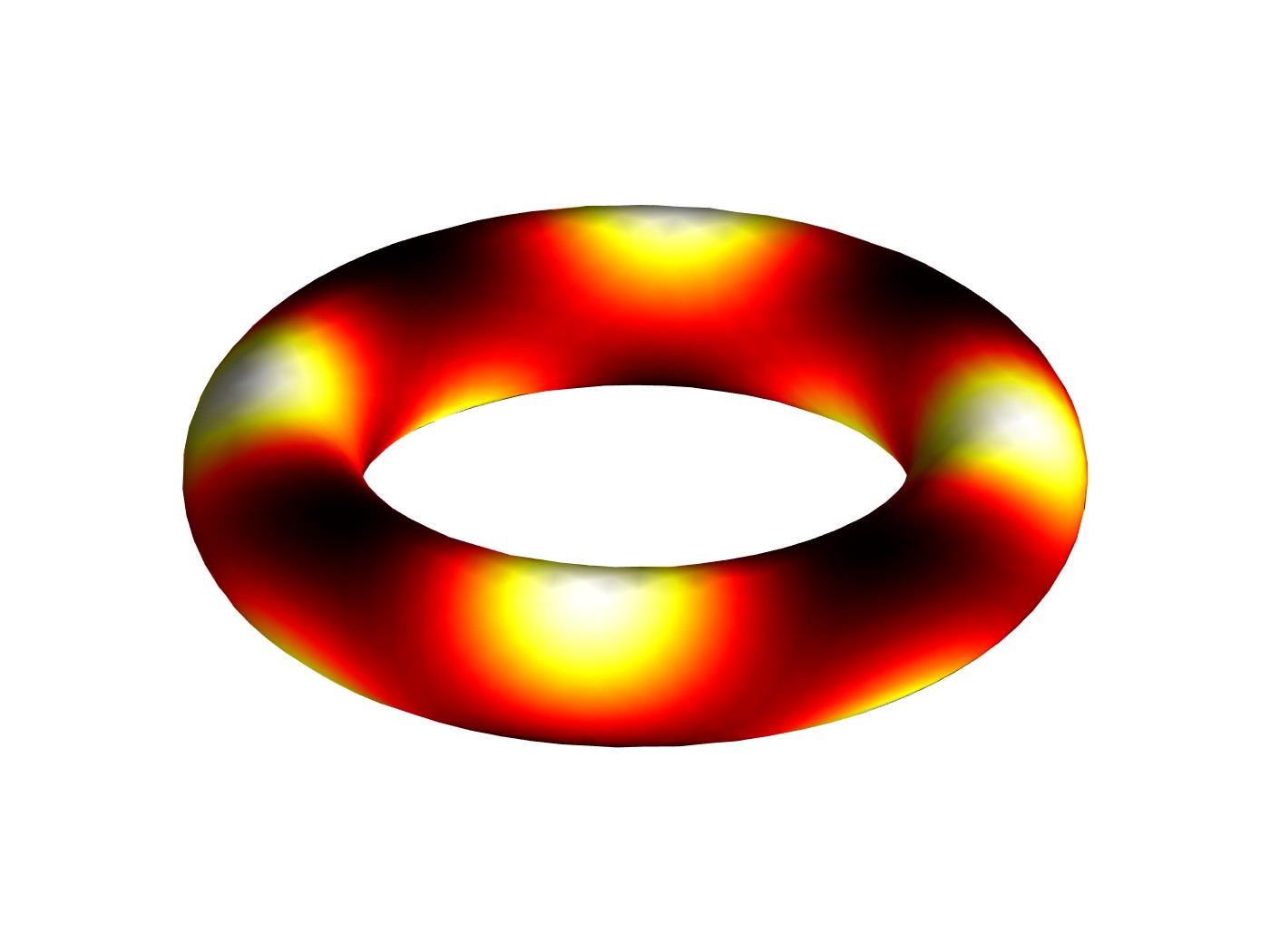}
\put(75,2){$\yellow{{\raisebox{-.5ex}{\scalebox{2.5}{$\bullet$}}} }\!\times\!8$}
\end{overpic}
\begin{overpic}[width=0.24\textwidth,trim=130 150 130 150, clip=true,tics=10]{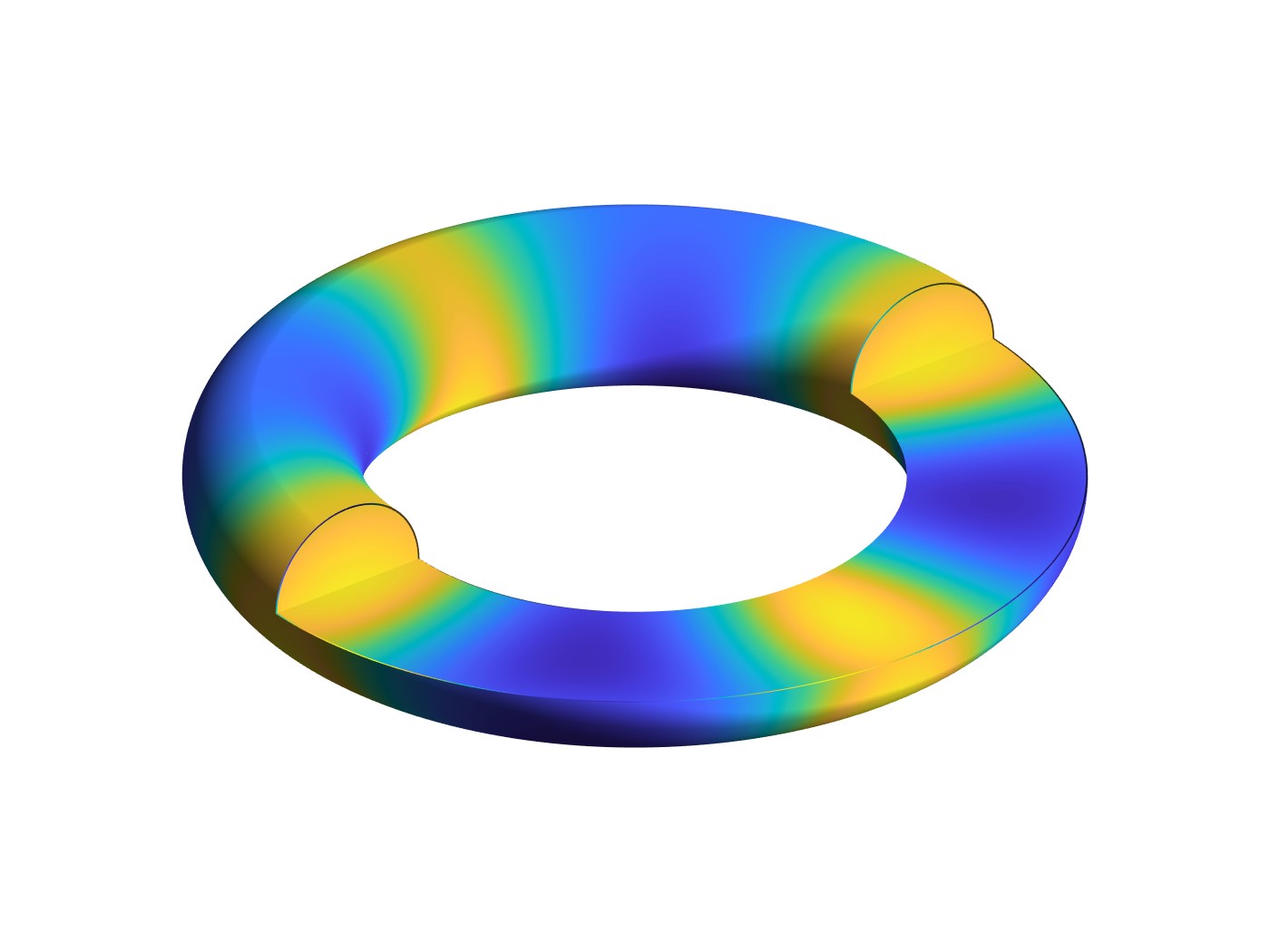}
\put(75,2){$\yellow{{\raisebox{-.5ex}{\scalebox{2.5}{$\bullet$}}} }\!\times\!4$}
\end{overpic}
\begin{overpic}[width=0.24\textwidth,trim=130 150 130 150, clip=true,tics=10]{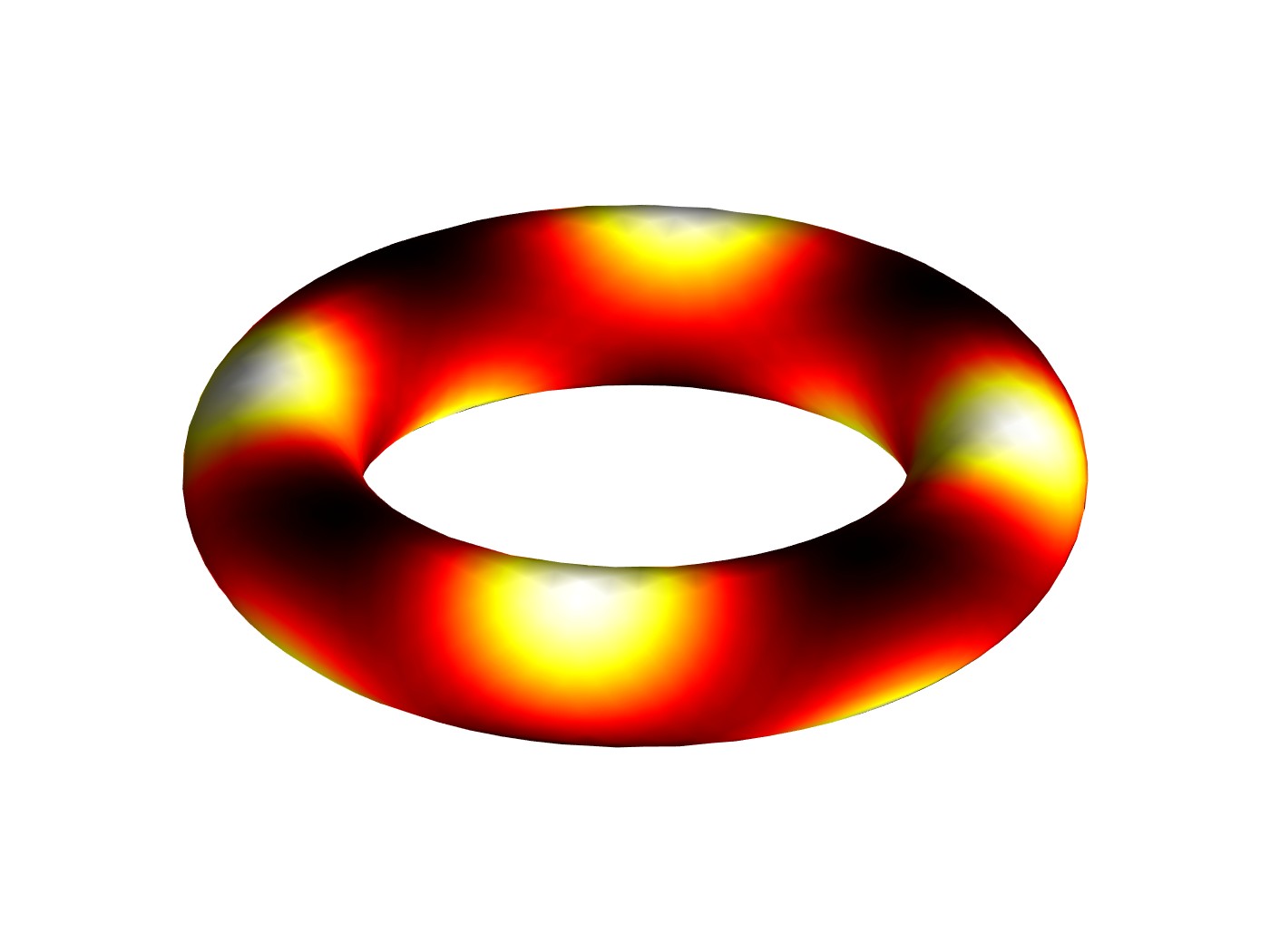}
\put(75,2){$\yellow{{\raisebox{-.5ex}{\scalebox{2.5}{$\bullet$}}} }\!\times\!8$}
\end{overpic}
\\
\begin{overpic}[width=0.24\textwidth,trim=130 150 130 150, clip=true,tics=10]{Figures/BS3D/CBStoruseps4dt1noGdyU.jpg}
    \put(75,2){$\yellow{{\raisebox{-.5ex}{\scalebox{2.5}{$\bullet$}}} }\!\times\!4$}
    \put(-25,-35){\rotatebox{90}{\textbf{With Greedy}}}
    \put(-8,10){\rotatebox{90}{$\dt=0.01$}}
\end{overpic}
\begin{overpic}[width=0.24\textwidth,trim=130 150 130 150, clip=true,tics=10]{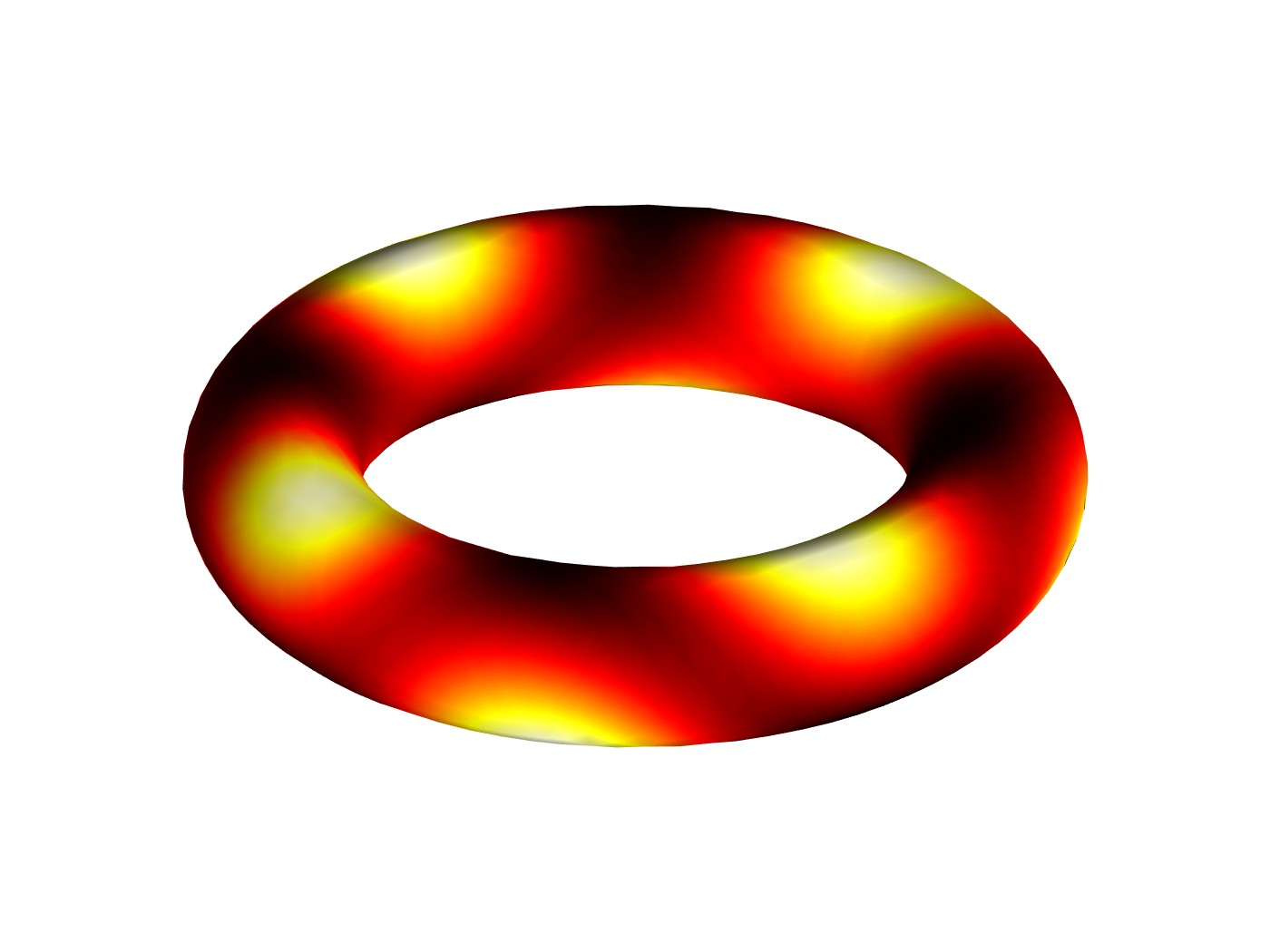}
\put(75,2){$\yellow{{\raisebox{-.5ex}{\scalebox{2.5}{$\bullet$}}} }\!\times\!8$}
\end{overpic}
\begin{overpic}[width=0.24\textwidth,trim=130 150 130 150, clip=true,tics=10]{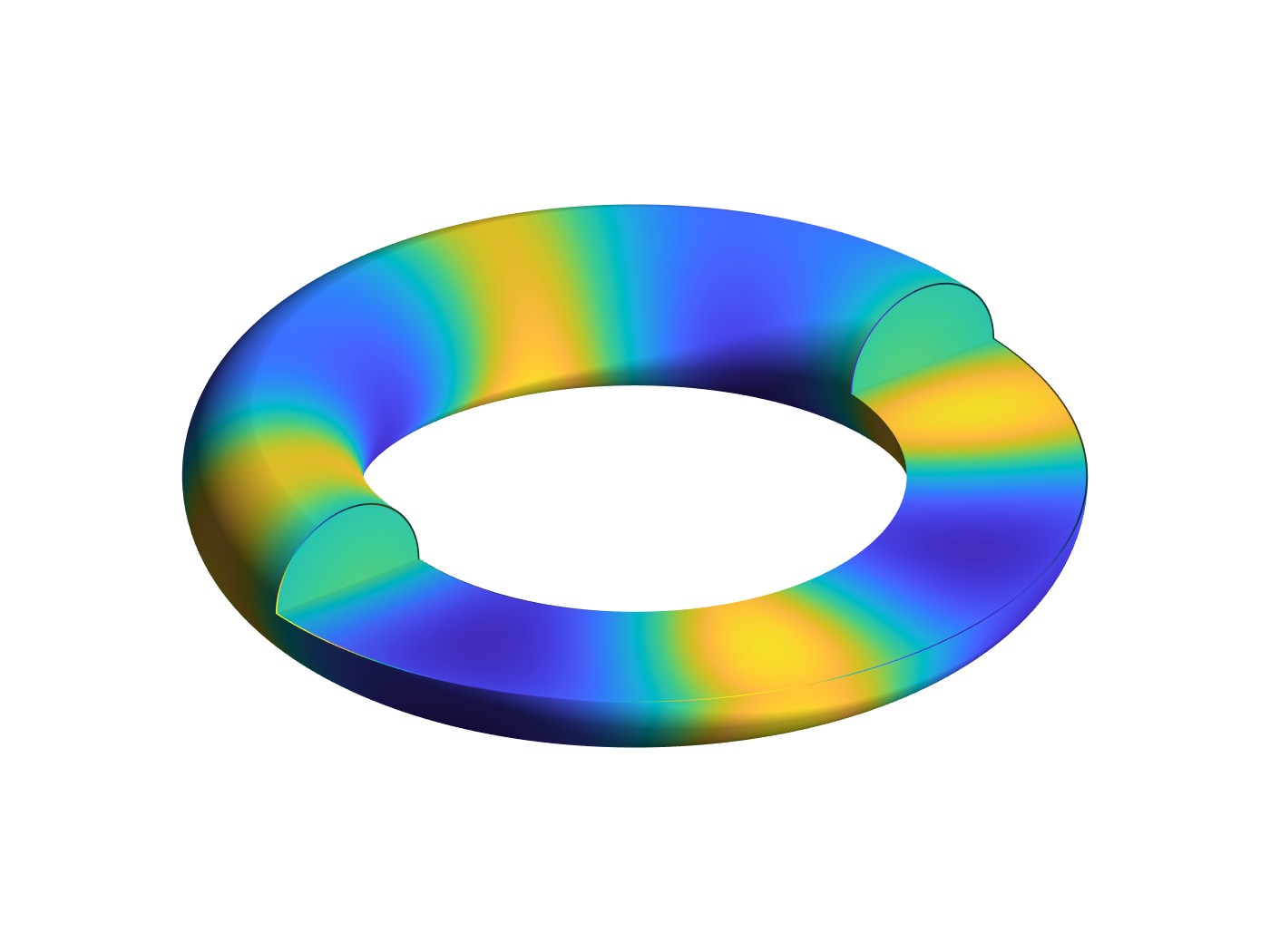}
\put(75,2){$\yellow{{\raisebox{-.5ex}{\scalebox{2.5}{$\bullet$}}} }\!\times\!4$}
\end{overpic}
\begin{overpic}[width=0.24\textwidth,trim=130 150 130 150, clip=true,tics=10]{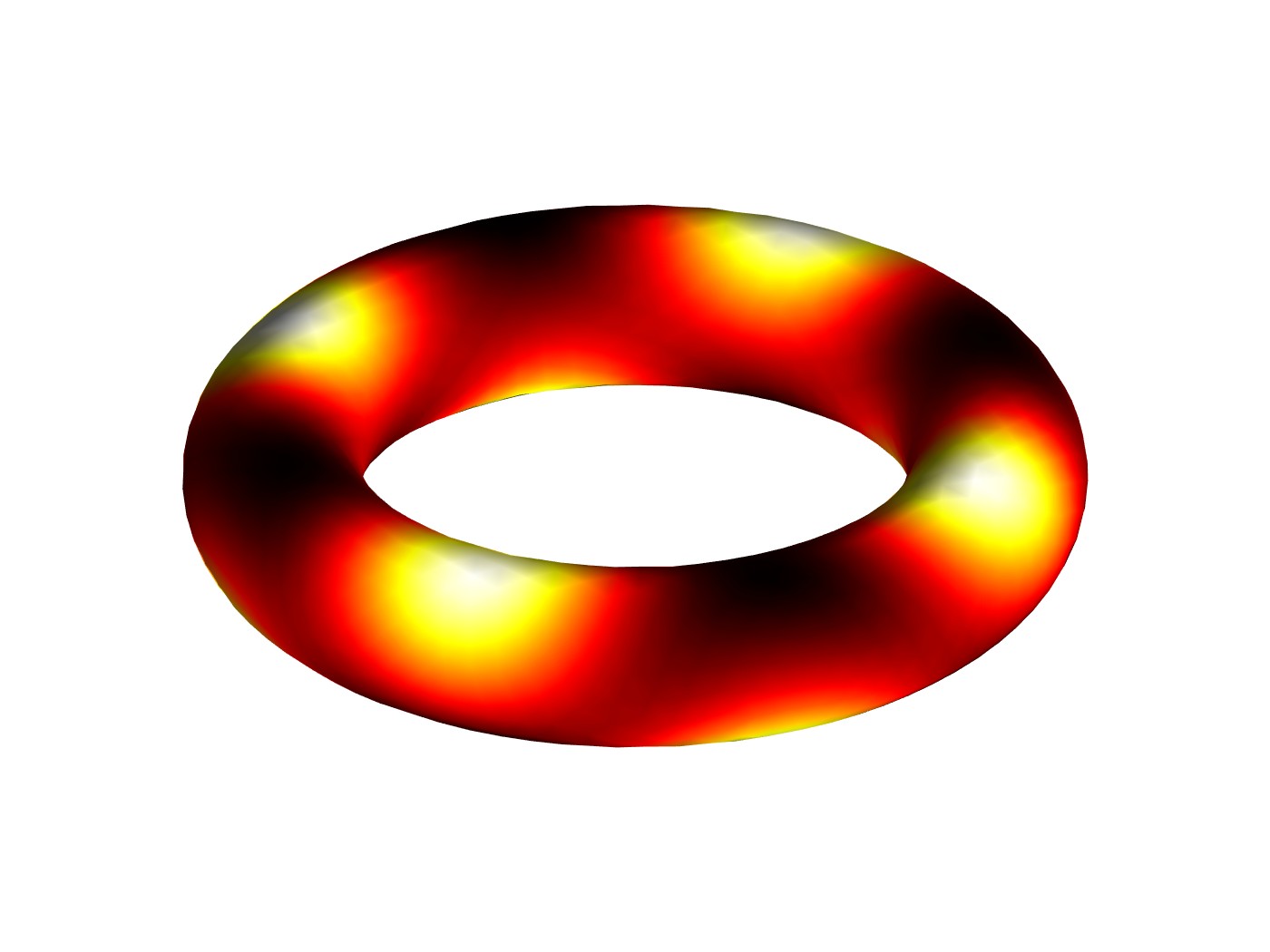}
\put(75,2){$\yellow{{\raisebox{-.5ex}{\scalebox{2.5}{$\bullet$}}} }\!\times\!8$}
\end{overpic}
\\
\begin{overpic}[width=0.24\textwidth,trim=130 150 130 150, clip=true,tics=10]{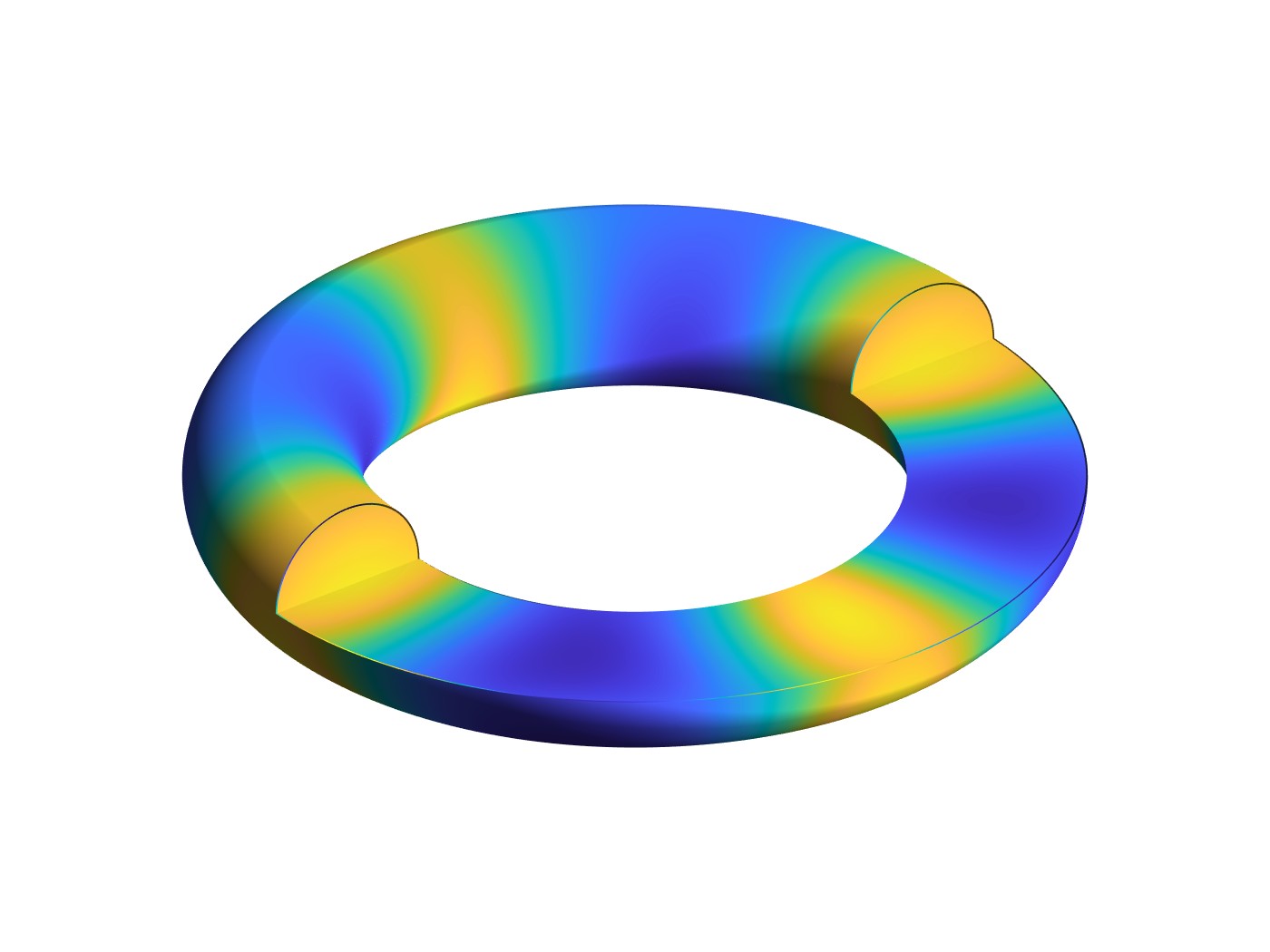}
\put(75,2){$\yellow{{\raisebox{-.5ex}{\scalebox{2.5}{$\bullet$}}} }\!\times\!4$}
\put(-8,10){\rotatebox{90}{$\dt=0.005$}}
\end{overpic}
\begin{overpic}[width=0.24\textwidth,trim=130 150 130 150, clip=true,tics=10]{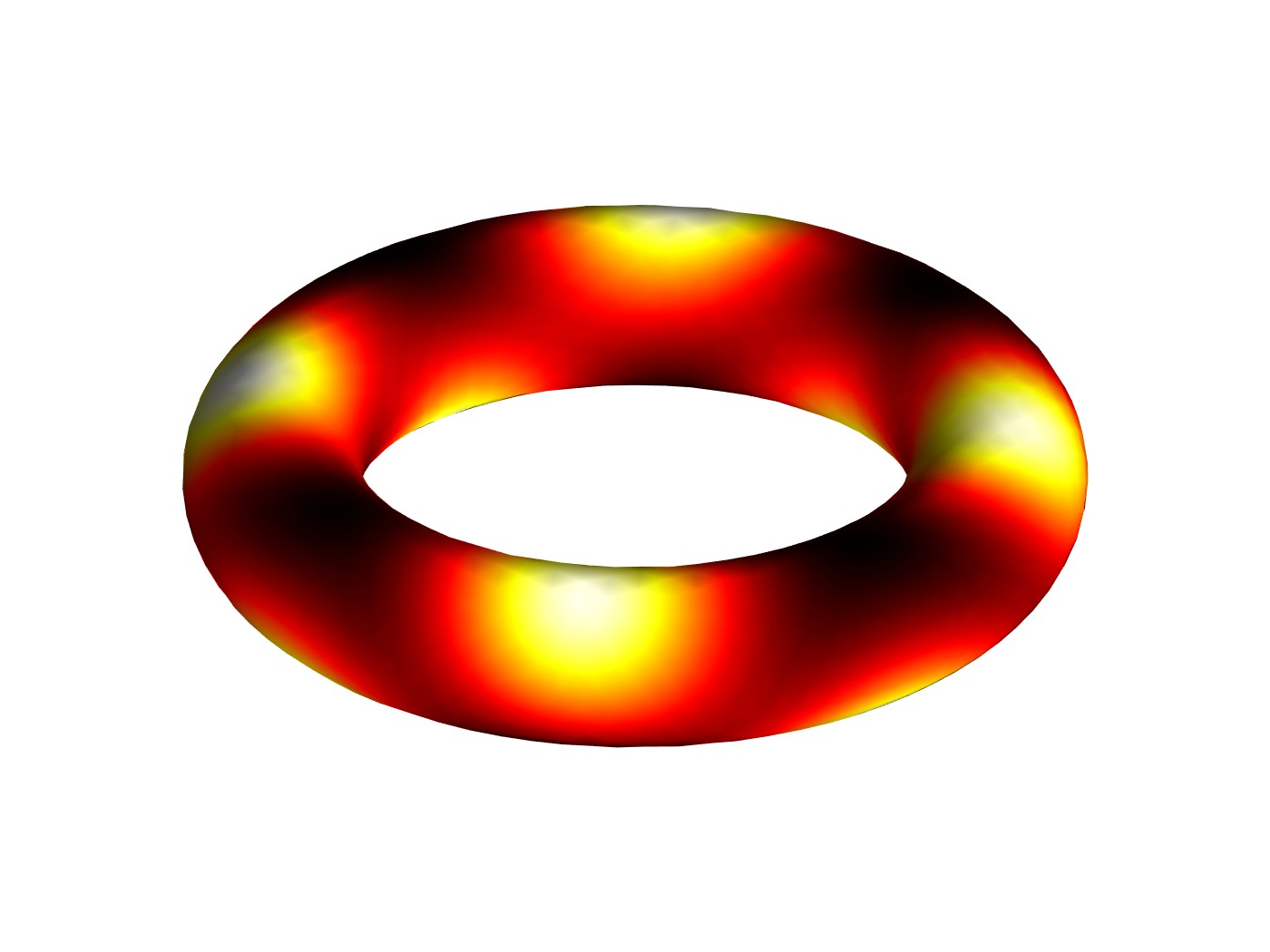}
\put(75,2){$\yellow{{\raisebox{-.5ex}{\scalebox{2.5}{$\bullet$}}} }\!\times\!8$}
\end{overpic}
\begin{overpic}[width=0.24\textwidth,trim=130 150 130 150, clip=true,tics=10]{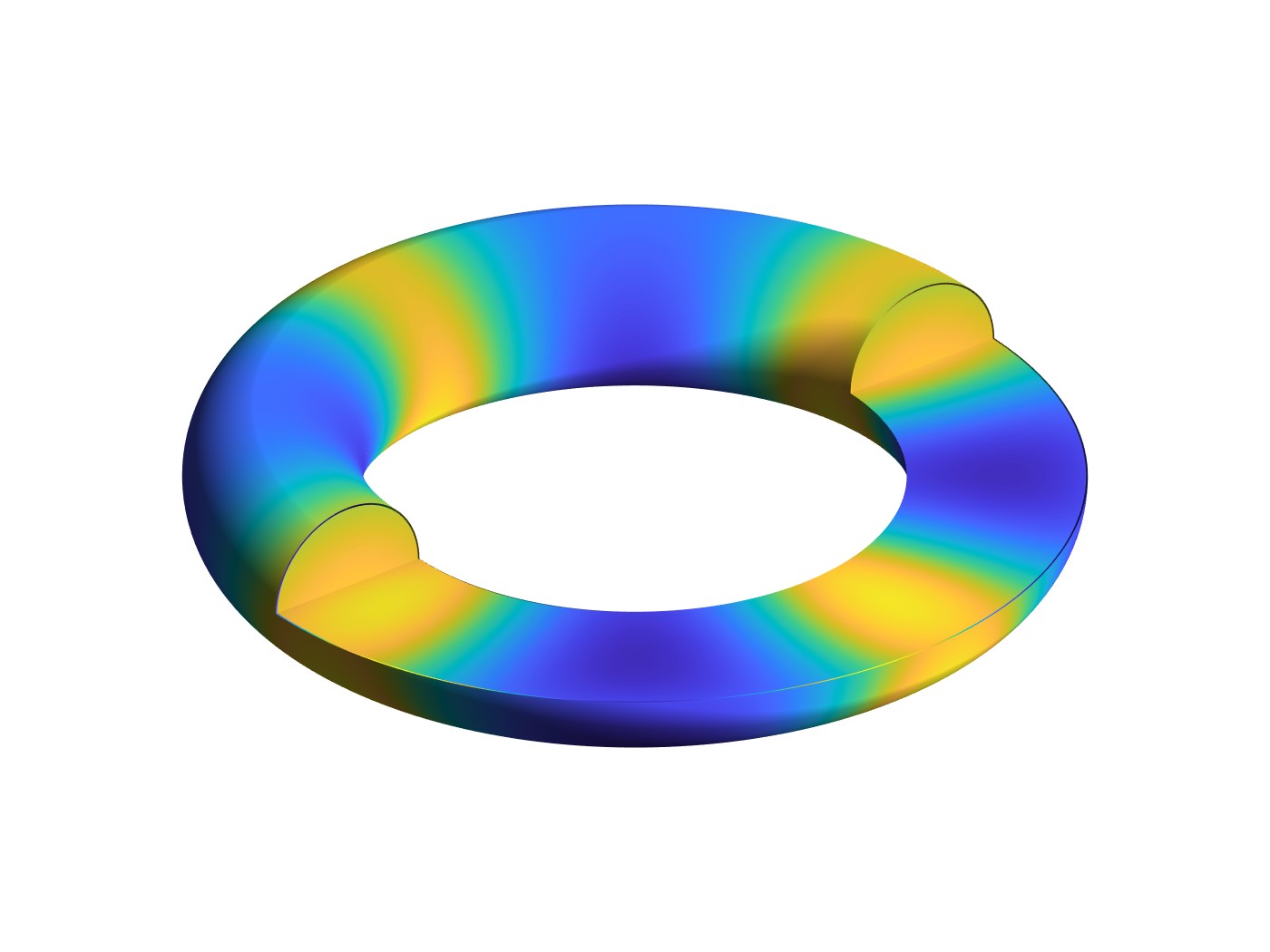}
\put(75,2){$\yellow{{\raisebox{-.5ex}{\scalebox{2.5}{$\bullet$}}} }\!\times\!4$}
\end{overpic}
\begin{overpic}[width=0.24\textwidth,trim=130 150 130 150, clip=true,tics=10]{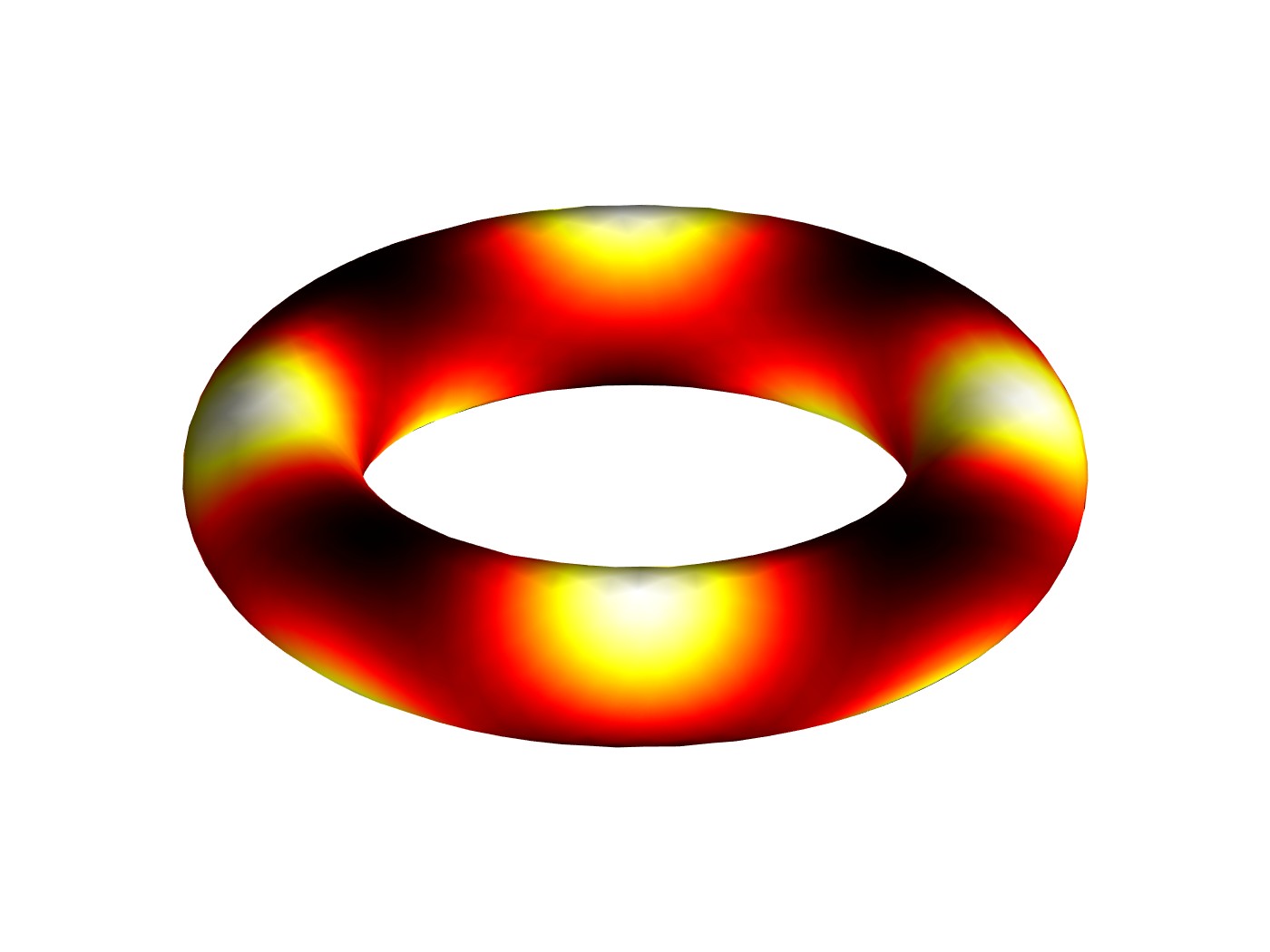}
\put(75,2){$\yellow{{\raisebox{-.5ex}{\scalebox{2.5}{$\bullet$}}} }\!\times\!8$}
\end{overpic}
\caption{Example~\ref{eg: CBSRD 3}: Bulk and surface solutions by solving the 3D coupled bulk-surface reaction-diffusion equation corresponding with Figure \ref{fig: CBS3D torus eps=1} with time steps $\triangle t =0.01$ and $\triangle t = 0.005$ and larger RBF shape parameters $\epsilon = 4$ and $\epsilon = 5$. The number of bases used and stopping criteria in the greedy cases are in Table \ref{ta: CBS3D eps and SC}.
   }\label{fig: CBS3D torus eps=4 5}
\end{figure}
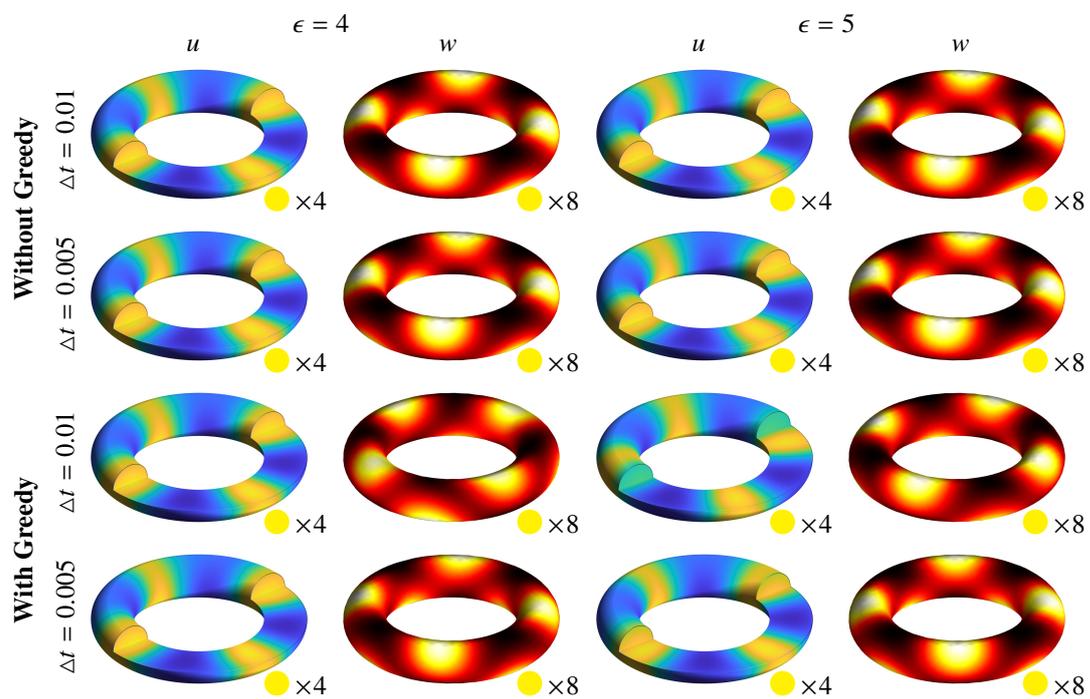

\begin{observation}[The greedy algorithm performs robustly in various domains]

In our continued exploration of the effectiveness of our approach, we further examine the solutions of (\ref{eq: CBS}) on additional 3D domains. Figures \ref{fig: CBS3D cycl} and \ref{fig: CBS3D ellip} present solutions of $u$ and $w$ in a Dupin's cyclide and an ellipsoid. For the cyclide, we use $\epsilon =4$, and the time stepping method with $\triangle t =0.01$ and $\triangle t = 0.005$. Due to the complex geometry, the spots formed in/on the cyclide are not as uniformly distributed as in the torus. Nonetheless, we observe that the greedy algorithm produces the same number of clear and separated spots as the non-greedy solutions, both in the bulk and on the surface. Specifically, there are consistently 6 spots in the bulk and 10 on the surface. However, under different settings, the distributions of spots may not always be identical. For example, when applying the greedy algorithm with $\triangle t = 0.01$, the solution in the cyclide domain appears different from other solutions. This difference is due to the fact that the positions of the spots are influenced by the number and position of the selected basis functions.

Likewise, the solutions on the ellipsoid with the greedy algorithm also exhibit slight differences when compared to those without. Specifically, when we apply the greedy algorithm with $\dt = 0.005$, the spots on the top of the ellipsoid appear different from the patterns in other solutions. Nonetheless, despite using or not using the greedy algorithm, the number of spots in the bulk and on the surface remains almost the same.

\end{observation}
\end{example}

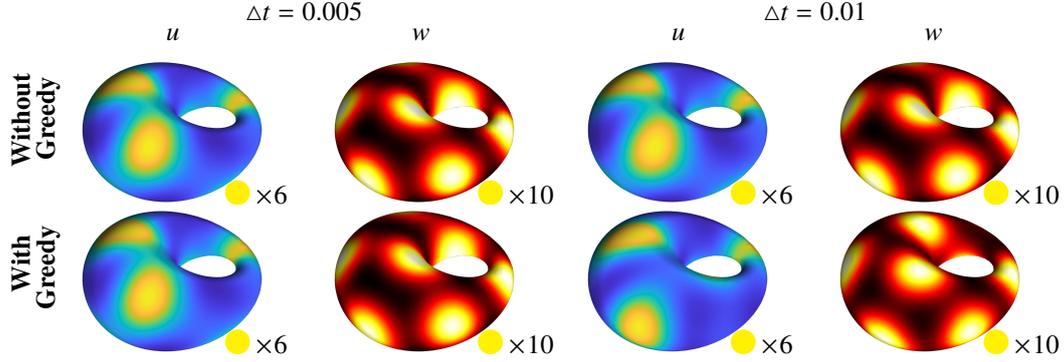
\begin{figure}
  \centering
  \begin{overpic}[width=0.24\textwidth,trim=130 180 130 180, clip=true,tics=10]{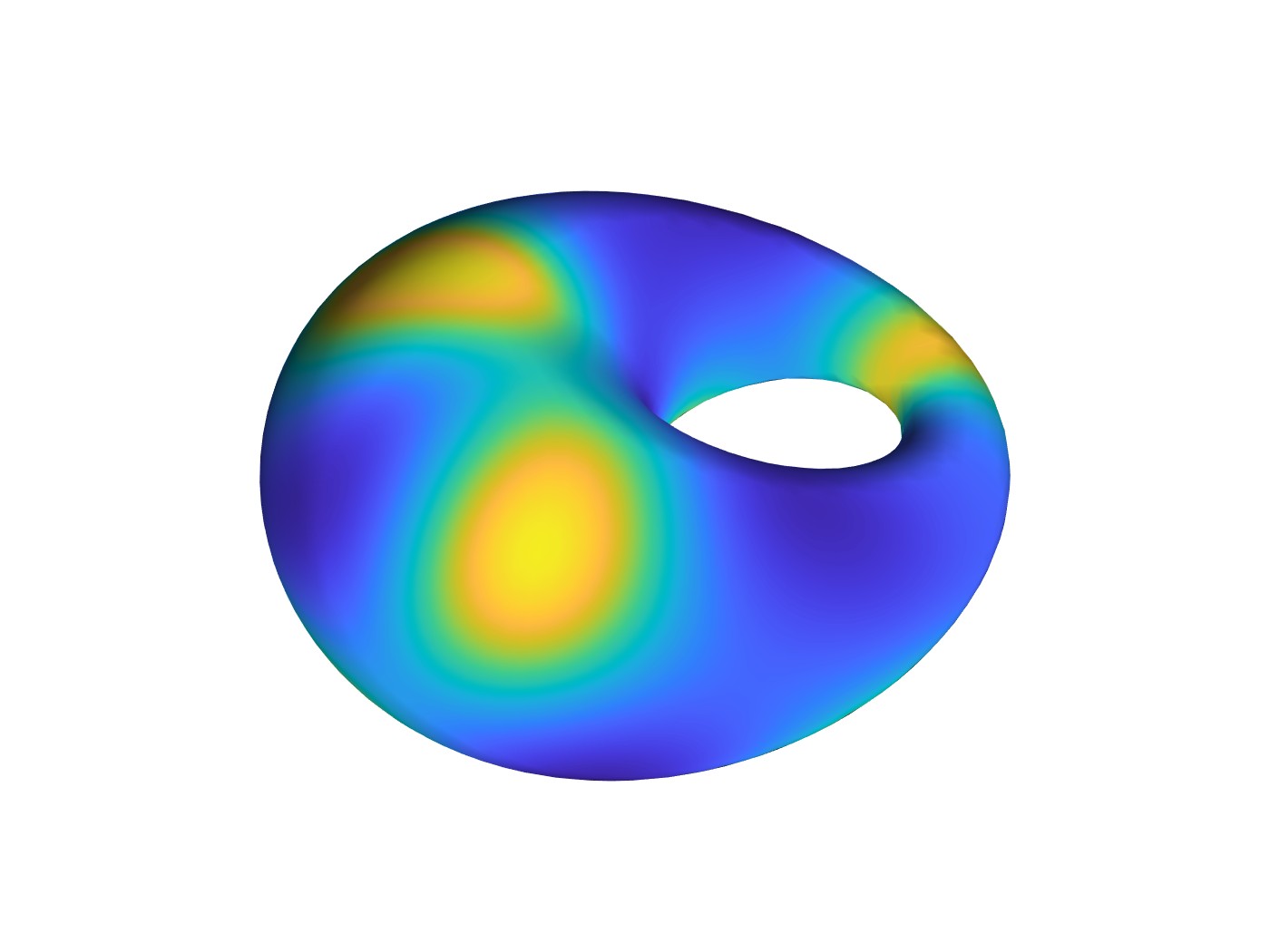}
    \put(-15,15){\rotatebox{90}{\textbf{Without}}}
    \put(-5,18){\rotatebox{90}{\textbf{Greedy}}}
    \put(80,75){$\dt=0.005$}
    \put(45,67){ $u$}
    \put(70,0){$\yellow{{\raisebox{-.5ex}{\scalebox{2.5}{$\bullet$}}} }\!\times\!6$}
\end{overpic}
\begin{overpic}[width=0.24\textwidth,trim=130 180 130 180, clip=true,tics=10]{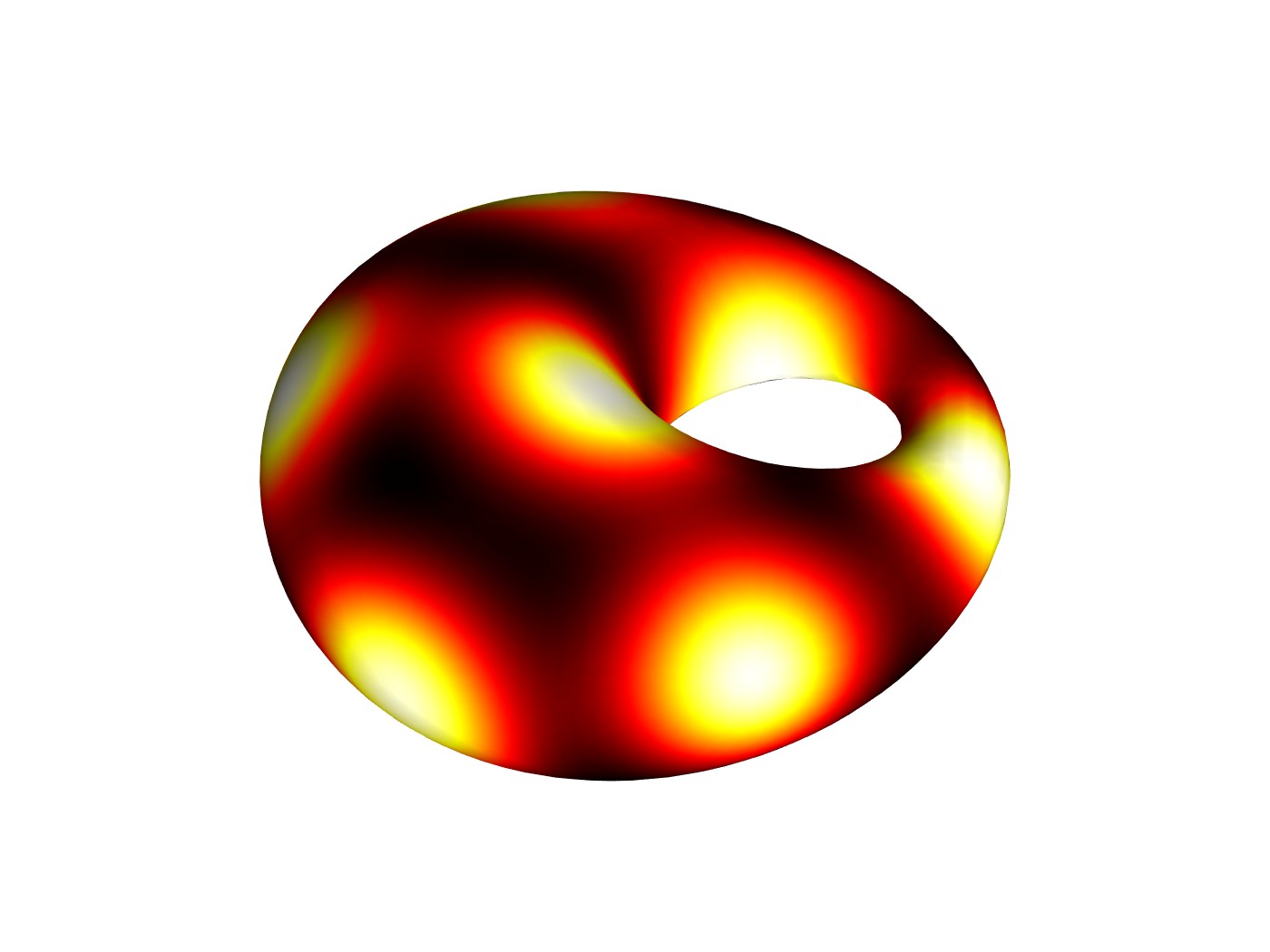}
\put(70,0){$\yellow{{\raisebox{-.5ex}{\scalebox{2.5}{$\bullet$}}} }\!\times\!10$}
\put(45,67){$w$}
\end{overpic}
\begin{overpic}[width=0.24\textwidth,trim=130 180 130 180, clip=true,tics=10]{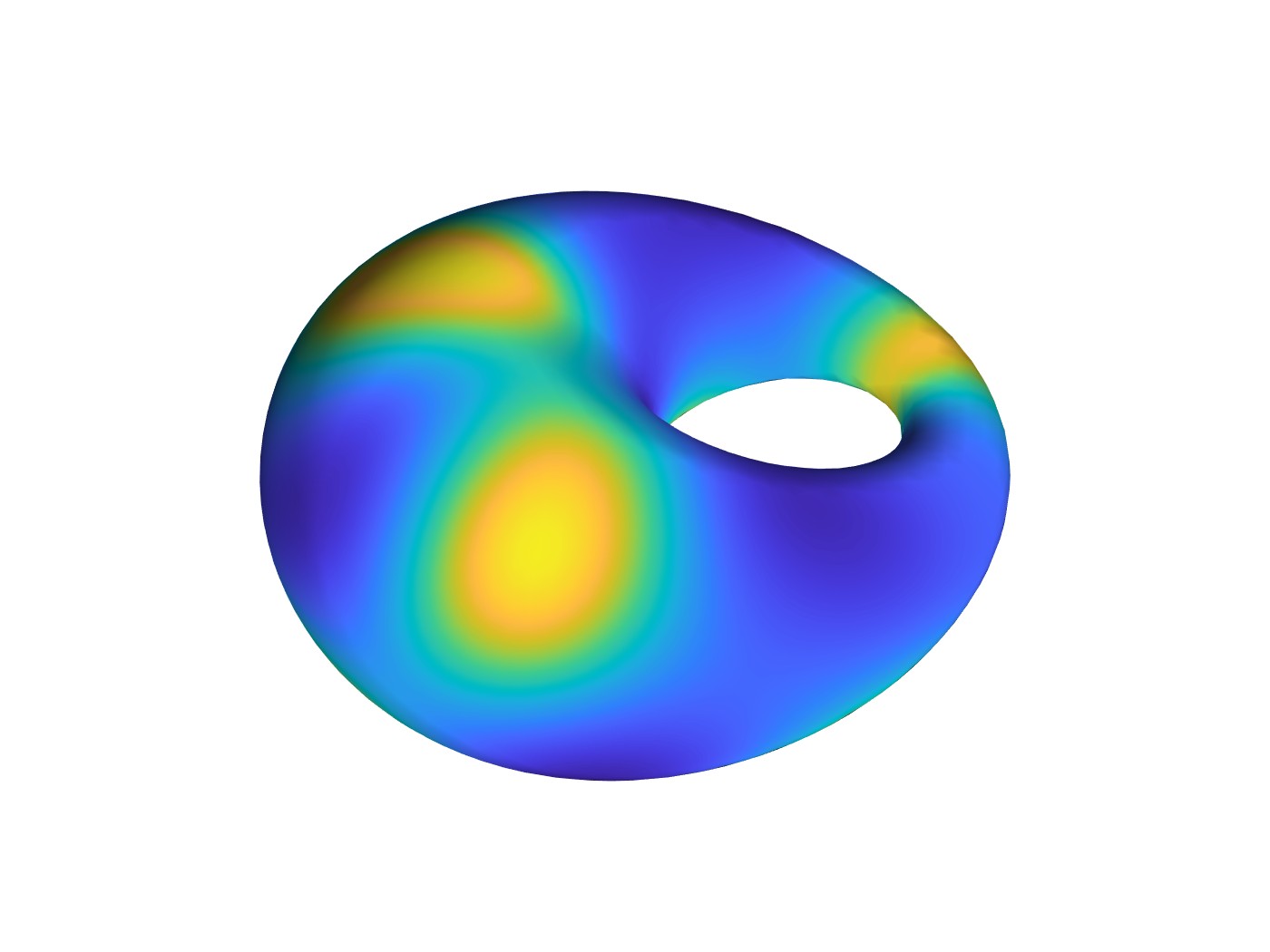}
\put(70,0){$\yellow{{\raisebox{-.5ex}{\scalebox{2.5}{$\bullet$}}} }\!\times\!6$}
    \put(85,75){$\dt=0.01$}
    \put(45,67){ $u$}
\end{overpic}
\begin{overpic}[width=0.24\textwidth,trim=130 180 130 180, clip=true,tics=10]{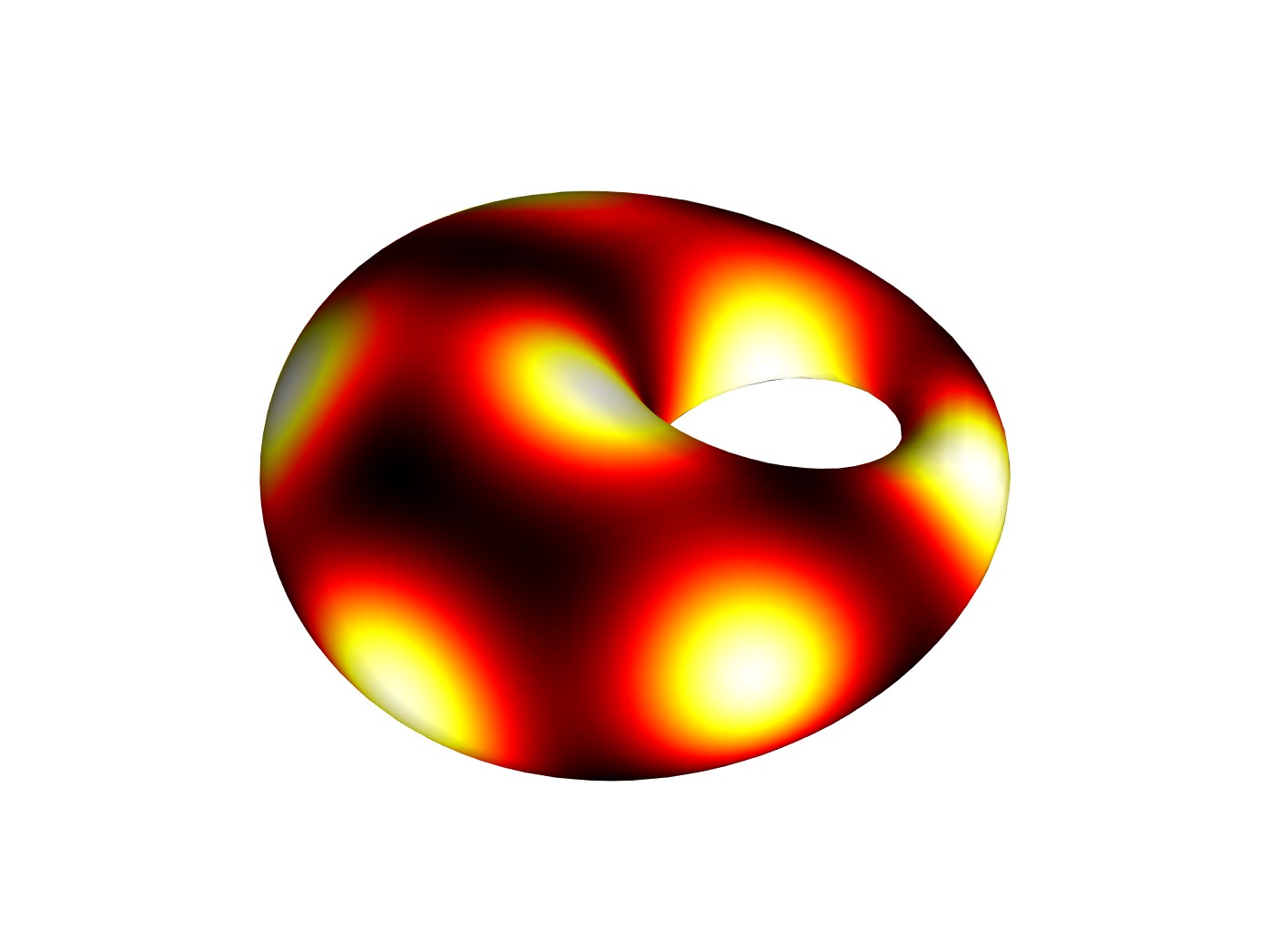}
\put(45,67){ $w$}
\put(70,0){$\yellow{{\raisebox{-.5ex}{\scalebox{2.5}{$\bullet$}}} }\!\times\!10$}
\end{overpic}
\\
  \begin{overpic}[width=0.24\textwidth,trim=130 180 130 180, clip=true,tics=10]{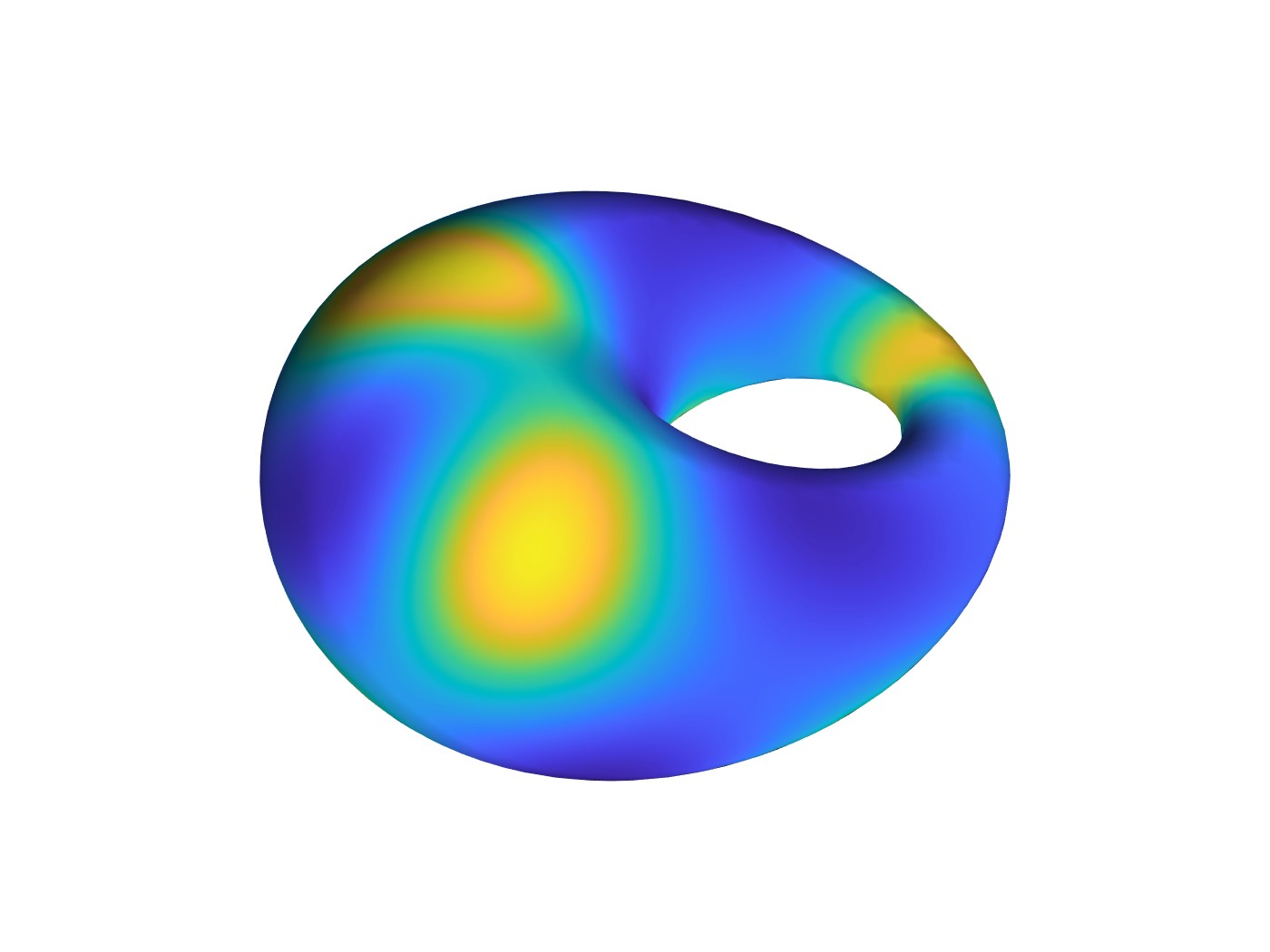}
    \put(-15,25){\rotatebox{90}{\textbf{With}}}
    \put(-5,18){\rotatebox{90}{\textbf{Greedy}}}
    \put(70,0){$\yellow{{\raisebox{-.5ex}{\scalebox{2.5}{$\bullet$}}} }\!\times\!6$}
\end{overpic}
\begin{overpic}[width=0.24\textwidth,trim=130 180 130 180, clip=true,tics=10]{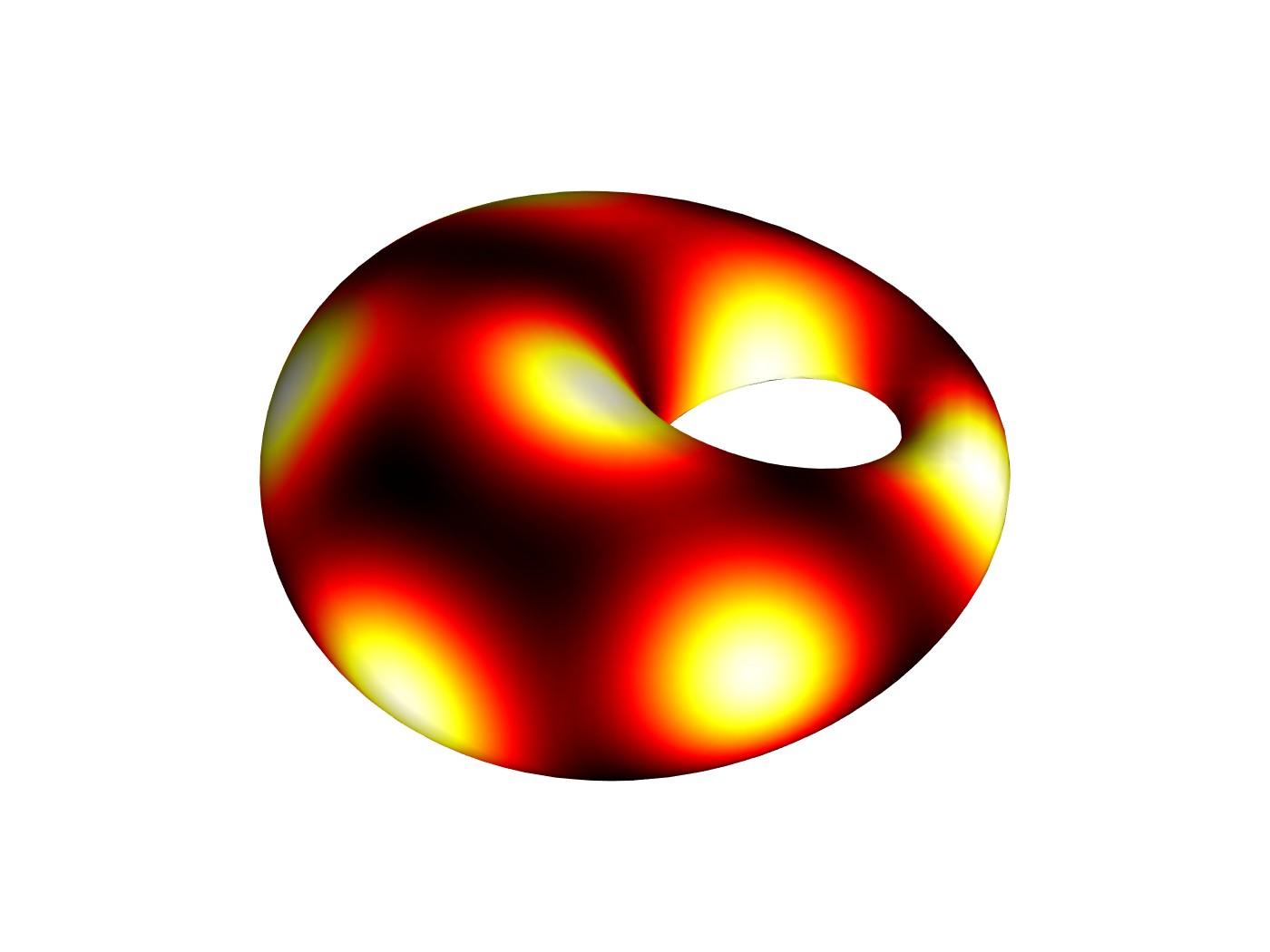}
\put(70,0){$\yellow{{\raisebox{-.5ex}{\scalebox{2.5}{$\bullet$}}} }\!\times\!10$}
\end{overpic}
\begin{overpic}[width=0.24\textwidth,trim=130 180 130 180, clip=true,tics=10]{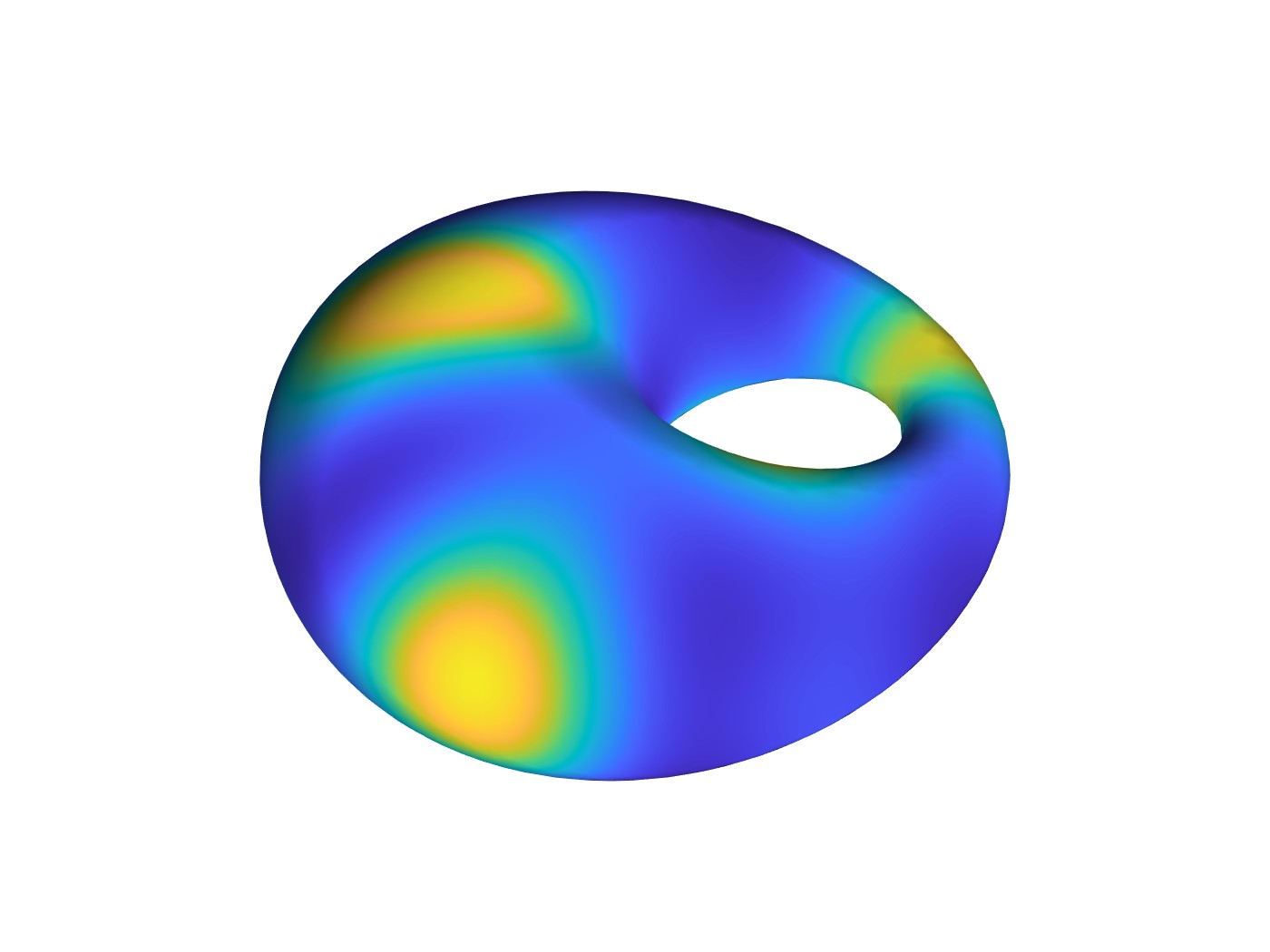}
\put(70,0){$\yellow{{\raisebox{-.5ex}{\scalebox{2.5}{$\bullet$}}} }\!\times\!6$}
\end{overpic}
\begin{overpic}[width=0.24\textwidth,trim=130 180 130 180, clip=true,tics=10]{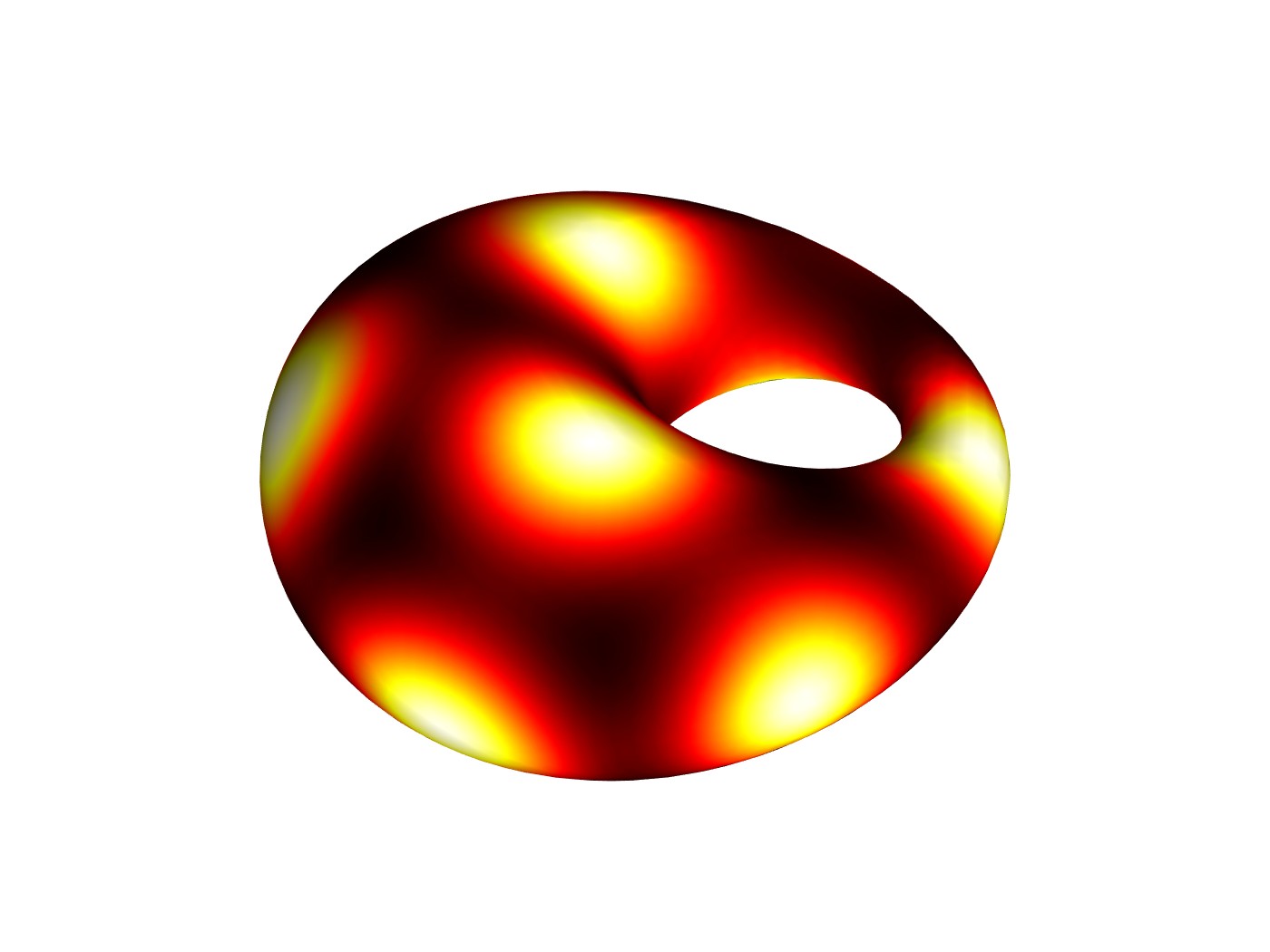}
\put(70,0){$\yellow{{\raisebox{-.5ex}{\scalebox{2.5}{$\bullet$}}} }\!\times\!10$}
\end{overpic}
\caption{Example~\ref{eg: CBSRD 3}: Solutions by solving the 3D coupled bulk-surface reaction-diffusion equation with and without the greedy algorithm in/on a Dupin's cyclide. The meshfree method uses $n = 6760+2956$ and $\epsilon = 4$. The time steps are taken as $\triangle t = 0.005$ and $\triangle t = 0.01$.
   }\label{fig: CBS3D cycl}
\end{figure}

\begin{figure}
  \centering
  \medskip
  \medskip
    \begin{overpic}[width=0.24\textwidth,trim=130 120 130 150, clip=true,tics=10]{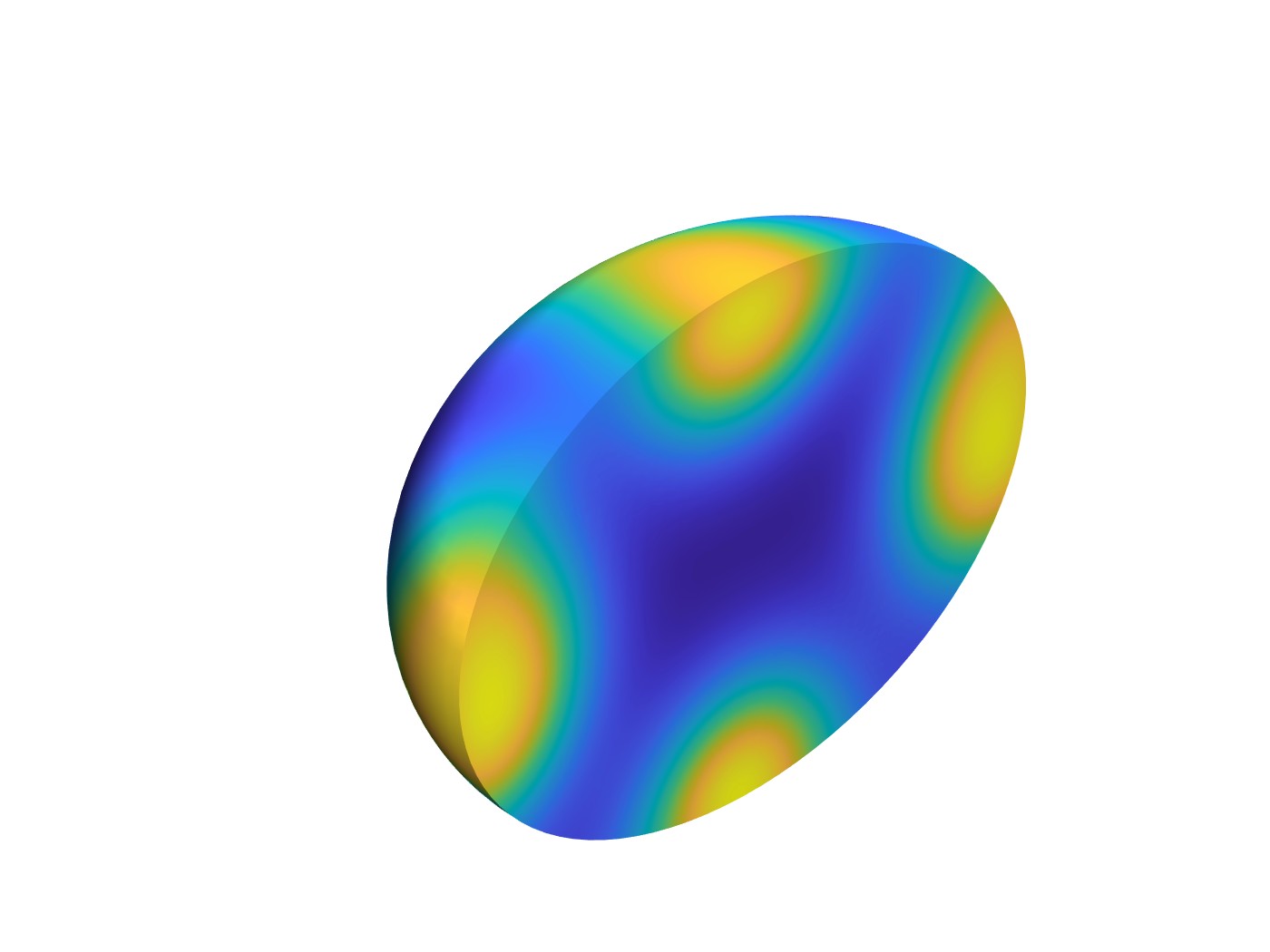}
    \put(-15,15){\rotatebox{90}{\textbf{Without}}}
    \put(-5,18){\rotatebox{90}{\textbf{Greedy}}}
    \put(80,75){$\dt=0.005$}
    \put(45,67){$u$}
    \put(73,2){$\yellow{{\raisebox{-.5ex}{\scalebox{2.5}{$\bullet$}}} }\!\times\!6$}
\end{overpic}
\begin{overpic}[width=0.24\textwidth,trim=130 150 130 150, clip=true,tics=10]{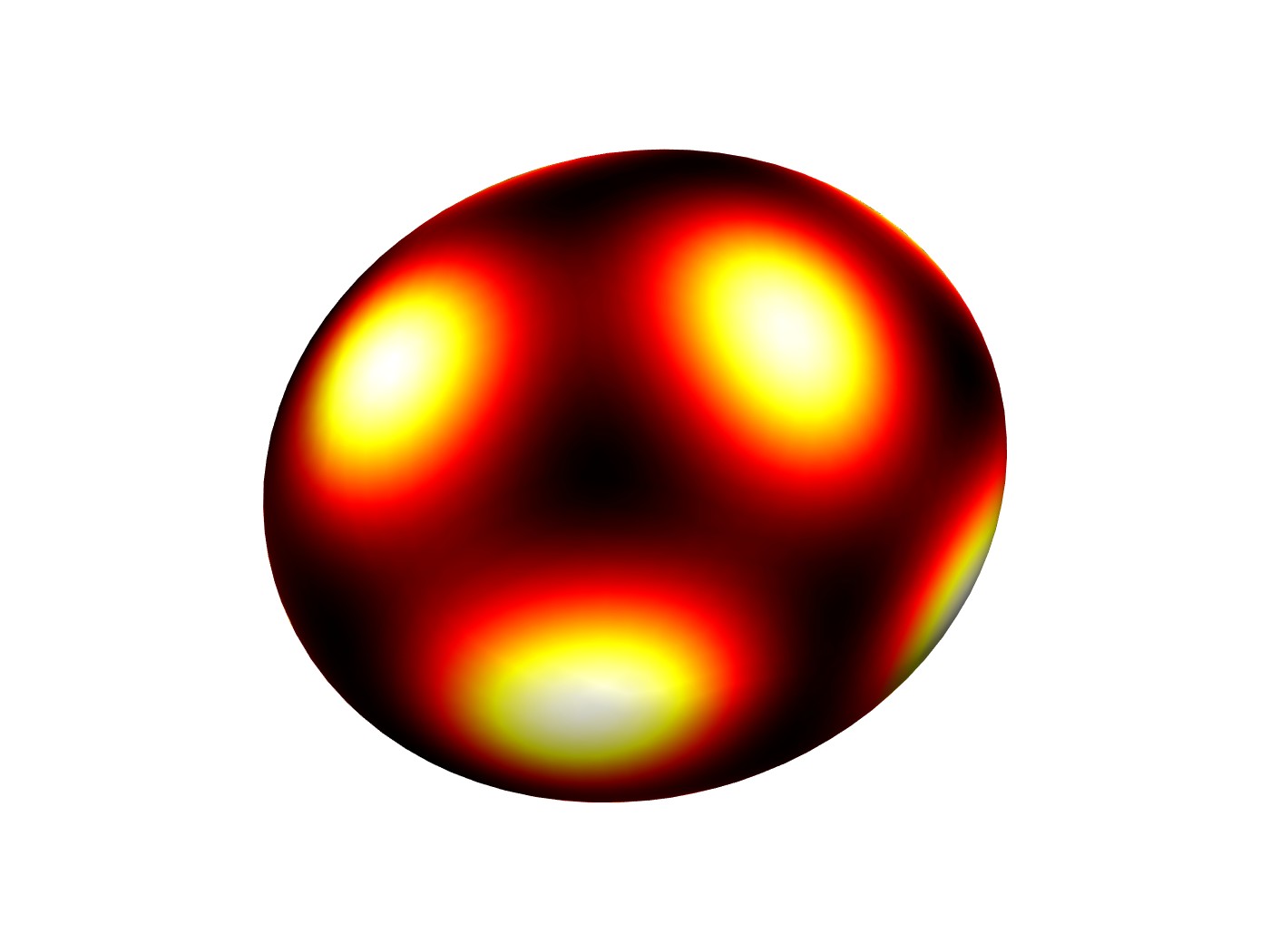}
\put(45,67){$w$}
\put(73,2){$\yellow{{\raisebox{-.5ex}{\scalebox{2.5}{$\bullet$}}} }\!\times\!9$}
\end{overpic}
  \begin{overpic}[width=0.24\textwidth,trim=130 120 130 150, clip=true,tics=10]{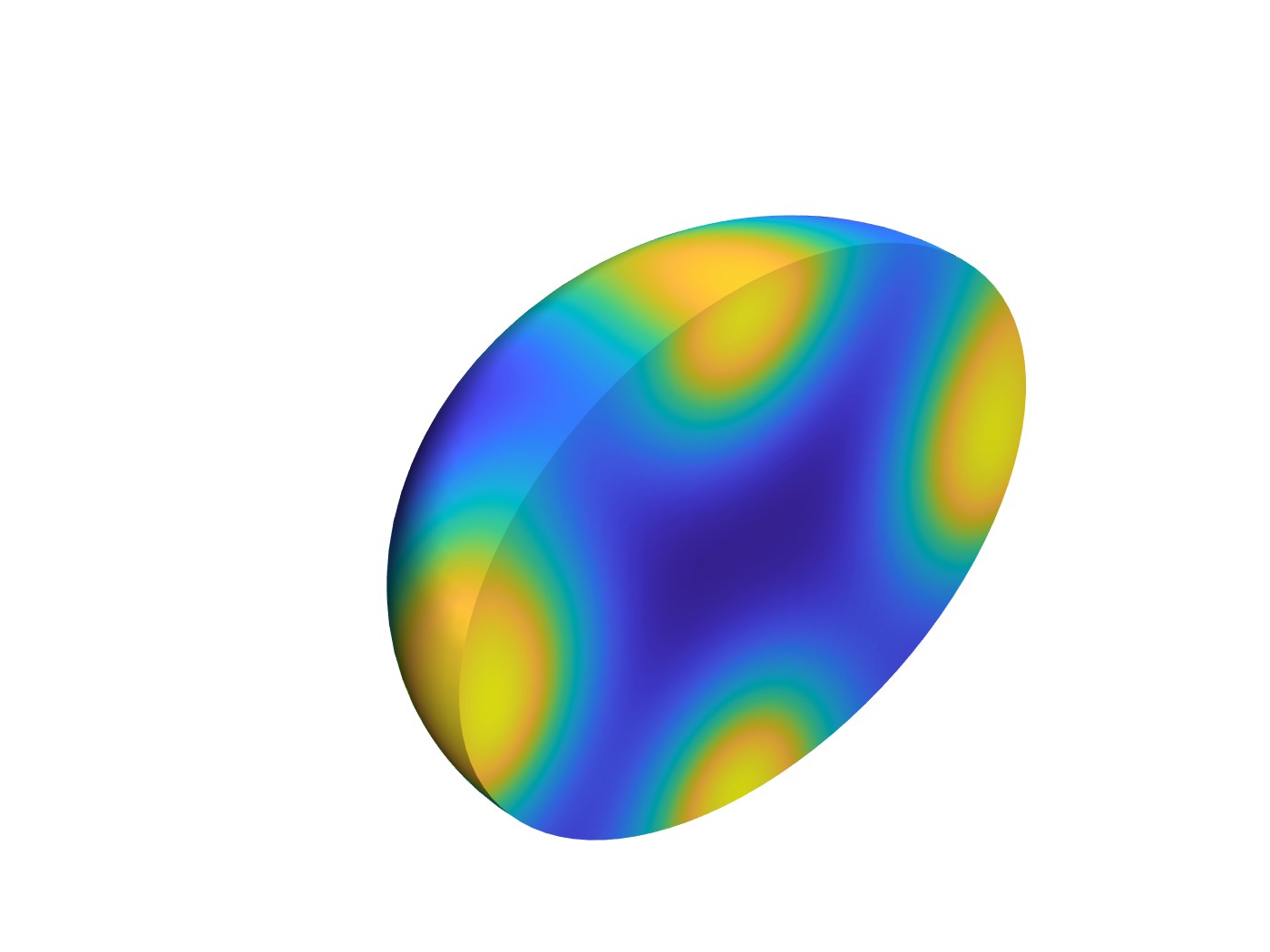}
\put(73,2){$\yellow{{\raisebox{-.5ex}{\scalebox{2.5}{$\bullet$}}} }\!\times\!6$}
    \put(80,75){$\dt=0.01$}
    \put(45,67){$u$}
\end{overpic}
\begin{overpic}[width=0.24\textwidth,trim=130 150 130 150, clip=true,tics=10]{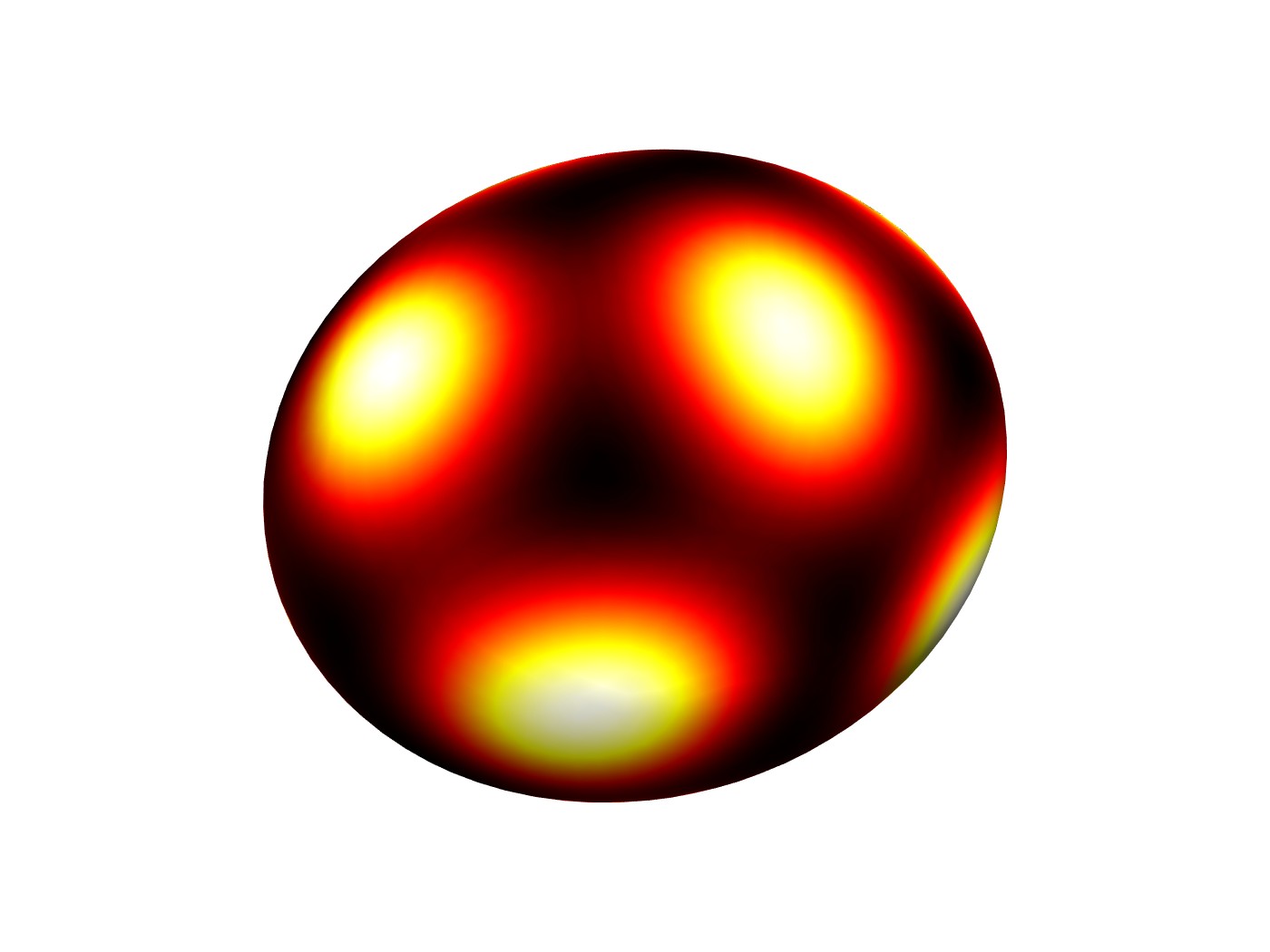}
\put(45,67){$w$}
\put(73,2){$\yellow{{\raisebox{-.5ex}{\scalebox{2.5}{$\bullet$}}} }\!\times\!9$}
\end{overpic}
\\
  \begin{overpic}[width=0.24\textwidth,trim=130 120 130 150, clip=true,tics=10]{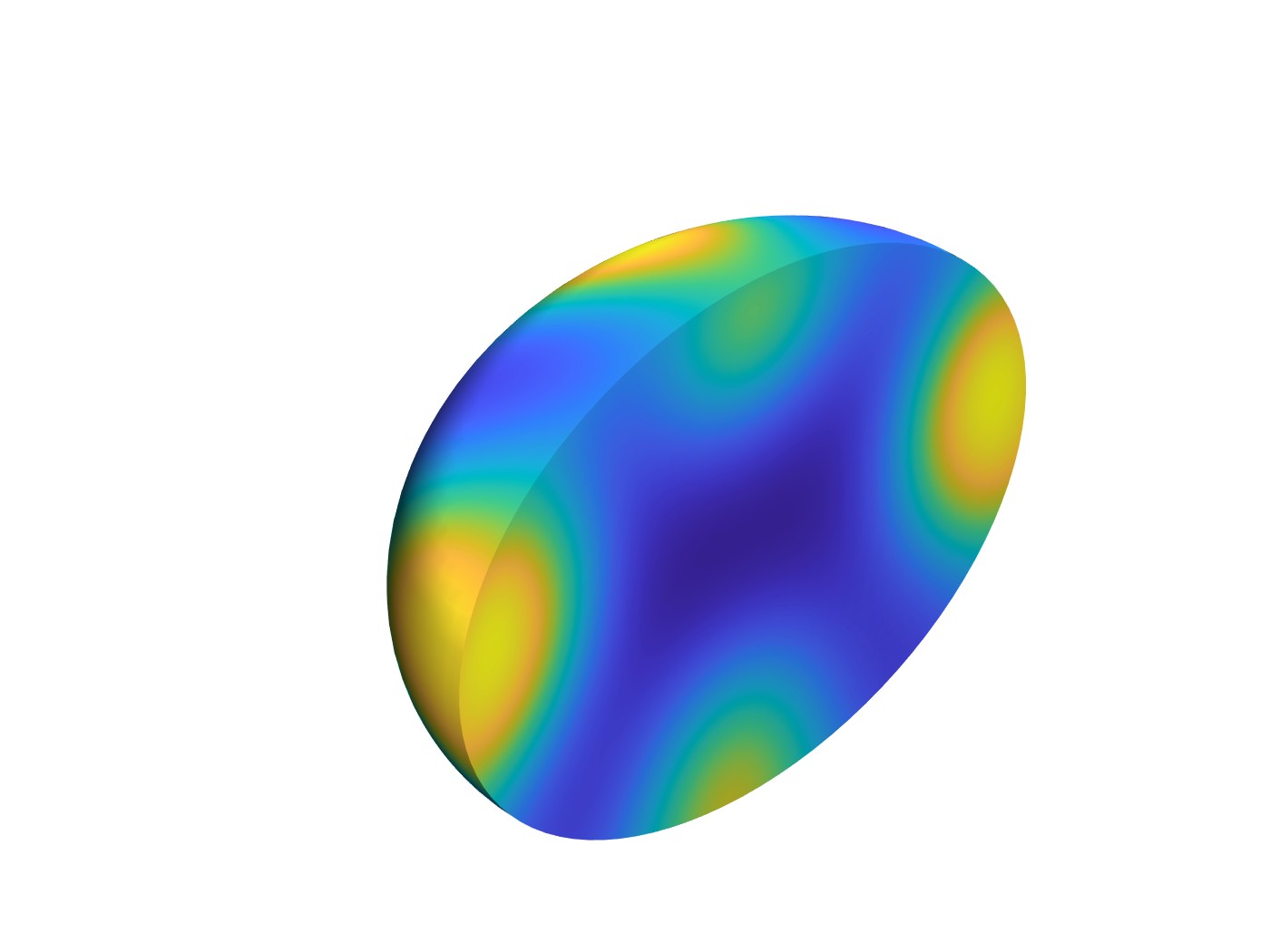}
    \put(-15,25){\rotatebox{90}{\textbf{With}}}
    \put(-5,18){\rotatebox{90}{\textbf{Greedy}}}
    \put(73,2){$\yellow{{\raisebox{-.5ex}{\scalebox{2.5}{$\bullet$}}} }\!\times\!6$}
\end{overpic}
\begin{overpic}[width=0.24\textwidth,trim=130 150 130 150, clip=true,tics=10]{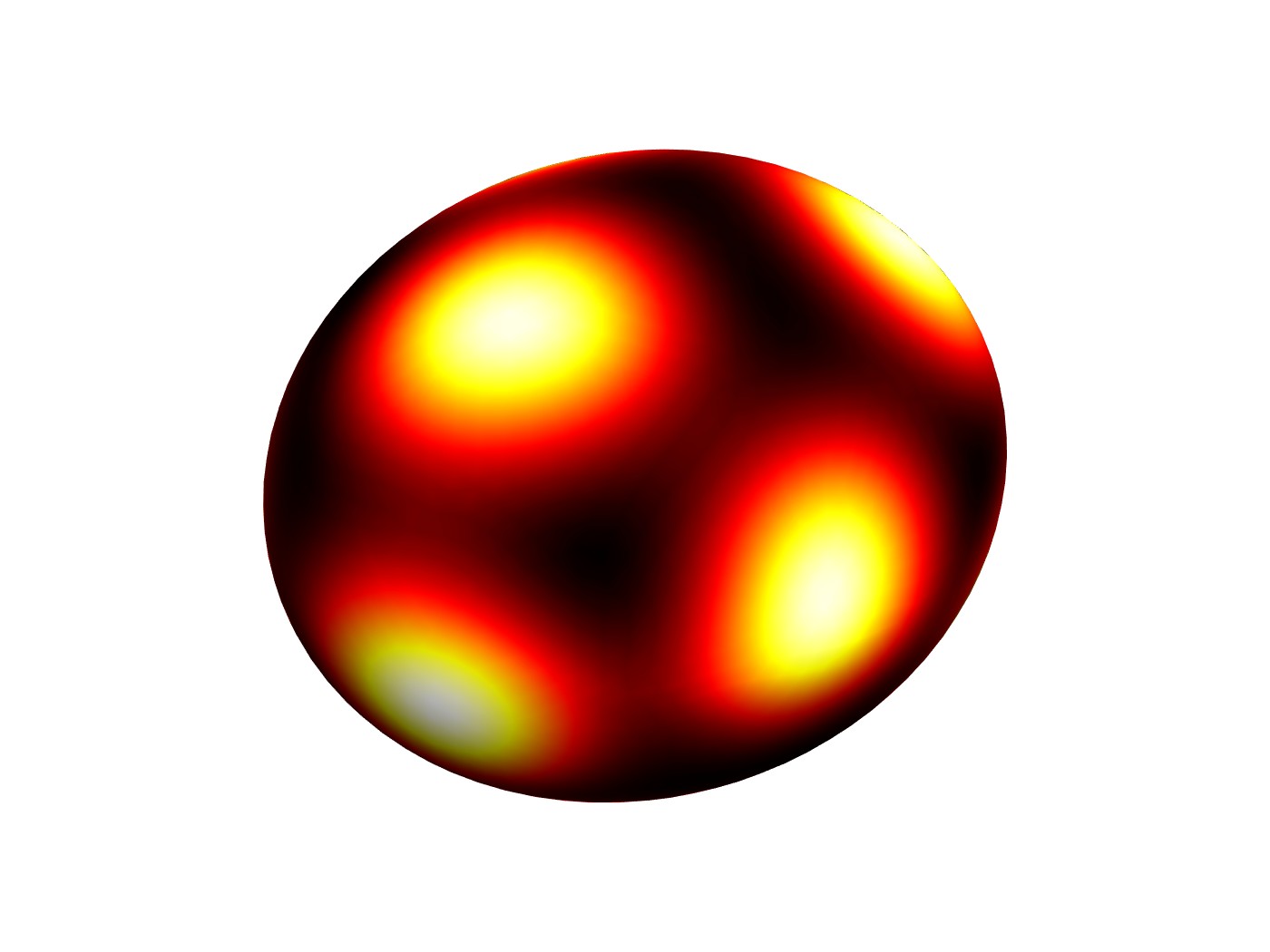}
\put(73,2){$\yellow{{\raisebox{-.5ex}{\scalebox{2.5}{$\bullet$}}} }\!\times\!9$}
\end{overpic}
\begin{overpic}[width=0.24\textwidth,trim=130 120 130 150, clip=true,tics=10]{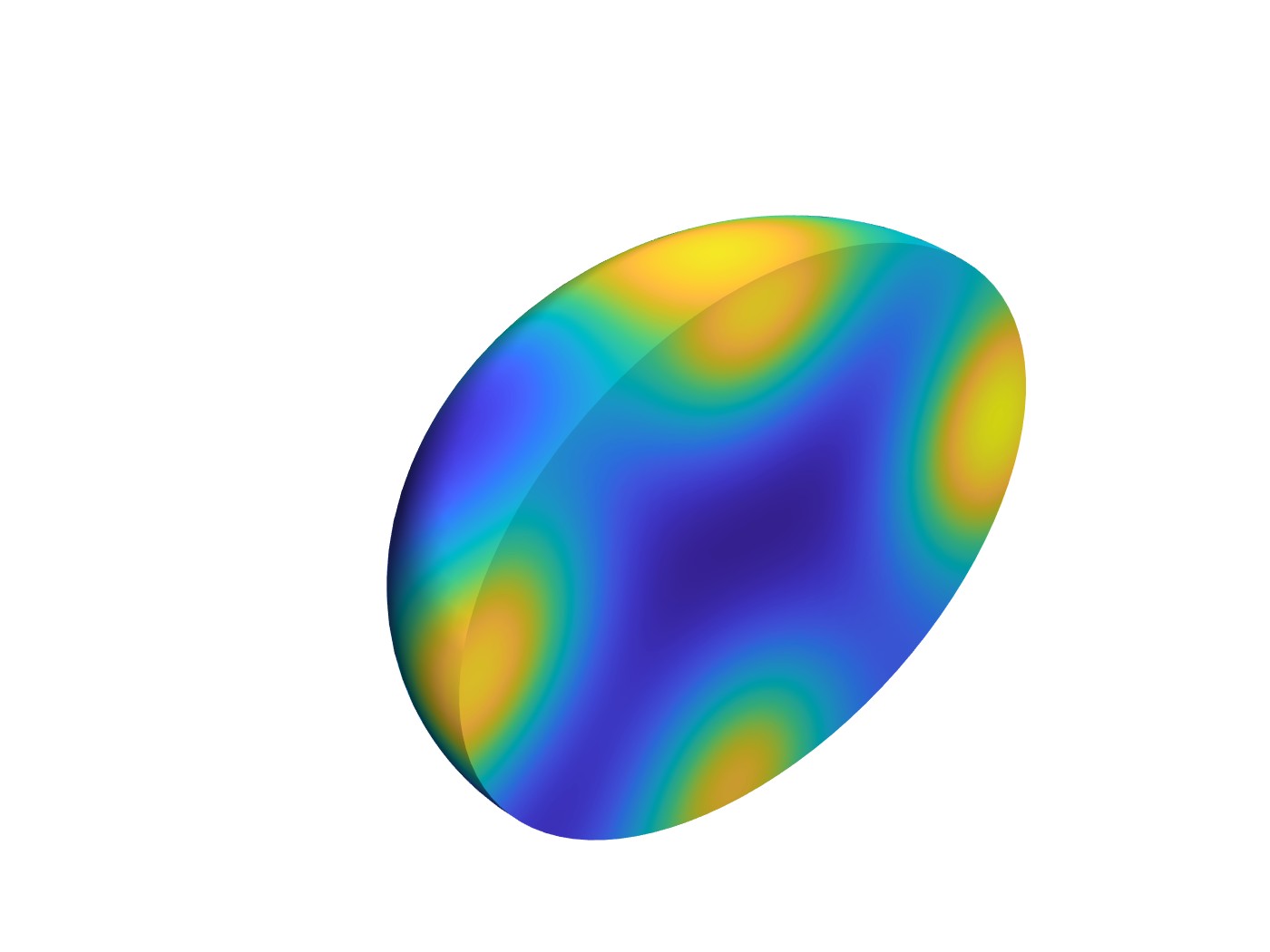}
\put(73,2){$\yellow{{\raisebox{-.5ex}{\scalebox{2.5}{$\bullet$}}} }\!\times\!6$}
\end{overpic}
\begin{overpic}[width=0.24\textwidth,trim=130 150 130 150, clip=true,tics=10]{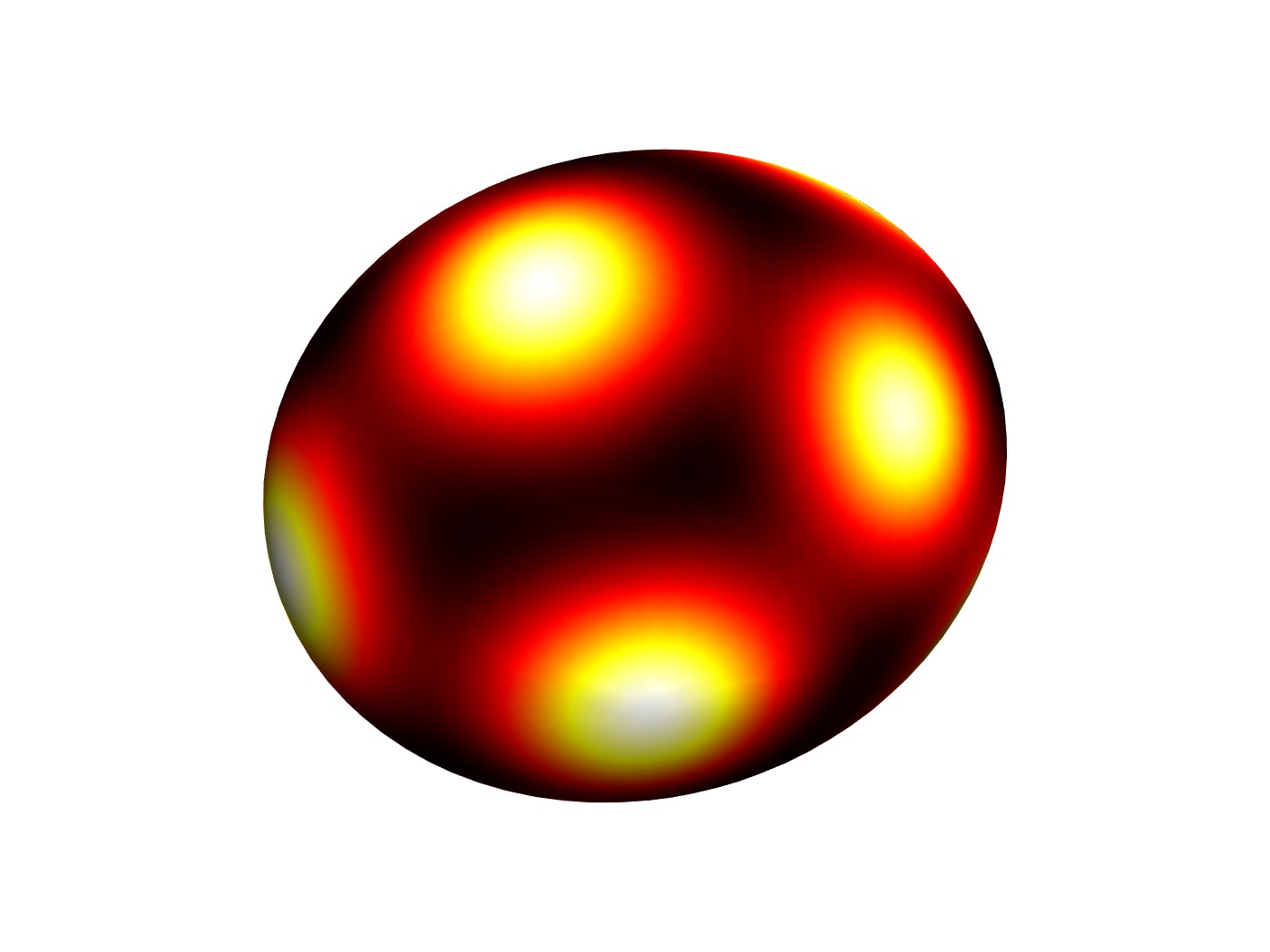}
\put(73,2){$\yellow{{\raisebox{-.5ex}{\scalebox{2.5}{$\bullet$}}} }\!\times\!9$}
\end{overpic}
\caption{Example~\ref{eg: CBSRD 3}: Solutions by solving the 3D coupled bulk-surface reaction-diffusion equation with and without the greedy algorithm in/on an ellipsoid. The meshfree method uses $n = 3395+1164$ and $\epsilon = 6$. The time steps are taken as $\triangle t = 0.005$ and $\triangle t = 0.01$.
   }\label{fig: CBS3D ellip}
\end{figure}

\section{Conclusion}
\label{sec: Conclusion}
We present two stopping criteria for the block greedy algorithm to improve stability and accuracy when solving fully-discrete matrix systems of a meshfree time-stepping method for solving parabolic PDEs.
When stopped by a large condition number, we propose a minimum residual search within newly added columns to maximize the approximation power. When stopped by a small residual, a backtracking process is used to reduce selected columns while keeping the residual below a threshold.
Together with suggested time step-dependent tolerance selection strategies, these stopping criteria enable the greedy algorithm to balance accuracy and stability. This is demonstrated in solving 2D and 3D heat equations, where the solutions' accuracy is somewhat independent of which stopping criteria is used.
We applied the greedy algorithm to coupled bulk-surface pattern formation as a practical first step. We use symmetry arguments on interfacial patterns to identify numerical defects. 

The numerical examples presented in both 2D and 3D demonstrate that, under appropriate temporal and spatial conditions, utilizing the greedy algorithm along with varying numbers of basis functions can lead to the emergence of reasonable patterns. We have made improvements to the greedy algorithm to enhance its stability when solving PDEs using meshfree and time-stepping methods. The patterns efficiently simulated using the greedy algorithm successfully maintain the general characteristics of the original configurations. In our future work, we will explore the possibility of combining our approach with variably scaled RBF \cite{Chen+SuETAL-Solvinteprobsurf:22, Marchetti+Leevan-stocexteRippalgo:22} to enhance the stability of the system even further. We may also consider combining our approach with \cite{ Petras_2019, Petras_2016} in solving the coupled bulk-surface reaction-diffusion equations on domains with moving surfaces.

\section*{Acknowledgements}
This work was supported by the General Research Fund (GRF No.
12303818, 12301419, 12301520) of Hong Kong Research Grant Council.



\begin{thebibliography}{10}
\expandafter\ifx\csname url\endcsname\relax
  \def\url#1{\texttt{#1}}\fi
\expandafter\ifx\csname urlprefix\endcsname\relax\def\urlprefix{URL }\fi
\expandafter\ifx\csname href\endcsname\relax
  \def\href#1#2{#2} \def\path#1{#1}\fi

\bibitem{AMIRFAKHRIAN2016278}
M.~Amirfakhrian, M.~Arghand, E.~Kansa, A new approximate method for an inverse time-dependent heat source problem using fundamental solutions and rbfs, Eng. Anal. Bound. Elem. 64 (2016) 278--289.

\bibitem{Vito2004SomePO}
E.~de~Vito, L.~Rosasco, A.~Caponnetto, M.~Piana, A.~Verri, Some properties of regularized kernel methods, J. Mach. Learn. Res. 5 (2004) 1363--1390.

\bibitem{Schaback+Wendland-Adapgreetechappr:00}
R.~Schaback, H.~Wendland, Adaptive greedy techniques for approximate solution of large {RBF} systems, Numer. Algorithms 24~(3) (2000) 239--254.

\bibitem{Campagna+DeMarchi-AdvComputMath:22}
R.~Campagna, S.~De~Marchi, E.~Perracchione, G.~Santin, Stable interpolation with exponential-polynomial splines and node selection via greedy algorithms, Adv. Comput. Math. 48~(6) (2022).

\bibitem{Wenzel-AnalTargDataGree:23}
T.~Wenzel, G.~Santin, B.~Haasdonk, Analysis of target data-dependent greedy kernel algorithms: Convergence rates for $f$-, $f\cdot p$- and $f/p$-greedy, Constr. Approx. 57~(1) (2023) 45--74.

\bibitem{Hon+SchabackETAL-adapgreealgosolv:03}
Y.~C. Hon, R.~Schaback, X.~Zhou, An adaptive greedy algorithm for solving large {RBF} collocation problems, Numer. Algorithms 32~(1) (2003) 13--25.

\bibitem{Ling+OpferETAL-Resumeshcolltech:06}
L.~Ling, R.~Opfer, R.~Schaback, Results on meshless collocation techniques, Eng Anal Bound Elem 30~(4) (2006) 247--253.

\bibitem{Ling+Schaback-imprsubsselealgo:09}
L.~Ling, R.~Schaback, An improved subspace selection algorithm for meshless collocation methods, Int J Numer Methods Eng 80~(13) (2009) 1623--1639.

\bibitem{Ling-fastblocalgoquas:16}
L.~Ling, A fast block-greedy algorithm for quasi-optimal meshless trial subspace selection, SIAM J. Sci. Comput. 38~(2) (2016) A1224--A1250.

\bibitem{Ling+Chiu-Fulladapkernmeth:18}
L.~Ling, S.~N. Chui, Fully adaptive kernel-based methods, Int J Numer Methods Eng 114~(4) (2018) 454--467.

\bibitem{chen2023exploring}
M.~Chen, L.~Ling, Exploring oversampling in rbf least-squares collocation method of lines for surface diffusion (2023).
\newblock \href {http://arxiv.org/abs/2203.08579} {\path{arXiv:2203.08579}}.

\bibitem{Chen_2023}
M.~Chen, K.~C. Cheung, L.~Ling, A kernel-based least-squares collocation method for surface diffusion, SIAM J. Numer. Anal. 61~(3) (2023) 1386–1404.

\bibitem{Hammarling+Lucas-UpdaFactLeasSqua:08}
S.~Hammarling, C.~Lucas, Updating the {QR} factorization and the least squares problem, MIMS EPrint 2008.111, Manchester Institute for Mathematical Sciences, School of Mathematics, University of Manchester, Manchester, UK (2008).

\bibitem{Ruuth-Implmethreacprob:95}
S.~J. Ruuth, Implicit-explicit methods for reaction-diffusion problems in pattern formation, J. Math. Biol. 34~(2) (1995) 148--176.

\bibitem{Chen+Ling-Extrmeshcollmeth:20}
M.~Chen, L.~Ling, Extrinsic meshless collocation methods for {PDE}s on manifolds, SIAM J. Numer. Anal. 58~(2) (2020) 988--1007.

\bibitem{Cheung+Ling-Kernembemethconv:18}
K.~C. Cheung, L.~Ling, A kernel-based embedding method and convergence analysis for surfaces {PDE}s, SIAM J. Sci. Comput. 40~(1) (2018) A266--A287.

\bibitem{Cheung+LingETAL-leaskerncollmeth:18}
K.~C. Cheung, L.~Ling, R.~Schaback, ${H}^2$--convergence of least-squares kernel collocation methods, SIAM J. Numer. Anal. 56~(1) (2018) 614--633.

\bibitem{Santin_2021}
G.~Santin, T.~Karvonen, B.~Haasdonk, Sampling based approximation of linear functionals in reproducing kernel {H}ilbert spaces, BIT Numer. Math. 62~(1) (2021) 279--310.

\bibitem{Symmetrybreaking2014}
A.~Rätz, M.~Röger, Symmetry breaking in a bulk-surface reaction-diffusion model for signaling networks, Nonlinearity 27 (2014).

\bibitem{Rtz2011TuringII}
A.~R{\"a}tz, M.~R{\"o}ger, Turing instabilities in a mathematical model for signaling networks, J. Math. Biol. 65 (2011) 1215 -- 1244.

\bibitem{2022IJMMM..29...32L}
Q.~{Li}, X.~{Lin}, Q.~{Luo}, Y.~{Chen}, J.~{Wang}, B.~{Jiang}, F.~{Pan}, {Kinetics of the hydrogen absorption and desorption processes of hydrogen storage alloys: A review}, Int. J. Miner. Metall. Mater. 29~(1) (2022) 32--48.

\bibitem{BSMz}
A.~Madzvamuse, H.~Chung, C.~Venkataraman, Stability analysis and simulations of coupled bulk-surface reaction-diffusion systems, Proc R Soc (London) A. 471 (2015).

\bibitem{Turing1}
J.~Murray, {Mathematical Biology {II}: Spatial Models and Biomedical Applications}, {Springer}, 1997.

\bibitem{Fuselier+Wright-ScatDataInteEmbe:12}
E.~J. Fuselier, G.~B. Wright, Scattered data interpolation on embedded submanifolds with restricted positive definite kernels: {S}obolev error estimates, SIAM J. Numer. Anal. 50~(3) (2012) 1753--1776.

\bibitem{CMBS}
M.~Chen, L.~Ling, Kernel-based meshless collocation methods for solving coupled bulk–surface partial differential equations, J. Sci. Comput. 81 (2019) 375–391.

\bibitem{CBSCM}
C.~B. Macdonald, B.~Merriman, S.~J. Ruuth, Simple computation of reaction–diffusion processes on point clouds, Proceedings of the National Academy of Sciences 110 (2013) 9209 -- 9214.

\bibitem{Chen+SuETAL-Solvinteprobsurf:22}
M.~Chen, L.~Ling, Y.~Su, Solving interpolation problems on surfaces stochastically and greedily, Dolomites Res. Notes Approx. 15~(3) (2022) 26--36.

\bibitem{Marchetti+Leevan-stocexteRippalgo:22}
F.~Marchetti, L.~Ling, A stochastic extended {R}ippa's algorithm for {LpOCV}, Appl. Math. Lett. 129 (2022).

\bibitem{Petras_2019}
A.~Petras, L.~Ling, C.~Piret, S.~Ruuth, A least-squares implicit {RBF}-{FD} closest point method and applications to {PDEs} on moving surfaces, J. Comput. Phys. 381 (2019) 146--161.

\bibitem{Petras_2016}
A.~Petras, S.~Ruuth, {PDEs} on moving surfaces via the closest point method and a modified grid based particle method, J. Comput. Phys. 312 (2016) 139--156.

\end{thebibliography}

\end{document}